\documentclass[12pt]{amsart}
\usepackage{amsmath,amssymb}
\newtheorem{theorem}{Theorem}
\newtheorem{lemma}{Lemma}
\newtheorem{proposition}{Proposition}
\newtheorem{corollary}{Corollary}
\begin{document}
\centerline{\Large Yang-Mills fields on CR manifolds} \vskip 0.3in

\centerline{\large Elisabetta Barletta \hspace{1cm} Sorin
Dragomir} {\small
\begin{center}
\par\noindent Universit\`a degli Studi della Basilicata,
\par\noindent
Dipartimento di Matematica, \par\noindent Campus Macchia Romana,
\par\noindent 85100 Potenza, Italy, \par\noindent
e-mail: {\tt barletta@unibas.it} $\;\;\;$ {\tt dragomir@unibas.it}
\end{center}} \vskip 0.1in \centerline{\large Hajime Urakawa}
{\small
\begin{center}
\par\noindent Division of Mathematics \par\noindent Graduate
School of Information Sciences \par\noindent Tohoku University
\par\noindent Aoba 09, Sendai, 980-8579, Japan \par\noindent
e-mail: {\tt urakawa@math.is.tohoku.ac.jp} \end{center}}
\title{}
\author{}
\begin{abstract}
We study {\em pseudo Yang-Mills fields} on a compact strict\-ly
pseudoconvex CR manifold $M$, i.e. the critical points of the
functional $\mathcal{PYM}(D) = \frac{1}{2} \int_M \| \pi_H R^D
\|^2 \theta \wedge (d \theta )^n$, where $D$ is a connection in a
Hermitian CR-ho\-lo\-mor\-phic vector bundle $(E,h) \to M$. Let
$\Omega = \{ \varphi < 0 \} \subset \mathbb{C}^n$ be a smoothly
bounded strictly pseuodoconvex domain and $g$ the Bergman metric
on $\Omega$. We show that boundary values $D_b$ of Yang-Mills
fields $D$ on $(\Omega , g)$ are pseudo Yang-Mills fields on
$\partial \Omega$, provided that $i_T R^{D_b} = 0$ and $i_N R^D =
0$ on $H(\partial \Omega )$. If $S^1 \to C(M)
\stackrel{\pi}{\rightarrow} M$ is the canonical circle bundle and
$\pi^* D$ is a Yang-Mills field with respect to the Fefferman
metric $F_\theta$ of $(M , \theta )$ then $D$ is a pseudo
Yang-Mills field on $M$. The Yang-Mills equations $\delta^{\pi^*
D} R^{\pi^* D} = 0$ project on the Euler-Lagrange equations
$\delta^D_b R^D = 0$ of the variational principle $\delta \;
\mathcal{PYM} (D) = 0$, provided that $i_T R^D = 0$. When $M$ has
vanishing pseudohermitian Ricci curvature the pullback $\pi^* D$
of the (CR invariant) Tanaka connection $D$ of $(E , h)$ is a
Yang-Mills field on $C(M)$. We derive the second variation formula
$\{ d^2 \, \mathcal{PYM}(D^t )/d t^2 \}_{t=0} = \int_M \langle
\mathcal{S}_b^D (\varphi ), \varphi \rangle \, \theta \wedge (d
\theta )^n$, $D^t = D + A^t$ (provided that $D$ is a pseudo
Yang-Mills field and $\varphi \equiv \{ d A^t /d t\}_{t=0} \in
{\rm Ker}(\delta^D )$), and show that $\mathcal{S}^D_b (\varphi )
\equiv \Delta^D_b \varphi + \mathcal{R}^D_b (\varphi )$, $\varphi
\in \Omega^{0,1} ({\rm Ad}(E))$, is a subelliptic operator.
\end{abstract}
\maketitle

\section{Introduction}
A series of papers published in the last decade (J. Lewandowski \&
P. Nurowski, \cite{kn:LeNu}, P. Nurowski, \cite{kn:Nur}, P.
Nurowski \& J. Tafel, \cite{kn:NuTa}) are devoted to exploring the
relationship among CR structures on $3$-dimensional ma\-ni\-folds
and null solutions to Ein\-stein equations, Max\-well equations,
and Yang-Mills equations (cf. also J. Tafel, \cite{kn:Taf}).
Specifically, if $(M, T_{1,0}(M))$ is a nondegenerate
$3$-dimensional CR manifold endowed with the contact form $\theta$
and with the (locally defined) complex $1$-form $\theta^1$ such
that $\theta^1 (T_1 ) = 1$, $\theta^1 (T_{\overline{1}}) = 0$
(where $T_1$ is a local generator of the CR structure
$T_{1,0}(M)$) let us consider the semi-Riemannian metric
\begin{equation}
F = 2 p^2 \{ (\pi^* \theta^1 ) \odot (\pi^* \theta^{\overline{1}})
- (\pi^* \theta ) \odot \sigma \} , \label{e:i1}
\end{equation}
on $M \times \mathbb{R}$, where $p$ is a real valued function on
$M \times \mathbb{R}$ and $\sigma$ is a real $1$-form on $M \times
\mathbb{R}$ such that
\[ \pi^* \left( \theta  \wedge \theta^1 \wedge
\theta^{\overline{1}} \right) \wedge \sigma \neq 0. \] Here $\pi :
M \times \mathbb{R} \to M$ is the projection. P. Nurowski has
determined (cf. \cite{kn:Nur}) local solutions to the Yang-Mills
equations on $(M \times \mathbb{R}, F)$, under the additional
assumption that the shear-free congruence of null geodesics
tangent to $\partial /\partial \gamma$ ($\gamma$ is the natural
coordinate function on $\mathbb{R}$) possesses $3$ linearly
independent symmetries $\{ X_i : 1 \leq i \leq 3 \}$. Let $D$ be a
$SU(3)$-connection in a vector bundle $\hat{E} \to M \times
\mathbb{R}$ locally described by a matrix of $1$-forms $A = b \;
\pi^* \theta^1 + \overline{b} \; \pi^* \theta^{\overline{1}} + c
\; \pi^* \theta + e \; \sigma$, where $a$, $b$, $c$ and $e$ are
$\mathcal{G} \otimes \mathbb{C}$-valued functions ($\mathcal{G} =
{\bf su}(3)$) on $M \times \mathbb{R}$. When the Lie group $G_3$
(whose Lie algebra is generated by the $X_i$'s) consists of
symmetries of $D$ (i.e. each element of $G_3$ induces a gauge
transformation of $A$) then (by a result of J. Harnad \& S.
Shnider \& L. Vinet, \cite{kn:HaShVi}) up to some gauge
transformation $A$ is strictly invariant under $G_3$. Then
$\mathcal{L}_{X_i} A = 0$, conditions which may be exploited to
show that locally $D$ may be looked for in the form $A = B \,
\pi^* \Omega^1 + \overline{B} \, \pi^* \Omega^{\overline{1}} + C
\, \pi^* \Omega$, with $B,C \in \mathcal{G} \otimes \mathbb{C}$.
Here $\Omega$ and $\Omega^1$ are a new contact form and a new
local coframe such that
\begin{equation}
\mathcal{L}_{\tilde{X}_i} \Omega = 0, \;\;
\mathcal{L}_{\tilde{X}_i} \Omega^1 = 0, \;\; d \Omega = 2
\sqrt{-1} \, \Omega^1 \wedge \Omega^{\overline{1}} ,
\label{e:i2}
\end{equation}
while $\tilde{X}_i$ are the projections on $M$ of the (nontrivial)
symmetries $X_i$. Finally the Yang-Mills equations (for
$SU(3)$-fields) on $M \times \mathbb{R}$ may be solved (together
with the condition that $D$ is null, i.e. $d A + A \wedge A =
(\pi^* \Omega ) \wedge (\Phi \; \pi^* \Omega^1 + \overline{\Phi}
\; \pi^* \Omega^{\overline{1}})$, for some $\mathcal{G} \otimes
\mathbb{C}$-valued function $\Phi$). For instance (cf. (5.10) in
\cite{kn:Nur}, p. 805)
\begin{equation}
A = \rho \; \vec{n} \cdot f \; \pi^* (e^{i \phi} \Omega^1 + e^{-i
\phi} \Omega^{\overline{1}}) \label{e:i3}
\end{equation}
is a solution, where $\vec{n} \in \mathbb{R}^3$ is a unit vector,
$\rho \in \mathbb{R}$, $\phi \in [0, 2 \pi ]$, and $f = (e_1 , e_2
, e_3 )$ is a basis in $\mathcal{G}$. It is noteworthy that
$\tilde{X}_i$ turn out to be symmetries of the CR structure
$T_{1,0}(M)$ (in a sense that will be explained in section 3) and
that the CR structures admitting the $3$-dimensional symmetry
group $G_3$ are fully classified in \cite{kn:LeNu}, according to
the Bianchi type of $G_3$. For instance, if $M = \{ (x,y,z) \in
\mathbb{R}^3 : y \neq 0 \}$ carries the CR structure $T_{1,0}(M) =
\mathbb{C} T_1$ with
\begin{equation}
T_1 = \frac{y}{1 + y^2} \; \frac{\partial}{\partial x} - \frac{i
\, y}{2} \; \frac{\partial}{\partial y} + \frac{1}{y(1+y^2 )} \;
\frac{\partial}{\partial z}
\label{e:gen}
\end{equation}
and the contact form $\theta = (1/y) \, d x - y \, d z$ (such $M$
possesses a symmetry group of Bianchi type $VI_0$) then a local
solution $A$ (to the Yang-Mills equations) of the form
(\ref{e:i3}) may be produced. The example (\ref{e:gen}) of a CR
structure on $\mathbb{R}^3 \setminus \{ y = 0 \}$ will be
encountered again in section 3. Let now $M$ be a compact strictly
pseudoconvex CR manifold, of arbitrary CR dimension $n$. Let $S^1
\to C(M) \stackrel{\pi}{\longrightarrow} M$ be the canonical
circle bundle (cf. section 3 for definitions). Note that $C(M)$
and $M \times \mathbb{R}$ are locally diffeomorphic. If $\theta$
is a contact form on $M$ then $C(M)$ carries a natural Lorentz
metric $F_\theta$ (the Fefferman metric) and a moment's thought
(compare to (\ref{e:4}) in section 5) shows that when $M$ is
$3$-dimensional the Fefferman metric $F_\theta$ is of the form
(\ref{e:i1}). Then, under the symmetry assumptions above,
(\ref{e:i3}) is a (local) solution to the Yang-Mills equations
\begin{equation}
\delta^{\mathbb{D}} R^{\mathbb{D}} = 0
\label{e:i4}
\end{equation}
on $(C(M), F_\theta )$ (with $n=1$), and in general it is
conceivable that when the CR structure $T_{1,0}(M)$ possesses a
symmetry group $G_{2n+1}$, Nurowski's scheme may produce local
symmetric null solutions to (\ref{e:i4}). A first step towards the
achievement of this goal is performed in section 3. Note that
(\ref{e:i3}) is the pullback (via $\pi$) to $M \times \mathbb{R}$
of a field on $M$. It is then a natural question whether given a
Yang-Mills field on $(C(M), F_\theta )$ of the form $\pi^* D$, it
follows that $D$ is a Yang-Mills field on $(M, g_\theta )$, where
$g_\theta$ is the Webster metric. This question is answered in
section 5, where we integrate along the fibre in the Yang-Mills
functional $\widehat{\mathcal{Y}\mathcal{M}}$ on $C(M)$ and
produce the new functional (\ref{e:funct}). As it turns out, $D$
is a {\em pseudo Yang-Mills} field (i.e. a critical point of
(\ref{e:funct})) rather than a Yang-Mills field on $(M , g_\theta
)$ (however, the two notions coincide in the special case $i_T R^D
= 0$). The converse (i.e. whether given a pseudo Yang-Mills field
$D$ on $M$ its pullback $\pi^* D$ is a Yang-Mills field on $C(M)$)
is examined in Theorem \ref{t:1}. Solving (\ref{e:i4}) on $C(M)$
is therefore closely related to solving the pseudo Yang-Mills
equations
\begin{equation} \delta^D_b R^D = 0
\label{e:i5}
\end{equation}
on $M$, and indeed (\ref{e:i4}) projects (under additional
conditions, cf. section 5) on $M$ to give (\ref{e:i5}). One of the
main results in this paper is that solutions to (\ref{e:i5}) occur
as boundary values of Yang-Mills fields on a strictly pseudoconvex
bounded domain $\Omega \subset \mathbb{C}^n$ endowed with the
Bergman metric $g$ (cf. Theorem \ref{t:3}). For the proof of
Theorem \ref{t:1} we draw inspiration from \cite{kn:GrLe} and make
use of their canonical connection $\nabla$ (the {\em Graham-Lee
connection}) whose pointwise restriction to a level set (near
$\partial \Omega$) of a defining function of $\Omega$ is the
better known Tanaka-Webster connection of the level set. Using the
fine asymptotic properties of the Bergman kernel of $\Omega$ we
may choose a defining function allowing an explicit relationship
among the Bergman metric $g$ and the Webster metric of each level
set, and therefore an explicit relationship among the Levi-Civita
connection of $(\Omega , g)$ and the Graham-Lee connection. In the
end, an elementary asymptotic analysis shows that boundary values
$D_b$ of Yang-Mills fields $D$ on $(\Omega , g)$ satisfy
(\ref{e:i5}) provided that $D_b$ satisfy certain compatibility
conditions along $\partial \Omega$ (cf. section 4). In sections 6
and 7 we obtain the first and second variation formulae for the
functional (\ref{e:funct}). The relevant operator occurring in the
second variation formula is shown to be subelliptic of order $1/2$
(cf. Theorem \ref{t:2}). The problem of building an appropriate
stability theory (along the lines of \cite{kn:BoLa}, yet relying
on the subelliptic rather than on the elliptic theory) remains
open. We feel that the importance of the Graham-Lee connection
$\nabla$ in applications deserves Appendix A: there we provide a
new axiomatic description of $\nabla$ together with a index-free
proof.

\section{Statement of main results}
Let $(M, T_{1,0}(M))$ be a compact strictly
pseudoconvex CR manifold, of CR dimension $n$, and $\theta$ a
contact form on $M$. Let $(E, \overline{\partial}_E ) \to M$ be a
CR-holomorphic vector bundle and $h$ a Hermitian metric in $E$.
Let ${\mathcal C}(E , h)$ be the affine space of all connections
$D$ in $E$ such that $D h = 0$. We consider the functional
\begin{equation}
\mathcal{PYM}(D) = \frac{1}{2} \int_M \| \pi_H R^D \|^2 \, \theta
\wedge (d \theta )^n .
\label{e:funct}
\end{equation}
Here $\pi_H : \Omega^2 ({\rm \, Ad} E) \to \Omega^2 ({\rm \, Ad}
E)/ {\mathcal J}_\theta^2$ is the natural projection and
${\mathcal J}^{\bullet}_\theta$ the ideal generated by $\theta$ in
$\Omega^{\bullet}({\rm \, Ad} E)$. A {\em pseudo Yang-Mills field}
on $M$ is a critical point of $\mathcal{PYM} : {\mathcal C}(E, h)
\to [0, + \infty )$. We shall show that
\begin{theorem} Let $\Omega = \{ z \in U : \varphi (z) < 0 \}$
be a smoothly bounded strictly pseudoconvex domain in
$\mathbb{C}^n$ and $g$ its Bergman metric. Let $\pi : F \to U$ be
a holomorphic vector bundle and $h$ a Hermitian metric on $F$. Let
$D_b \in \mathcal{C}(E, h)$ $(E = \pi^{-1}(\partial \Omega ))$ be
the boundary values of a Yang-Mills field $D \in \mathcal{C}(F ,
h)$ on $(\Omega , g)$. Assume that $i_T R^{D_b} = 0$. Then $D_b$
is a pseudo Yang-Mills field if and only if $i_N R^D = 0$ on
$H(\partial \Omega )$. \label{t:3}
\end{theorem}
\noindent
Here $T$ is the characteristic direction of $(\partial
\Omega , \theta )$, $\theta \equiv \frac{i}{2}
(\overline{\partial} -
\partial ) \varphi$, and $H(\partial \Omega )$ is the Levi
distribution. Also $N = - J T$ ($J$ is the complex structure on
$\mathbb{C}^n$). The proof relies on the explicit relationship
among the Levi-Civita connection $\nabla^g$ of $(\Omega  , g)$ and
the Graham-Lee connection $\nabla$ of $\varphi$ (cf.
\cite{kn:GrLe} and our Appendix A for the description and main
properties of $\nabla$).
\par Urakawa has started (cf. \cite{kn:Ura1}-\cite{kn:Ura3}) a
study of Yang-Mills fields on $M$, that is of critical points of
the functional
\[ \mathcal{YM}(D) = \frac{1}{2} \int_M \| R^D \|^2 \, d \, {\rm
vol}(g_\theta ), \] where $d \, {\rm vol} (g_\theta )$ is the
canonical volume form associated to the Webster metric $g_\theta$
of $(M , \theta )$. As it will be shortly shown, $\mathcal{YM}$
and $\mathcal{PYM}$ are related. To motivate the definition of
$\mathcal{PYM}$ let $F_\theta$ be the Fefferman metric of $(M ,
\theta )$ (a Lorentz metric on $C(M)$, the total space of the
canonical circle bundle $\pi : C(M) \to M$ (cf. e.g. J.M. Lee,
\cite{kn:Lee})). By a result of E. Barletta et alt.,
\cite{kn:BaDrUr}, the base map $\phi : M \to N$ corresponding to
any smooth $S^1$-invariant harmonic map $\Phi : C(M) \to N$ from
$(C(M), F_\theta )$ into a Riemannian manifold $(N, g_N )$ is
locally a subelliptic harmonic map (in the sense of J. Jost \&
C-J. Xu, \cite{kn:JoXu}). Also $\phi$ is a critical point of the
functional $E(\phi ) = \frac{1}{2} \int_M trace_{G_\theta} \left(
\pi_H  \phi^* g_N \right)\, \theta \wedge (d \theta )^n$, where
$G_\theta$ is the Levi form. Here, if $B$ is a bilinear form on
$T(M)$ then $\pi_H B$ denotes the restriction of $B$ to $H(M)$,
the Levi distribution of $(M, T_{1,0}(M))$. The functional $E$
itself is obtained by integration along the fibre in the Dirichlet
functional $\mathbb{E}(\Phi ) = \frac{1}{2} \int_{C(M)}
trace_{F_\theta} (\Phi^* g_N ) \, d {\rm vol} (F_\theta )$, where
$\Phi = \phi \circ \pi$. Then perhaps subelliptic harmonic maps
(rather than harmonic maps, with respect to the Webster metric)
are the natural objects of study in CR geometry. Another example
of the sort is the {\em CR Yamabe problem}, i.e. given a contact
form $\theta$ on $M$ such that $G_\theta$ is positive definite,
find a contact form $\hat{\theta} = e^{u} \theta$, $u \in C^\infty
(M)$, such that the pseudohermitian scalar curvature $\hat{\rho}$
of $(M , \hat{\theta})$ is a constant $\lambda$. By a result of
J.M. Lee, \cite{kn:Lee}, the Fefferman metric changes conformally
$F_{\hat{\theta}} = e^{u \circ \pi} F_\theta$. Also the scalar
curvature $K : C(M) \to \mathbb{R}$ of $(C(M) , F_\theta )$ is
$S^1$-invariant and the corresponding base function $\pi_* K : M
\to \mathbb{R}$ is, up to a constant, the pseudohermitian scalar
curvature $\rho$ of $(M , \theta )$ (precisely $\pi_* K =
\frac{2n+1}{n+1} \rho$). Therefore, the CR Yamabe problem is
nothing but the Yamabe problem for the Fefferman metric and the
relevant equation (the Yamabe equation on $(C(M) , F_\theta )$)
projects on $c_n \Delta_b u + \rho u = \lambda u^{p-1}$ (the {\em
CR Yamabe equation}), a nonlinear subelliptic equation on $M$
(which may be analyzed with the techniques in \cite{kn:FoSt}, cf.
D. Jerison \& J.M. Lee, \cite{kn:JL1}-\cite{kn:JL2}, and N. Gamara
\& R. Yacoub, \cite{kn:GaYa}, for a complete solution to the CR
Yamabe problem). The common feature of the two examples above is
that both provide natural objects on $M$, as projections of
($S^1$-invariant) geometric quantities on $C(M)$, associated to
the Fefferman metric. A more refined statement is that both
examples lead to nonlinear subelliptic problems on $M$. This has
been already emphasized for the CR Yamabe problem. As to the
example of $S^1$-invariant harmonic maps $\Phi : C(M) \to N$, the
base map is a solution to $\Delta_b \phi^i +
g^{\alpha\overline{\beta}} T_\alpha (\phi^j )T_{\overline{\beta}}
(\phi^k ) ( {(\Gamma_N ) }^i_{jk} \circ \phi ) = 0$, where
${(\Gamma_N )}_{jk}^i$ are the Christoffel symbols of the second
kind of $g_N$. On the same line of thought, we may state the
following
\begin{theorem} Let $M$ be a compact strictly pseudoconvex CR manifold, of
CR dimension $n$. Let $\theta$ be a contact form on $M$ with
$G_\theta$ positive definite. Let $(E , \overline{\partial}_E )
\to M$ be a CR-holomorphic vector bundle and $h$ a Hermitian
metric in $E$. $i)$ There is a constant $c_n$ depending only on
the dimension and the orientation of $M$ such that
\begin{equation}
c_n \, \mathcal{YM}(D) = \mathcal{PYM}(D) + 2 \int_M \| i_T R^D
\|^2 \, \theta \wedge (d \theta )^n , \;\; D \in \mathcal{C}(E ,
h). \label{e:2}\end{equation} Consequently, given a Hermitian
connection $D$ in $E$ whose curvature $R^D$ is of type $(1,1)$,
$D$ is a pseudo Yang-Mills field on $M$ if and only if $D$ is the
Tanaka connection of $(E, \overline{\partial}_E , h)$. $ii)$ Let
$\widehat{\mathcal{YM}}(\mathbb{D}) = \frac{1}{2} \int_{C(M)}
\langle R^{\mathbb{D}} , R^{\mathbb{D}} \rangle \; d \, {\rm
vol}(F_\theta )$ be the Yang-Mills functional on $C(M)$, for
$\mathbb{D} \in \mathcal{C}(\pi^* E , \pi^* h)$. Then
\begin{equation}
\widehat{\mathcal{YM}}( \pi^* D ) = 2 \pi \mathcal{PYM} (D), \;\;
D \in \mathcal{C}(E , h). \label{e:1}\end{equation} Consequently,
if $\pi^* D$ is a Yang-Mills field on $(C(M), F_\theta )$ then $D$
is a pseudo Yang-Mills field on $M$. Viceversa, let $D$ be a
pseudo Yang-Mills field on $M$ such that $i_T R^D = 0$. Then
$\pi^*D$ is a Yang-Mills field on $C(M)$ if and only if
\begin{equation}
(R^{\alpha\overline{\beta}} - \frac{\rho}{2(n+1)} \,
g^{\alpha\overline{\beta}} ) R^D (T_\alpha , T_{\overline{\beta}}
) u  = 0, \label{e:01}
\end{equation}
for some local frame $\{ T_\alpha : 1 \leq \alpha \leq n \}$ of
$T_{1,0}(M)$ at any point $x \in M$, and
\begin{equation}
\Lambda_\theta R^D = 0. \label{e:02}
\end{equation}
In particular, if $M$ is $($pseudohermitian$)$ Ricci flat then the
pullback $\pi^* D$ of the canonical Tanaka connection $D$ of
$(E,h)$ is a Yang-Mills field. \label{t:1}
\end{theorem}
The main ingredients in the proof of Theorem \ref{t:1} are a local
coordinate calculation of the Fefferman metric of $(M , \theta )$,
the explicit relationship among the Levi-Civita connection
$\nabla^{C(M)}$ of $(C(M), F_\theta )$ and the Tanaka-Webster
connection $\nabla$ of $(M , \theta )$ (cf. Lemma \ref{l:3}), and
Theorem 2.3 in \cite{kn:Ura1}, p. 551. We may also state
(delegating the definitions to section 2)
\begin{theorem} Let $D$ be a pseudo Yang-Mills field and $D^t = D
+ A^t$, $|t| < \epsilon$, a smooth variation of $D$ whose first
order part $\varphi \equiv \{ d A^t /d t\}_{t=0}$ satisfies $i_T
\varphi = 0$ and $\delta^D_b \varphi = 0$. Then
\begin{equation}
\frac{d^2}{d t^2} \{ \mathcal{PYM}(D^t )\}_{t=0} = \int_M \langle
\mathcal{S}^D_b (\varphi ) \, , \, \varphi \rangle \, \theta
\wedge (d \theta )^n \, \label{e:t2}
\end{equation}
where $\mathcal{S}^D_b (\varphi ) \equiv \Delta^D_b \varphi +
\mathcal{R}^D_b (\varphi )$ and $\Delta^D_b \varphi \equiv d^D_b
\delta^D_b \varphi + \delta^D_b d^D_b \varphi$ is the generalized
sublaplacian.  The operator $\mathcal{S}^D_b : \Omega^{0,1}({\rm
Ad}(E)) \to \Omega^{0,1}({\rm Ad}(E))$ is subelliptic of order
$1/2$. \label{t:2}
\end{theorem}
As $\mathcal{R}^D_b$ is a zero order operator, the crucial point
in the proof of Theorem \ref{t:2} is to show that
\begin{equation}
(\Delta^D_b \varphi ) \otimes e_j = 2 \{ \square_b \varphi^i_j +
\label{e:t2.2}
\end{equation}
\[ +  (n-1) (\nabla_T \varphi^i_j +
\varphi^i_j \circ \tau ) \circ J \} \otimes e_i + lower \; order
\; terms, \] for any $\varphi \in \Omega^{0,1}({\rm Ad}(E))$,
$\varphi e_j = \varphi^i_j \otimes e_i$, and then exploit the
subellipticity of the Kohn-Rossi operator $\square_b$ on scalar
$(0,1)$-forms.

\section{CR and pseudohermitian geometry}
\subsection{Basic definitions and results}
Let $M$ be a $C^\infty$ manifold, of real dimension $(2n+1)$. A
complex subbundle $T_{1,0}(M) \subset T(M) \otimes \mathbb{C}$, of
complex rank $n$, is a {\em CR structure} on $M$ (of {\em CR
dimension} $n$) if
\[ T_{1,0}(M) \cap T_{0,1}(M) = (0), \]
\[ Z, W \in \Gamma^\infty (T_{1,0}(M)) \Longrightarrow [Z,W] \in
\Gamma^\infty (T_{1,0}(M)). \] Here $T_{0,1}(M) =
\overline{T_{1,0}(M)}$ is the complex conjugate of $T_{1,0}(M)$.
Also, if $E \to M$ is a vector bundle then $\Gamma^\infty (E)$
denotes the space of $C^\infty$ sections in $E$ (eventually
defined on some open set $U \subseteq M$, to be understood from
the context). The {\em tangential Cauchy-Riemann operator}
\[ \overline{\partial}_b : C^\infty (M) \to \Gamma^\infty
(T_{0,1}(M)^* ) \] is given by $(\overline{\partial}_b f)
\overline{Z} = \overline{Z}(f)$, for any $C^\infty$ function $f :
M \to \mathbb{C}$ and any $Z \in T_{1,0}(M)$. Let $E \to M$ be a
complex vector bundle over a CR manifold. A {\em
pre-$\overline{\partial}$-operator} is a first order differential
operator
\[ \overline{\partial}_E  : \Gamma^\infty (E) \to \Gamma^\infty
(T_{0,1}(M)^* \otimes E) \] such that
\[ \overline{\partial}_E (f u) = f \overline{\partial}_E u +
(\overline{\partial}_b f) \otimes u, \] for any $f \in C^\infty
(M)$ and any $u \in \Gamma^\infty (E)$. A pair $(E ,
\overline{\partial}_E )$ consisting of a complex vector bundle and
a pre-$\overline{\partial}$-operator is a {\em CR-holomorphic}
vector bundle if $\overline{\partial}_E$ satisfies the {\em
integrability condition}
\[ [\overline{Z} , \overline{W}] \cdot u = \overline{Z} \cdot
\overline{W} \cdot u - \overline{W} \cdot \overline{Z} \cdot u, \]
for any $u \in \Gamma^\infty (E)$, $Z, W \in T_{1,0}(M)$. Here
$\overline{Z} \cdot u$ is short for $(\overline{\partial}_E u)
\overline{Z}$.
\par
Let $H(M) = {\rm Re} \{ T_{1,0}(M) \oplus T_{0,1}(M) \}$ be the
{\em Levi distribution} and $J : H(M) \to H(M)$, $J(Z +
\overline{Z}) = i (Z - \overline{Z})$, $Z \in T_{1,0}(M)$, its
complex structure ($i = \sqrt{-1}$). When $M$ is oriented, which
is assumed throughout this paper, the conormal bundle $H(M)^\bot_x
= \{ \omega \in T^*_x (M) : {\rm Ker}(\omega ) \supseteq H(M)_x
\}$, $x \in M$, is an oriented real line bundle, hence trivial
($H(M)^\bot \approx M \times \mathbb{R}$, a vector bundle
isomorphism). Therefore $H(M)^\bot \to M$ admits globally defined
nowhere zero sections $\theta \in \Gamma^\infty (H(M)^\bot )$,
each of which is referred to as a {\em pseudohermitian structure}
on $M$. The {\em Levi form} is
\[ L_\theta (Z , \overline{W}) = - i (d \theta )(Z ,
\overline{W}), \] for any $Z, W \in T_{1,0}(M)$. $(M ,
T_{1,0}(M))$ is {\em nondegenerate} if $L_\theta$ is nondegenerate
for some $\theta$.  If this is the case, each pseudohermitian
structure $\theta$ is a {\em contact form}, i.e. $\theta \wedge (d
\theta )^n$ is a volume form on $M$. Two pseudohermitian
structures $\theta , \hat{\theta} \in \Gamma^\infty (H(M)^\bot )$
are related by $\hat{\theta} = f \theta$, for some $C^\infty$
function $f : M \to \mathbb{R} \setminus \{ 0 \}$. Then
$L_{\hat{\theta}} = f L_\theta$, hence nondegeneracy is a CR
invariant notion (i.e. invariant under a transformation $\theta
\mapsto f \theta$ of the pseudohermitian structure). Let $T$ is
the unique nowhere zero globally defined tangent vector field on
$M$, transverse to the Levi distribution, determined by $\theta
(T) = 1$ and $i_T \, d \theta = 0$ (the {\em characteristic
direction} of $d \theta$). Also, let us consider the
semi-Riemannian metric $g_\theta$  (the {\em Webster metric} of
$(M , \theta )$) given by
\[ g_\theta (X, Y) = G_\theta (X,Y),
\;\; g_\theta (X , T) = 0, \;\; g_\theta (T,T) = 1, \] where
$G_\theta (X,Y) = (d \theta )(X, J Y)$, $X,Y \in H(M)$, is the
(real) {\em Levi form} (note that $L_\theta$ and (the
$\mathbb{C}$-linear extension of) $G_\theta$ coincide on
$T_{1,0}(M) \otimes T_{0,1}(M)$). $(M , T_{1,0}(M))$ is {\em
strictly pseudoconvex} if $L_\theta$ is positive definite for some
$\theta$. For instance, if $M = \{ (x,y,u) \in \mathbb{R}^3 : y
\neq 0 \}$ is endowed with the CR structure given by (\ref{e:gen})
in the Introduction then a calculation shows that the
characteristic direction (corresponding to the contact form
$\theta = (1/y) \, d x - y \, d z$) is
\[ T = \frac{y(3 + y^2 )}{4(1 + y^2 )} \, \frac{\partial}{\partial
x} - \frac{i}{8} \,  y (1 - y^2 ) \, \frac{\partial}{\partial y} -
\frac{1 + 3 y^2}{4 y (1 + y^2 )} \, \frac{\partial}{\partial z} \]
so that
\[ [T_1 , T_{\overline{1}}] = \frac{i}{2} \, T_1 - \frac{i}{2} \,
T_{\overline{1}} - \frac{2i}{1 + y^2} \, T \] (where
$T_{\overline{1}} = \overline{T}_1$). Consequently \[ L_\theta
(T_1 , T_{\overline{1}}) = (i/2) \theta ([T_1 , T_{\overline{1}}])
= 1/(1 + y^2 ), \] hence $M$ is strictly pseudoconvex. A
fundamental result in pseudohermitian geometry (established
independently by N. Tanaka, \cite{kn:Tan}, and S. Webster,
\cite{kn:Web}) is that on any nondegenerate CR manifold on which a
contact form $\theta$ has been fixed there is a unique linear
connection $\nabla$ (the {\em Tanaka-Webster connection} of $(M ,
\theta )$) such that i) $H(M)$ is parallel with respect to
$\nabla$, ii) $\nabla g_\theta = 0$, $\nabla J = 0$, and iii) the
torsion $T_\nabla$ of $\nabla$ is {\em pure}, i.e.
\[ T_\nabla (Z , W) = 0, \;\; T_\nabla (Z ,
\overline{W}) = 2 i L_\theta (Z , \overline{W}) T, \;\; Z,W \in
T_{1,0}(M), \]
\[ \tau \circ J + J \circ \tau = 0, \]
where $\tau (X) = T_\nabla (T , X)$, $X \in T(M)$, is the {\em
pseudohermitian torsion}. If $M$ is $3$-dimensional ($n = 1$) and
$T_1$ is a local generator of the CR structure we set \[
\nabla_{T_1} T_1 = \Gamma^1_{11} T_1 \, , \;\;
\nabla_{T_{\overline{1}}} T_1 = \Gamma^1_{\overline{1}1} T_1 \, ,
\;\; \nabla_{T} T_1  = \Gamma^1_{01} T_1 \, . \] A calculation
(based on (i)-(iii)) shows that \begin{equation}\label{e:web1}
\Gamma^1_{11} = g^{1\overline{1}} \{ T_1 ( g_{1\overline{1}}) -
g_\theta (T_1 , [T_1 , T_{\overline{1}}]) \} , \end{equation}
\begin{equation}
\Gamma^1_{\overline{1}1} = g^{1\overline{1}} g_\theta
([T_{\overline{1}} , T_1 ] , T_{\overline{1}}),
\end{equation}
\begin{equation}
\label{e:web3} \Gamma^1_{01} = g^{1\overline{1}} g_\theta ([T,T_1
] , T_{\overline{1}}).
\end{equation}
Here $g_{1\overline{1}} = L_\theta (T_1 , T_{\overline{1}})$ and
$g^{1\overline{1}} = 1/g_{1\overline{1}}$. Going back to the
example $\mathbb{R}^3 \setminus \{ y = 0 \}$ with the CR structure
(\ref{e:gen}) we have
\[ [T, T_1 ] = \frac{i}{8} \, (1 - y^2 ) \, T_1 + \frac{i}{8} \,
(1 + y^2 ) \, T_{\overline{1}} \] hence (by
(\ref{e:web1})-(\ref{e:web3}))
\[ \Gamma^1_{11} = i \left( \frac{1}{2} + \frac{y^2}{1 + y^2}
\right) , \;\; \Gamma^1_{\overline{1}1} = - \frac{i}{2} \, (1 +
y^2 ), \;\; \Gamma^1_{01} = \frac{i}{8} \, (1 - y^2 ). \] We
assume from now on that, unless otherwise stated, $M$ is strictly
pseudoconvex. A complex valued differential $p$-form $\eta$ on $M$
is {\em of type} $(p,0)$ (or a $(p,0)$-{\em form} on $M$) if
$T_{0,1}(M) \, \rfloor \, \eta = 0$. Let $\theta$ be a contact
form on $M$ and $T$ the characteristic direction of $d \theta$.
 Let $\{ T_\alpha : 1 \leq \alpha \leq n
\}$ be a local frame in $T_{1,0}(M)$, defined on an open set $U
\subseteq M$. Let $\{ \theta^\alpha : 1 \leq \alpha \leq n \}$ be
the corresponding {\em admissible coframe}, i.e, the (locally
defined) complex $1$-forms determined by $\theta^\alpha (T_\beta )
= \delta^\alpha_\beta$, $\theta^\alpha (T_{\overline{\beta}} ) =
0$, and $\theta^\alpha (T) = 0$. Here $T_{\overline{\beta}} =
\overline{T_\beta}$. Then $\{ \theta^\alpha ,
\theta^{\overline{\alpha}} , \theta \}$ is a (local) frame of $T^*
(M) \otimes \mathbb{C}$ on $U$ and a $(p,0)$-form $\eta$ on $M$
may be locally expressed as sums of monomials of the form
$\theta^{\alpha_1} \wedge \cdots \wedge \theta^{\alpha_p}$ or
$\theta \wedge \theta^{\alpha_1} \wedge \cdots \wedge
\theta^{\alpha_{p-1}}$ (with $C^\infty (U)$-coefficients).
Therefore, the top degree complex forms $\eta$ such that
$T_{1,0}(M) \, \rfloor \, \eta = 0$ are (unlike the case of
complex manifolds, where the top degree is the complex dimension)
the forms of type $(n+1,0)$ (where $n$ is the CR dimension). Let
$K(M) = \Lambda^{n+1,0}(M) \to M$ be the bundle of $(n+1,0)$-forms
on $M$ (the {\em canonical line bundle}). There is a natural
action of $\mathbb{R}^+ = (0, + \infty )$ on $K(M) \setminus \{
zero \; section \}$. Let $C(M)$ be the quotient space and $\pi :
C(M) \to M$ the projection. Then $C(M) \to M$ is a principal
$S^1$-bundle (the {\em canonical circle bundle}). Its locally
trivial structure is described by
\[ \pi^{-1}(U) \to U \times S^1 , \;\; [\omega ] \mapsto (x \; ,
\; \lambda /|\lambda |) , \]
\[ \omega = \lambda (\theta \wedge \theta^1 \wedge \cdots \wedge
\theta^n )_x , \; x \in U, \; \lambda \in \mathbb{C}^* =
\mathbb{C} \setminus \{ 0 \} . \] We shall need the local fibre
coordinate \[ \gamma : \pi^{-1}(U) \to \mathbb{R}, \;\; \gamma
([\omega ]) = \arg (\lambda /|\lambda | ), \] where $\arg : S^1
\to [0, 2 \pi ).$ Let $(E , \overline{\partial}_E ) \to M$ be a
CR-holomorphic vector bundle. Let $h$ be a Hermitian metric in
$E$. Let $\mathcal{C}(E , h)$ be the affine space of all
connections $D$ in $E$ such that $D h = 0$, i.e.
\[ X(h(u,v)) = h(D_X u , v)  + h(u , D_{\overline{X}} v), \]
for any $X \in T(M) \otimes \mathbb{C}$ and any $u,v \in
\Gamma^\infty (E)$. A connection $D \in \mathcal{C}(E , h)$ is
{\em Hermitian} if $D^{0,1} = \overline{\partial}_E$. Here
$D^{0,1} u$ is the restriction of $D u$ to $T_{0,1}(M)$. Let ${\rm
Ad} (E) \to M$ be the subbundle of ${\rm End}(E) \to M$ consisting
of all skew-symmetric endomorphisms $S$, i.e. $h(S u , v) + h(u ,
S v) = 0$, for any $u, v \in \Gamma^\infty (E)$. By a result in
\cite{kn:DrUr}, p. 43, given a contact form $\theta$ and an
endomorphism $S \in \Gamma^\infty ({\rm Ad}(E))$ there is a unique
Hermitian connection $D = D(h, \theta , S)$ in $E$ (the {\em
canonical $S$-connection}) such that
\begin{equation}
\Lambda_\theta \; R^D = 2 n S.
\label{e:axS}
\end{equation}
Here $R^D = D \circ D : \Omega^0 (E) \to \Omega^2 (E)$ is the
curvature $2$-form of $D$. Also we set $\Omega^k (E) =
\Gamma^\infty (\Lambda^k T^* (M) \otimes E)$, $k \geq 0$. If $F
\to M$ is a vector bundle and $\varphi \in \Gamma^\infty (T^* (M)
\otimes T^* (M) \otimes F)$ the trace $\Lambda_\theta \varphi$ of
$\varphi$ is given by
\[ i (\Lambda_\theta \varphi )_x = \sum_{\alpha =1}^n \varphi
(Z_\alpha , Z_{\overline{\alpha}})_x \, , \] where $\{ Z_\alpha
\}$ is a (local) orthonormal (i.e. $L_\theta (Z_\alpha ,
Z_{\overline{\beta}}) = \delta_{\alpha\beta}$) frame of
$T_{1,0}(M)$ on $U \ni x$. Therefore $\Lambda_\theta \varphi \in
\Gamma^\infty (F)$. When $S = 0$ the canonical $S$-connection is
the {\em Tanaka connection} $D(h, \theta , 0)$ in $E \to M$ (cf.
\cite{kn:Tan}). $D(h , \theta , 0)$ is a CR invariant. Assume $M$
to be compact. The {\em Yang-Mills functional} $\mathcal{YM} :
\mathcal{C}(E,h) \to [0, + \infty )$ is given by
\[ \mathcal{YM}(D) = \frac{1}{2} \int_M \| R^D \|^2 \, \theta
\wedge (d \theta )^n . \] A {\em Yang-Mills field} on $M$ is a
critical point $D \in \mathcal{C}(E , h)$ of $\mathcal{YM}$, i.e.
a solution to the Yang-Mills equations
\begin{equation}
\delta^D R^D = 0.
\label{e:3}
\end{equation}
Let $\Omega$ be a differential $2$-form on $M$. Then $\Omega$ is
{\em of type} $(1,1)$ if $\Omega (Z , W) = 0$, $\Omega
(\overline{Z}, \overline{W}) = 0$, for any $Z, W \in T_{1,0}(M)$,
and $i_T \, \Omega = 0$. Let $D \in \mathcal{C}(E, h)$ be a
Hermitian connection such that its curvature $R^D$ is a form of
type $(1,1)$. By a result in \cite{kn:Ura1}, $D$ is a Yang-Mills
field if and only if $D$ is the Tanaka connection $D(h, \theta ,
0)$. In general, canonical $S$-connections solve the inhomogeneous
Yang-Mills equations $\delta^D R^D = f$, in the presence of
suitable compatibility conditions satisfied by $f$ (cf. Theorem 2
in \cite{kn:DrUr}, p. 44-45).

\subsection{Symmetric CR structures}
 \par The CR structure
$T_{1,0}(M)$ is {\em symmetric} if there is $X \in \mathcal{X}(M)$
such that \[ \mathcal{L}_X \theta = t \, \theta , \;\;
\mathcal{L}_X \theta^\alpha = w^\alpha_\beta \, \theta^\beta +
\ell^\alpha \, \theta , \] for some functions $t, \,
w^\alpha_\beta , \, \ell^\alpha$ on $M$ ($t$ real valued) and $X$
is a {\em symmetry} of $T_{1,0}(M)$. If \[ \hat{\theta} = e^u \,
\theta , \;\; \hat{\theta}^\alpha = U^\alpha_\beta \, \theta^\beta
+ v^\alpha \, \theta \, , \] (where $[U^\alpha_\beta ]$ is ${\rm
GL}(n, \mathbb{C})$-valued) and $X$ is a symmetry of the CR
structure then
\begin{equation} \mathcal{L}_X \hat{\theta} = \hat{t} \;
\hat{\theta} , \;\; \hat{t} \equiv t + X(u), \label{e:n4}
\end{equation}
\begin{equation}
\mathcal{L}_X \, \hat{\theta}^\alpha = \hat{w}^\alpha_\beta \,
\hat{\theta}^\beta + \hat{\ell}^\alpha \, \hat{\theta} \, ,
\label{e:n5}
\end{equation}
\[ \hat{w}^\alpha_\beta \equiv (U^{-1})^\gamma_\beta \{ X
(U^\alpha_\gamma ) + U^\alpha_\rho w^\rho_\gamma \} , \]
\[ \hat{\ell}^\alpha \equiv e^{-u} \{ X(v^\alpha ) + U^\alpha_\beta \ell^\beta
+ v^\alpha \, t  - (U^{-1})^\gamma_\rho v^\rho [X(U^\alpha_\gamma
) + U^\alpha_\beta w^\beta_\gamma ] \} .
\] In particular (\ref{e:n4})-(\ref{e:n5}) show that the notion of
symmetric CR structure is globally defined. Assume from now on
that the CR structure $T_{1,0}(M)$ admits $2n+1$ linearly
independent symmetries $\tilde{X}_1 , \cdots , \tilde{X}_{2n+1}
\in \mathcal{X}(M)$ such that $[\tilde{X}_i , \tilde{X}_j ] =
c^k_{ij} \tilde{X}_k$, for some $c^k_{ij} \in \mathbb{R}$.
\begin{proposition} $(${\rm P. Nurowski, \cite{kn:Nur}}$)$
\par\noindent Let $M$ be a strictly pseudoconvex CR manifold with
$H^1 (M ; \mathbb{R}) = 0$. There is a transformation $\{ \theta ,
\theta^\alpha \} \mapsto \{ \Omega , \Omega^\alpha \}$ of the form
\begin{equation}
\Omega = e^u \theta , \;\; \Omega^\alpha = U^\alpha_\beta
\theta^\beta + v^\alpha \theta ,
\label{e:trans}
\end{equation}
where $[U^\alpha_\beta ]$ is ${\rm GL}(n, \mathbb{C})$-valued,
such that
\[ \mathcal{L}_{\tilde{X}_i} \Omega = 0, \;\;
\mathcal{L}_{\tilde{X}_i} \Omega^\alpha = 0, \;\; 1 \leq i \leq
2n+1. \] \label{p:n2}
\end{proposition}
\noindent Here $H^1 (M ; \mathbb{R})$ is the first de Rham
cohomology group. Its vanishing guarantees that the solution $u$
to (\ref{e:n6}) is globally defined. \par\noindent {\em Proof of
Proposition $\ref{p:n1}$}. As $\tilde{X}_i$ are symmetries of the
CR structure
\[ \mathcal{L}_{\tilde{X}_i} \theta = t_i \theta , \;\;
\mathcal{L}_{\tilde{X}_i} \theta^\alpha = w^\alpha_{i\beta}
\theta^\beta + \ell^\alpha_i \theta . \] We must solve the system
of first order linear PDEs
\begin{equation}
t_i + \tilde{X}_i (u) = 0, \label{e:n6}
\end{equation}
\begin{equation}
\tilde{X}_i (U^\alpha_\beta ) + U^\alpha_\gamma w^\gamma_{i\beta}
= 0 , \label{e:n7}
\end{equation}
\begin{equation}
\tilde{X}_i (v^\alpha ) + U^\alpha_\beta \ell^\beta_i + v^\alpha
t_i = 0, \label{e:n8}
\end{equation}
with the unknowns $u$, $U^\alpha_\beta$ and $v^\alpha$. Let $\eta
\in \Omega^1 (M)$ be defined by $\eta (\tilde{X}_i ) = t_i$, $1
\leq i \leq 2n+1$. Then (\ref{e:n6}) may be written $d u + \eta =
0$. We have
\[ \mathcal{L}_{\tilde{X}_i} \mathcal{L}_{\tilde{X}_j} \theta =
\mathcal{L}_{\tilde{X}_j} \mathcal{L}_{\tilde{X}_i} \theta +
\mathcal{L}_{[\tilde{X}_i , \tilde{X}_j ]} \theta \] hence
\[ \tilde{X}_i (t_j ) - \tilde{X}_j (t_i ) - c^k_{ij} t_k = 0
\] that is $d \eta = 0$. Thus there is a globally defined real valued function
$g \in C^\infty (M)$ such that $\eta = d g$ and $u \equiv - g$
solves (\ref{e:n6}). Next, we consider the (locally defined)
$1$-forms $\eta^\alpha_\beta$ and $\eta^\alpha$ given by \[
\eta^\alpha_\beta (\tilde{X}_i ) = w^\alpha_{i\beta} \, , \;\;
\eta^\alpha (\tilde{X}_i ) = \ell^\alpha_i \, , \;\; 1 \leq i \leq
2n+1. \] Then (\ref{e:n7})-(\ref{e:n8}) may be written
\begin{equation}
d U^\alpha_\beta + U^\alpha_\gamma \, \eta^\gamma_\beta = 0,
\label{e:n9}
\end{equation}
\begin{equation}
d v^\alpha + v^\alpha \, \eta + U^\alpha_\beta \, \eta^\beta  = 0.
\label{e:n10}
\end{equation}
Assuming that (\ref{e:n9}) has been solved in a neighborhood $U$
of each point, let us solve (\ref{e:n10}). Multiplying in both
sides by $e^{-u}$ (where $u$ is a solution to (\ref{e:n6})) leads
to
\[ d(e^{-u} v^\alpha ) + e^{-u} U^\alpha_\beta \eta^\beta = 0. \]
Therefore, to prove existence of a (local) solution $v^\alpha$ to
(\ref{e:n10}) it suffices to show that $e^{-u} U^\alpha_\beta
\eta^\beta$ is exact (in a neighborhood of a point). The identity
\[ \mathcal{L}_{\tilde{X}_i} \mathcal{L}_{\tilde{X}_j} \theta^\alpha =
\mathcal{L}_{\tilde{X}_j} \mathcal{L}_{\tilde{X}_i} \theta^\alpha
+ \mathcal{L}_{[\tilde{X}_i , \tilde{X}_j ]} \theta^\alpha \]
yields
\begin{equation}
d \eta^\alpha_\beta = \eta^\alpha_\gamma \wedge \eta^\gamma_\beta
\, , \label{e:n11}
\end{equation}
\begin{equation}
d \eta^\alpha = \eta^\alpha_\beta \wedge \eta^\beta + \eta^\alpha
\wedge \eta . \label{e:n12}
\end{equation} Let $U^\alpha_\beta$ be a solution to (\ref{e:n9}).
Then (by (\ref{e:n12}))
\[ d(e^{-u} U^\alpha_\beta \eta^\beta ) = e^{-u} \{ d
U^\alpha_\beta \wedge \eta^\beta + U^\alpha_\beta \, d \eta^\beta
- U^\alpha_\beta \, d u \wedge \eta^\beta \} = \]
\[ = e^{-u} U^\alpha_\beta \{ d \eta^\beta - \eta^\beta_\gamma
\wedge \eta^\gamma - \eta^\beta \wedge \eta \} = 0. \] Thus there
is a function $f^\alpha \in C^\infty (U)$ such that $e^{-u}
U^\alpha_\beta \eta^\beta = d f^\alpha$ and $v^\alpha \equiv - e^u
f^\alpha$ solves (\ref{e:n10}). To solve (\ref{e:n9}) let $(U, x^i
)$ be a normal coordinate neighborhood at a point $x_0 \in M$ (we
think of $M$ as a Riemannian manifold with the Webster metric
$g_\theta$). We shall show that for any $c^\alpha_\beta \in
\mathbb{C}$ there is a unique solution to (\ref{e:n9}) with the
initial condition $U^\alpha_\beta (x_0 ) = c^\alpha_\beta$. Let $a
= (a^1 , \cdots , a^{2n+1}) \in U$ be an arbitrary point and let
us consider the geodesic $a_t = (a^1 t , \cdots , a^{2n+1} t)$.
Let $f^\alpha_\beta (t)$ be the solution to the Cauchy problem for
the system of ODEs
\[
\frac{d f^\alpha_\beta}{d t} + f^\alpha_\gamma (t) \,
\eta^\gamma_\beta (\dot{a}_t ) = 0, \] with the initial condition
$f^\alpha_\beta (0) = c^\alpha_\beta$, where $\dot{a}_t$ is the
tangent vector at $a_t$. We define $U^\alpha_\beta \in C^\infty
(U)$ by setting $U^\alpha_\beta (a) = f^\alpha_\beta (1)$. Of
course, if we start with $\det (c^\alpha_\beta ) \neq 0$ then
$[U^\alpha_\beta ]$ is ${\rm GL}(m, \mathbb{C})$-valued on a
neighborhood of $x_0$. We wish to show that $U^\alpha_\beta$
satisfies (\ref{e:n9}), i.e.
\begin{equation}
Y (U^\alpha_\beta ) + U^\alpha_\gamma (a) \, \eta^\gamma_{\beta \,
, \, a} (Y) = 0, \label{e:n14}
\end{equation}
for any $Y = b^j (\partial /\partial x^j )_a \in T_a (M)$. We
start by extending $Y$ to the vector field $Y = b^j \; \partial
/\partial x^j$ with constant components $b^j$ on $U$. Similarly,
let us extend the vector field $\dot{a}_t$ along the geodesic
$a_t$ to the vector field $X = a^i \; \partial /\partial x^i$. We
shall show that along $a_t$
\begin{equation}
X\left( Y(U^\alpha_\beta ) + U^\alpha_\rho \, \eta^\rho_\beta (Y)
\right) + \left( Y(U^\alpha_\gamma ) + U^\alpha_\rho \,
\eta^\rho_\gamma (Y) \right) \; \eta^\gamma_\beta (X) = 0.
\label{e:n15}
\end{equation}
When this is done, we see that $Y(U^\alpha_\beta ) + U^\alpha_\rho
\, \eta^\rho_\beta (Y)$, clearly satisfying $Y(U^\alpha_\beta ) +
U^\alpha_\rho \, \eta^\rho_\beta (Y) = 0$ at $x_0$, must be the
zero function (which satisfies (\ref{e:n15}) with the same initial
condition), i.e. (\ref{e:n14}) is verified. It remains that we
prove (\ref{e:n15}). This follows from (\ref{e:n11}). Indeed (as
$[X,Y] = 0$)
\[ X(\eta^\alpha_\beta (Y)) = Y(\eta^\alpha_\beta (X)) +
\eta^\alpha_\gamma (X) \eta^\gamma_\beta (Y) - \eta^\alpha_\gamma
(Y) \eta^\gamma_\beta (X) \] and the proof of (\ref{e:n15}) is
straightforward. Proposition \ref{p:n2} is a first step towards
recovering the methods of P. Nurowski,\cite{kn:Nur} [eventually
leading to local solutions of the Yang-Mills equations on $(C(M),
F_\theta )$] as mentioned in the Introduction. The result in
Proposition \ref{p:n2} may be refined to show that there is a
coframe $\{ \Omega , \Omega^\alpha \}$ such that
\[ \mathcal{L}_{\tilde{X}_i} \Omega = 0, \;\;
\mathcal{L}_{\tilde{X}_i} \Omega^\alpha = 0, \;\; d \Omega = 2 i
\sum_{\alpha =1}^n \Omega^\alpha \wedge \Omega^{\overline{\alpha}}
\]
(compare to (\ref{e:i2}) in the Introduction). The proof is
illustrative of the local methods in pseudohermitian geometry. Let
$\{ \Omega , \Omega^\alpha \}$ be the $1$-forms furnished by
Proposition \ref{p:n2}, given by a transformation of the form
(\ref{e:trans}). If $\{ T , T_\alpha \}$ is such that $\theta (T)
= 1$, $i_T \; d \theta = 0$ and $\theta^\alpha (T_\beta ) =
\delta^\alpha_\beta$, $\theta^\alpha (T_{\overline{\beta}}) = 0$,
$\theta^\alpha (T) = 0$, let us set $W_\alpha =
(U^{-1})^\beta_\alpha T_\beta$. One may easily show that
\begin{equation}
d \Omega = 2 i G_{\alpha\overline{\beta}} \; \Omega^\alpha \wedge
\Omega^{\overline{\beta}} +  \Phi \wedge \Omega ,\label{e:loc27}
\end{equation}
where $G_{\alpha\overline{\beta}} = e^u (U^{-1})^\gamma_\alpha
U^{\overline{\rho}}_{\overline{\beta}} g_{\gamma\overline{\rho}}$
and $\Phi = e^{-u} \{ W_\alpha (u) \Omega^\alpha +
W_{\overline{\alpha}} (u) \Omega^{\overline{\alpha}} \}$. By
$\mathcal{L}_X = d \circ i_X + i_X \circ d$ it follows that
$\mathcal{L}_{\tilde{X}_i} d \Omega = 0$. Taking the Lie
derivative of (\ref{e:loc27}) gives
\[ 0 = \mathcal{L}_{\tilde{X}_i} d \Omega = 2 \sqrt{-1}
\mathcal{L}_{\tilde{X}_i} (G_{\alpha\overline{\beta}})
\Omega^\alpha \wedge \Omega^{\overline{\beta}} +
(\mathcal{L}_{\tilde{X}_i} \Phi ) \wedge \Omega \] hence (as
$\mathcal{L}_{\tilde{X}_i} \Phi \equiv 0$, $\bmod \; \Omega$) on
one hand $\tilde{X}_i (G_{\alpha\overline{\beta}}) = 0$, i.e.
$G_{\alpha\overline{\beta}} = a_{\alpha\overline{\beta}} \in
\mathbb{C}$, and on the other $\mathcal{L}_{\tilde{X}_i} \Phi =
0$. The latter may be written $d u_\alpha - u_\alpha \; d u = 0$
(where $u_\alpha = W_\alpha (u)$). Hence $d(e^{-u} u_\alpha ) =
0$, i.e. $u_\alpha = c_\alpha e^u$, for some $c_\alpha \in
\mathbb{C}$. Therefore $\Phi = c_\alpha \Omega^\alpha +
c_{\overline{\alpha}} \Omega^{\overline{\alpha}}$. Finally, let
$[b^\alpha_\beta ]$ be a square root of
$[a_{\alpha\overline{\beta}}]$ (as $[a_{\alpha\overline{\beta}}]$
is positive definite) and consider the transformation
\[ \hat{\Omega}^\alpha = b^\alpha_\beta \Omega^\beta + \frac{i}{2}
\; c_{\overline{\beta}}
(b^{-1})^{\overline{\beta}}_{\overline{\alpha}} \Omega . \] Then
$d \Omega = 2 i \sum_{\alpha =1}^n \hat{\Omega}^\alpha \wedge
\hat{\Omega}^{\overline{\alpha}}$ and $\mathcal{L}_{\tilde{X}_i}
\hat{\Omega}^\alpha = 0$. Q.e.d.

\section{Boundary values of Yang-Mills fields}
Let $\Omega \subset \mathbb{C}^n$ be a bounded domain with smooth
boundary $\partial \Omega$, i.e. there is a neighborhood $U
\supset \overline{\Omega}$ and a real valued function $\varphi \in
C^\infty (U)$ such that $\Omega = \{ z \in U : \varphi (z) < 0
\}$, $\partial \Omega = \{ z \in U : \varphi (z) = 0 \}$, and
$\nabla \varphi (z) \neq 0$, for any $z \in \partial \Omega$. We
assume that $\Omega$ is strictly pseudoconvex, i.e. $\partial
\Omega$ is a strictly pseudoconvex CR manifold (with the natural
CR structure $T_{1,0}(\partial \Omega ) = T^{1,0}(\mathbb{C}^n )
\cap [T(\partial \Omega )\otimes \mathbb{C}]$ induced by the
complex structure of the ambient space). \par Let $\pi : F \to U$
be a holomorphic vector bundle. The portion $E = \pi^{-1}(\partial
\Omega )$ of $F$ over the boundary of $\Omega$ is CR-holomorphic.
Indeed, as $F$ is holomorphic, there is a natural differential
operator
\[ \overline{\partial}_F : \Gamma^\infty (F) \rightarrow
\Gamma^\infty (T^{0,1}(U)^* \otimes F) \] where $T^{0,1}(U)$ is
the anti-holomorphic tangent bundle over $U$. Given $u \in
\Gamma^\infty (E)$ let $\tilde{u} \in \Gamma^\infty (F)$ be a
$C^\infty$ extension of $u$ as a cross-section in $F$ and set
$\left( \overline{\partial}_E u \right)_z = \left(
\overline{\partial}_F \tilde{u}\right)_z$ for any $z \in
\partial \Omega$. The definition of $\left( \overline{\partial}_E
u\right)_z$ does not depend upon the choice of extension
$\tilde{u}$ of $u$ because $\left. (\overline{\partial}
f)\right|_{T_{0,1}(\partial \Omega )} = \left.
\overline{\partial}_b (f\right|_{\partial \Omega})$ for any
$C^\infty$ function $f : U \rightarrow \mathbb{C}$. Let $\{
\Phi_\alpha : \pi^{-1} (\Omega_\alpha ) \rightarrow \Omega_\alpha
\times \mathbb{C}^m : \alpha \in I \}$ be a trivialization atlas
for $F$ and $G_{\beta\alpha} : \Omega_\beta \cap \Omega_\alpha
\rightarrow {\rm GL}(m , \mathbb{C})$ the corresponding transition
functions. Set $U_\alpha = \Omega_\alpha \cap \partial \Omega$ and
$\left. g_{\beta\alpha} = {G_{\beta\alpha}}\right|_{U_{\alpha}
\cap U_{\beta}}$. As $G_{\beta\alpha}$ are holomorphic, it follows
that $E \rightarrow \partial \Omega$ is a peculiar type of
CR-holomorphic vector bundle (called {\em locally trivial} by C.
Le Brun, \cite{kn:Leb1}) in that its transition functions
$g_{\beta\alpha}$ are matrix valued CR functions on $\partial
\Omega$.
\par
Let $K(\zeta , z)$ be the Bergman kernel of $\Omega$. By a
classical result in \cite{kn:Fef}
\begin{equation}
K(\zeta , z) = c_\Omega |\nabla \varphi (z)|^2 \cdot \det
L_\varphi (z) \cdot \Psi (\zeta , z)^{-(n+1)} + H(\zeta , z),
\label{e:b1}
\end{equation}
(the {\em Fefferman asymptotic expansion formula} for the Bergman
kernel) where $H \in C^\infty (\overline{\Omega} \times
\overline{\Omega} \setminus \Delta )$, $\Delta$ is the diagonal of
$\partial \Omega \times
\partial \Omega$, and $H$ satisfies the estimate
\begin{equation}
|H(\zeta , z)| \leq c^\prime_\Omega |\Psi (\zeta , z)|^{-(n+1) +
1/2} \cdot |\log |\Psi (\zeta , z)||. \label{e:b2}
\end{equation}
Here $L_\varphi
=
\partial \overline{\partial} \varphi$. Also we set
\[ \Psi (\zeta , z) = (F(\zeta , z) - \varphi (z)) \chi (|\zeta -
z|) + (1 - \chi (|\zeta - z|)) |\zeta - z|^2 \] where
\[ F(\zeta , z) = - \sum_{j=1}^n \frac{\partial \varphi}{\partial
z^j} (z) (\zeta^j - z^j ) - \frac{1}{2} \sum_{j,k=1}^n
\frac{\partial^2 \varphi}{\partial z^j \partial z^k}(z) (\zeta^j -
z^j )(\zeta^k - z^k ) \] and $\chi (t)$ is a $C^\infty$ cut-off
function with $\chi (t) = 1$ for $|t| < \epsilon_0 /2$ and $\chi
(t) = 0$ for $|t| \geq 3 \epsilon_0 /4$. As a consequence of
(\ref{e:b1})
\[ K(z,z)^{-1/(n+1)} = |\varphi (z)| \left( \Phi (z) + H(z,z)
|\varphi (z)|^{n+1} \right)^{-1/(n+1)} \] where $\Phi (z) \equiv
c_\Omega |\nabla \varphi (z)|^2 \det L_\varphi (z)$ stays finite
near $\partial \Omega$ and (by (\ref{e:b2})) \[ |H(z,z)|\,
|\varphi (z)|^{n+1} \leq c^\prime_\Omega |\varphi (z)|^{1/2} |\log
|\varphi (z)|| \to 0, \;\;\; as \;\; z \to \partial \Omega . \]
Therefore $K(z,z)^{-1/(n+1)}$ vanishes at $\partial \Omega$. Also,
as $\Phi (z) \neq 0$ near the boundary, $\nabla K(z,z)^{-1/(n+1)}
\neq 0$ along $\partial \Omega$, hence $K(z,z)^{-1/(n+1)}$ may be
used as a defining function for $\Omega$. \par For the rest of
this section we assume that $\varphi (z) \equiv -
K(z,z)^{-1/(n+1)}$ and set $\theta \equiv
\frac{i}{2}(\overline{\partial} - \partial )\varphi$. Then $d
\theta = i \, \partial \overline{\partial} \varphi$. Let us
differentiate $\log |\varphi | = - (1/(n+1)) \log K$ (where $K$ is
short for $K(z,z)$) so that to obtain
\[ \frac{1}{\varphi} \, \overline{\partial} \varphi = -
\frac{1}{n+1} \, \overline{\partial} \log K. \] Applying the
operator $i \, \partial$ leads to
\begin{equation} \frac{1}{\varphi} \; d \theta -
\frac{i}{\varphi^2} \; \partial \varphi \wedge \overline{\partial}
\varphi = - \frac{i}{n+1} \, \partial \overline{\partial} \log K.
\label{e:b3}
\end{equation}
We shall need the Bergman metric
\[ g = \frac{\partial^2 \log K}{\partial z^j \partial
\overline{z}^k} \; d z^j \odot d \overline{z}^k \, . \] As well
known, $g$ is a K\"ahler metric on $\Omega$ (K\"ahler-Einstein
when $\Omega$ is homogeneous). Here $\odot$ denotes the symmetric
tensor product, i.e. $\alpha\odot \beta = \frac{1}{2} (\alpha
\otimes \beta + \beta \otimes \alpha )$. Let us set $\omega (X,Y)
= g(X , J Y)$ (the K\"ahler $2$-form of $(\Omega ,J, g)$, where
$J$ is the underlying complex structure). Then $\omega = - i \,
\partial \overline{\partial} \log K$ and (\ref{e:b3}) may be
written
\begin{equation}
g(X,Y) = \frac{n+1}{\varphi} \{ \frac{i}{\varphi} \, (\partial
\varphi \wedge \overline{\partial} \varphi )(X , J Y ) - d \theta
(X , J Y) \} , \label{e:b4}
\end{equation}
for any $X,Y \in \mathcal{X}(\Omega )$. \par We denote by
$M_\delta = \{ z \in \Omega : \varphi (z) = -\delta \}$ $(\delta >
0)$ the level sets of $\varphi$. For $\delta$ sufficiently small
$M_\delta$ is still a strictly pseudoconvex CR manifold (of CR
dimension $n-1$). Therefore, there is a one-sided neighborhood $V$
of $\partial \Omega$ which is foliated by the (strictly
pseudoconvex) level sets of $\varphi$. Let $\mathcal{F}$ be the
relevant foliation and let us denote by $H(\mathcal{F}) \to V$
(respectively by $T_{1,0}(\mathcal{F}) \to V$) the bundle whose
portion over $M_\delta$ is the Levi distribution $H(M_\delta )$
(respectively the CR structure $T_{1,0}(M_\delta )$). Note that
\[ T_{1,0}(\mathcal{F}) \cap T_{0,1}(\mathcal{F}) = (0), \]
\[ [\Gamma^\infty (T_{1,0}(\mathcal{F})) , \Gamma^\infty (T_{1,0}(\mathcal{F}))]
\subseteq \Gamma^\infty (T_{1,0}(\mathcal{F})). \] Here
$T_{0,1}(\mathcal{F}) = \overline{T_{1,0}(\mathcal{F})}$. By a
result in \cite{kn:LeMe}, there is a unique complex vector field
$\xi$ on $V$, of type $(1,0)$, such that $\partial \varphi (\xi )
= 1$ and $\xi$ is orthogonal to $T_{1,0}(\mathcal{F})$ with
respect to $\partial \overline{\partial} \varphi$, i.e. $\partial
\overline{\partial} \varphi (\xi , \overline{Z}) = 0$, for any $Z
\in T_{1,0}(\mathcal{F})$. We set $r \equiv 2 \, \partial
\overline{\partial} \varphi (\xi , \overline{\xi})$ ($r$ is the
{\em transverse curvature} of $\varphi$).  Let $\xi =
\frac{1}{2}(N - i T)$ be the real and imaginary parts of $\xi$.
Then
\[ (d \varphi )(N) = 2, \;\;\; (d \varphi )(T) = 0, \]
\[ \theta (N) = 0, \;\;\; \theta (T) = 1, \]
\[ \partial \varphi (N) = 1, \;\;\; \partial \varphi (T) = i. \]
In particular, $T$ is tangent to (the leaves of) $\mathcal{F}$.
$\mathcal{F}$ carries the tangential Riemannian metric $g_\theta$
(defined by (\ref{e:A.1}) in Appendix A). Note that the pullback
of $g_\theta$ to each leaf $M_\delta$ of $\mathcal{F}$ is the
Webster metric of $M_\delta$ (associated to the contact form
$j_\delta^* \theta$, where $j_\delta : M_\delta \subset V$). As a
consequence of (\ref{e:b4})
\begin{equation}
g(X,Y) = - \frac{n+1}{\varphi} \, g_\theta (X,Y), \;\;\; X,Y \in
H(\mathcal{F}). \label{e:b5}
\end{equation}
Also (by $J T = - N$ and (\ref{e:A.4}))
\begin{equation}
g(X,T) = 0, \;\; g(X, N) = 0, \;\;\; X \in H(\mathcal{F}),
\label{e:b6}
\end{equation}
\begin{equation}
g(T , N ) = 0, \;\; g(T,T) = g(N,N) = \frac{n+1}{\varphi} \left(
\frac{1}{\varphi} - r \right) . \label{e:b7}
\end{equation}
In particular $1 - r \varphi > 0$ everywhere in $\Omega$.  Using
(\ref{e:b5})-(\ref{e:b7}) we may relate the Levi-Civita connection
$\nabla^g$ of $(V , g)$ to the Graham-Lee connection $\nabla$ (cf.
Appendix A). By (\ref{e:b5}) (as $X(\varphi ) = 0$, $X \in
T(\mathcal{F})$)
\begin{equation}
g(\nabla^g_X Y , Z) = g (\nabla_X Y , Z), \;\; X,Y,Z \in
H(\mathcal{F}). \label{e:b8}
\end{equation}
Note that any tangent vector field $X \in T(V)$ decomposes as
\[ X = \pi_H X + \theta (X) T + \frac{1}{2} (d \varphi )(X) N, \]
($\pi_H : T(V) \to H(\mathcal{F})$ is the projection). By
(\ref{e:A.3}) $\theta ([T,X]) = 0$, $X \in H(\mathcal{F})$. Also
$[T, X] \in T(\mathcal{F})$, hence $[T,X] \in H(\mathcal{F})$, for
any $X \in H(\mathcal{F})$. Taking into account the identity
\begin{equation} 2g(\nabla^g_X Y , Z) = X(g(Y,Z)) + Y(g(X,Z)) -
Z(g(X,Y)) + \label{e:b9}
\end{equation}
\[ + g([X,Y],Z) + g([Z,X],Y) + g(X,[Z,Y]), \]
for any $X,Y,Z \in T(V)$, one has (by (\ref{e:b6}))
\[ 2 g(\nabla^g_X Y , T) = - T(g(X,Y)) + \]
\[ + g([X,Y],T) + g([T,X],Y) + g(X,[T,Y]) =  \]
\[ = \frac{n+1}{\varphi} \, \left\{ T(g_\theta (X,Y)) -
g_\theta ([T,X], Y) - g_\theta (X , [T, Y]) \right\} +  \] \[ +
\frac{n+1}{\varphi} \left( \frac{1}{\varphi} - r \right) \theta
([X,Y]), \] for any $X,Y \in H(\mathcal{F})$. By
(\ref{e:A.13})-(\ref{e:A.14}) and $\nabla_X Y \in H(\mathcal{F})$
it follows that
\[ T(g_\theta (X,Y)) - g_\theta ([T,X], Y) - g_\theta (X , [T,Y])
= 2 g_\theta (\tau X , Y) \] (note that one makes use of the fact
that $\tau : H(\mathcal{F}) \to H(\mathcal{F})$ is self-adjoint,
i.e. $g_\theta (\tau X , Y) = g_\theta (X , \tau Y)$, $X,Y \in
H(\mathcal{F})$) hence
\[ g(\nabla^g_X Y , T) = - g(\tau X , Y) - \frac{n+1}{\varphi} \left(
\frac{1}{\varphi} - r \right) \, (d \theta )(X,Y) \] or
\begin{equation}
g(\nabla^g_X Y , T) = - g(\tau X , Y) - \left( \frac{1}{\varphi} -
r \right) g(X , \phi Y), \label{e:b10}
\end{equation}
for any $X, Y \in H(\mathcal{F})$. Exploiting again $\nabla^g g =
0$ (and $g([X,Y], N) = 0$) we get
\[ 2 g(\nabla^g_X Y , N) = - (\mathcal{L}_N g)(X,Y), \;\; X, Y \in
H(\mathcal{F}). \] Hence (by (\ref{e:A.23}) in Lemma \ref{l:A.2})
\[ 2 g(\nabla^g_X Y , N) = - \frac{n+1}{\varphi^2} \, N(\varphi )
\, g_\theta (X,Y) + \frac{n+1}{\varphi} \, (\mathcal{L}_N g_\theta
)(X,Y) = \]
\[ = 2 \left( \frac{1}{\varphi} - r \right) g(X,Y) +
\frac{2(n+1)}{\varphi} \, (d \theta )(X , \tau Y) \] that is
\begin{equation}
g(\nabla^g_X Y , N) = \left( \frac{1}{\varphi} - r \right) g(X,Y)
+ g(X , \phi \, \tau \, Y), \label{e:b11}
\end{equation}
for any $X, Y \in H(\mathcal{F})$. Note that (\ref{e:b11}) may be
also derived from (\ref{e:b10}) by using the fact that $g$ is a
K\"ahler metric. Indeed
\[ g(\nabla^g_X Y , N) = g(J \nabla^g_X Y , J N) = g(\nabla^g_X J Y
, T) = \] \[ = - g(\tau X , \phi Y) - \left( \frac{1}{\varphi} - r
\right) g(X , \phi^2 Y), \] etc. For further use, let us also
retain that
\begin{equation}
(\mathcal{L}_N g)(X,Y) = - 2 \left( \frac{1}{\varphi} - r \right)
g(X,Y) - 2 g(X , \phi \tau Y), \label{e:b12}
\end{equation} for any $X,Y \in H(\mathcal{F})$.
At this point, the identities (\ref{e:b8}) and
(\ref{e:b10})-(\ref{e:b11}) lead to
\begin{equation} \nabla^g_X Y = \nabla_X Y +
\label{e:b13}
\end{equation}
\[ + \left\{ \frac{\varphi}{1  - \varphi r} \, g_\theta (\tau X ,
Y) + g_\theta (X , \phi Y) \right\} T - \] \[ - \left\{ g_\theta
(X,Y) + \frac{\varphi}{1  - \varphi r} \, g_\theta (X , \phi \,
\tau \, Y) \right\} N, \] for any $X,Y \in H(\mathcal{F})$. To
compute $\nabla^g_X T$ we use (\ref{e:b13}) and
\[ g(\nabla^g_X T , Y) = - g(T , \nabla^g_X Y ) \]
so that
\begin{equation}
g(\nabla_X^g T , Y) = g(\tau X , Y) + \left( \frac{1}{\varphi} - r
\right) g(X , \phi Y). \label{e:b14}
\end{equation}
The component along $T$ is $\frac{1}{2} X( \| T \|^2 )$ hence
\begin{equation}
g(\nabla^g_X T , T) = - \frac{n+1}{2 \varphi} \, X(r).
\label{e:b15}
\end{equation} Moreover (by (\ref{e:A.5}) in Appendix A)
\[ 2 g(\nabla^g_X T , N) = g(X , [N,T]) = - g (X , \phi \nabla^H r
) \] that is
\begin{equation}
g(\nabla^g_X T , N) = - \frac{n+1}{2\varphi} (\phi X)(r).
\label{e:b16}
\end{equation}
Summing up (by (\ref{e:b14})-(\ref{e:b16}))
\begin{equation}
\nabla^g_X T = \tau X - \left( \frac{1}{\varphi} - r \right) \phi
X - \frac{\varphi}{2(1 - r \varphi )} \left\{ X(r) T + (\phi X)(r)
N \right\} , \label{e:b17}
\end{equation}
for any $X \in H(\mathcal{F})$. Again by (\ref{e:b13}) and
\[ g(\nabla^g_X N , Y) = - g(N , \nabla_X^g Y) \]
we get
\begin{equation}
g(\nabla^g_X N , Y) = - \left( \frac{1}{\varphi} - r \right) \,
g(X,Y) - g(X , \phi\, \tau\, Y). \label{e:b18}
\end{equation}
Next (by (\ref{e:A.5}))
\begin{equation}
g(\nabla^g_X N , T) = \frac{n+1}{2 \varphi} \, (\phi X)(r).
\label{e:b19}
\end{equation}
Finally, the component along $N$ is $\frac{1}{2} X(\| N\|^2 )$
hence
\begin{equation}
g(\nabla^g_X N , N) = - \frac{n+1}{2 \varphi} X(r).
\label{e:b20}
\end{equation}
Summing up (by (\ref{e:b18})-(\ref{e:b20}))
\begin{equation}
\nabla^g_X N = - \left( \frac{1}{\varphi} - r \right) X + \tau \,
\phi\, X +  \label{e:b21}
\end{equation}
\[ + \frac{\varphi}{2(1 - r \varphi )} \{ (\phi X)(r) \, T -
X(r) \, N \} , \] for any $X \in H(\mathcal{F})$. We wish to
compute $\nabla^g_T X$. To this end (by $\nabla g_\theta = 0$)
\[ 2 g_\theta (\nabla_T X , Y) = T(g_\theta (X,Y)) + g_\theta
([T,X] , Y) + g_\theta ([Y,T], X) + \]
\[ + g_\theta (T , [Y,X]) + g_\theta (\tau X, Y) - g_\theta (\tau
Y , X) - 2 (d \theta )(X,Y) \] yielding (upon multiplication by
$-(n+1)/\varphi$)
\[ T(g(X,Y)) + g([T,X] , Y) + g([Y,T], X) = 2 g(\nabla_T X , Y).
\] Therefore (by $\nabla^g g = 0$)
\[ 2 g(\nabla^g_T X , Y) = T(g(X,Y)) + g([T,X],Y) + g([Y,T],X) +
\]
\[ + g(T , [Y,X]) = 2 g(\nabla_T X , Y) - \theta ([X,Y]) \| T \|^2 \]
or
\begin{equation}
g(\nabla^g_T X , Y) = g(\nabla_T X , Y) + \left( \frac{1}{\varphi}
- r \right) \, g(X , \phi Y). \label{e:b22}
\end{equation}
Similar to the above
\begin{equation}
g(\nabla^g_T X , T) = - \frac{n+1}{2 \varphi} \, X(r),
\label{e:b23}
\end{equation}
\begin{equation}
g(\nabla_T^g X , N) = - \frac{n+1}{2 \varphi} \, (\phi X)(r).
\label{e:b24}
\end{equation} Collecting the information in
(\ref{e:b22})-(\ref{e:b24}), we have proved
\begin{equation}
\nabla^g_T X = \nabla_T X - \left( \frac{1}{\varphi} - r \right)
\phi X - \frac{\varphi}{2(1 - r \varphi )} \{ X(r) T + (\phi X)(r)
N \} ,
\label{e:b25}
\end{equation}
for any $X \in H(\mathcal{F})$. Let us compute $\nabla^g_N X$. We
have
\[ 2 g(\nabla^g_N X , Y) = N(g(X,Y)) + g([N,X], Y) + g([Y,N], X) =
\]
\[ = 2 g([N,X] , Y) + (\mathcal{L}_N g)(X,Y). \]
Using (\ref{e:b12}) and
\[ [N,X] = \nabla_N X - r X - \tau (\phi X) \]
(cf. Appendix A) one shows that
\begin{equation}
g(\nabla^g_N X , Y) = g(\nabla_N X , Y) - \frac{1}{\varphi} \;
g(X,Y).
\label{e:b26}
\end{equation}
Calculations similar to the above also furnish
\begin{equation}
g(\nabla^g_N X , T) = \frac{n+1}{2\varphi} \; (\phi X)(r),
\label{e:b27}
\end{equation}
\begin{equation}
g(\nabla^g_N X , N) = - \frac{n+1}{2 \varphi} \; X(r).
\label{e:b28}
\end{equation}
Using (\ref{e:b26})-(\ref{e:b28}) we may now conclude that
\begin{equation}
\nabla^g_N X = \nabla_N X - \frac{1}{\varphi} \, X +
\frac{\varphi}{2(1 - r \varphi )} \{ (\phi X)(r) T - X(r) N \} ,
\label{e:b29}
\end{equation}
for any $X \in H(\mathcal{F})$. Moreover (omitting the details)
\[ g(\nabla^g_N T , X) = - \frac{n+1}{2 \varphi} \, (\phi X)(r),
\]
\[ g(\nabla^g_N T , T) = - \frac{n+1}{2 \varphi} \left\{ N(r) +
\frac{4}{\varphi^2} - \frac{2 r}{\varphi} \right\} , \]
\[ g(\nabla^g_N T , N) = - \frac{n+1}{2 \varphi} \; T(r), \]so
that
\begin{equation}
\nabla^g_N T = - \frac{1}{2} \, \phi \, \nabla^H r -
 \label{e:b30}
\end{equation}
\[ - \frac{\varphi}{2(1 - r \varphi )} \left\{ \left( N(r) +
\frac{4}{\varphi^2} - \frac{2 r}{\varphi} \right) T + T(r) N
\right\} . \] Similarly we find
\begin{equation}
\nabla^g_T N = \frac{1}{2} \, \phi \nabla^H r -
\label{e:b31}\end{equation}
\[ - \frac{\varphi}{2(1 - r \varphi )} \left\{ \left( N(r) +
\frac{4}{\varphi^2} - \frac{6 r}{\varphi} + 4 r^2 \right) T + T(r)
N \right\} , \]
\begin{equation}
\nabla^g_T T = - \frac{1}{2} \; \nabla^H r - \label{e:b32}
\end{equation}
\[ - \frac{\varphi}{2(1 - r \varphi )} \left\{ T(r) T - \left(
N(r) + \frac{4}{\varphi^2} - \frac{6 r}{\varphi} + 4 r^2 \right) N
\right\} , \]
\begin{equation}
\nabla^g_N N = - \frac{1}{2} \; \nabla^H r + \label{e:b33}
\end{equation}
\[ + \frac{\varphi}{2(1 - r \varphi )} \left\{ T(r) T - \left(
N(r) + \frac{4}{\varphi^2} - \frac{2 r}{\varphi} \right) N
\right\} . \] Let us consider a holomorphic vector bundle $\pi : F
\to U$, carrying the Hermitian metric $h$, and set $E_\delta =
\pi^{-1} (M_\delta )$ (the portion of $F$ over a leaf of
$\mathcal{F}$). A connection $D \in \mathcal{C}(F , h)$ induces a
connection $D^\delta \in \mathcal{C}(E_\delta , h_\delta )$ (where
$h_{\delta , z} = h_z$, $z \in M_\delta$). $D^\delta$ is most
easily described with respect to a local trivialization $\Phi :
\pi^{-1}(O) \to O \times \mathbb{C}^m$ of $F$, for some open
subset $O \subseteq U$. Let us set $\sigma_i (z) = \Phi^{-1} (z,
e_i )$, $z \in O$, $1 \leq i \leq m$, where $\{ e_1 , \cdots , e_m
\}$ is the canonical linear basis in $\mathbb{C}^m$. If $\left.
u_i \equiv \sigma_i \right|_{O \cap M_\delta}$ then $D^\delta$ is
given by
\[ (D^\delta_X u)_z = X(f^i )_z u_i (z) + f^i (z) (D_{(d\; j_\delta
) X} \sigma_i )_z \, , \;\; z \in O \cap M_\delta , \] for any
section $u = f^i u_i$, $f^i \in C^\infty (O \cap M_\delta )$, and
any $X \in \mathcal{X}(M_\delta )$. It is easily shown that the
definition of $(D^\delta_X u)_z$ doesn't depend upon the local
trivialization chart $\Phi$ at $z$ (i.e. if $g = [g_{ij}] : O \cap
O^\prime \to {\rm GL}(m, \mathbb{C})$, $g(z) = \Phi^\prime_z \circ
\Phi_z^{-1}$, are the transition functions of $F$ then
$(D^\delta_X u)_z$ is invariant under the transformation $\sigma_j
(z) = g^i_j (z) \sigma^\prime_i (z)$). Let $R^D \in \Omega^2 ({\rm
Ad}(F))$ and $\omega^i_j$ be the curvature tensor field and
connection $1$-forms of $D$ ($D \sigma_j = \omega^i_j \otimes
\sigma_i$), so that $R^D \sigma_j = 2(d \omega^i_j - \omega^i_k
\wedge \omega^k_j ) \otimes \sigma_i$. Also, let $R^\delta \in
\Omega^2 ({\rm Ad}(E_\delta ))$ and $(\omega_\delta )^i_j$ be the
curvature tensor field and the connection $1$-forms of $D^\delta$,
respectively. Then $(\omega_\delta )^i_j = j_\delta^* \omega^i_j$
yields
\begin{equation}
R^\delta u_i = (j^*_\delta R^D ) \sigma_i \, , \;\; 1 \leq i \leq
m. \label{e:b34}
\end{equation}
Let $\{ W_\alpha \}$ be a local orthonormal ($g_\theta (W_\alpha ,
W_{\overline{\beta}}) = \delta_{\alpha\beta}$) frame of
$T_{1,0}(\mathcal{F})$ and set
\[ E_\alpha \equiv   \sqrt{-\frac{\varphi}{n+1}} \; W_\alpha \, ,
\; \; 1 \leq \alpha \leq n-1 , \;\; E_n \equiv
\sqrt{\frac{2f\varphi}{n+1}} \; \xi \, ,  \] where $f \equiv
\varphi /(1 - r \varphi )$. Then, given a connection $D$ in $F \to
U$, for any $X \in H(M_\delta )$
\[ (\delta^D R^D ) X = - \sum_{a=1}^n \{ (D_{E_a} R^D
)(E_{\overline{a}} , X) + (D_{E_{\overline{a}}} R^D )(E_a , X) \}
= \] \[ = \frac{\varphi}{n+1} \sum_{\alpha =1}^{n-1} \{
(D_{W_\alpha} R^D )(W_{\overline{\alpha}} , X) +
(D_{W_{\overline{\alpha}}} R^D )(W_\alpha , X) \} - \]
\[ - \frac{2f\varphi}{n+1} \{ (D_{\xi} R^D )(\overline{\xi} , X) +
(D_{\overline{\xi}} R^D )(\xi , X) \}  \] and
\[ \sum_{\alpha =1}^{n-1} (D_{W_\alpha} R^D )(W_{\overline{\alpha}} , X)
\sigma_j  =
\] \[ = \sum_\alpha \{
D_{W_\alpha} \left( R^D (W_{\overline{\alpha}} , X)
\sigma_j \right ) - R^D (W_{\overline{\alpha}} , X) D_{W_\alpha}
\sigma_j - \]
\[ - R^D (\nabla^g_{W_\alpha} W_{\overline{\alpha}} ,  X) \sigma_j -
R^D (W_{\overline{\alpha}} , \nabla^g_{W_\alpha} X ) \sigma_j \} =
\;\; (by \; (\ref{e:b13})) \] \[ = \sum_\alpha \left\{
(D^\delta_{W_\alpha} R^\delta )(W_{\overline{\alpha}} , X) u_j -
\right. \]
\[ - [ f \; g_\theta (\tau \, W_\alpha , W_{\overline{\alpha}}) +
g_\theta (W_\alpha , \phi \, W_{\overline{\alpha}}) ] R^D (T , X)
\sigma_j + \] \[ + [ g_\theta (W_\alpha , W_{\overline{\alpha}}) +
f \; g_\theta (W_\alpha , \phi \, \tau\,  W_{\overline{\alpha}}) ]
R^D (N, X) \sigma_j - \] \[ - [ f \; g_\theta (\tau \, W_\alpha ,
X) + g_\theta (W_\alpha , \phi X) ] R^D (W_{\overline{\alpha}} ,
T) \sigma_j + \] \[ + \left. [g_\theta (W_\alpha  , X) + f \;
g_\theta (W_\alpha , \phi \, \tau \, X) ]  R^D
(W_{\overline{\alpha}} , N) \sigma_j \right\} .
\]
Therefore (by the purity axiom (\ref{e:A.8}))
\[ \sum_\alpha (D_{W_\alpha} R^D )(W_{\overline{\alpha}} , X)  =
\sum_\alpha (D^\delta_{W_\alpha} R^\delta )(W_{\overline{\alpha}}
, X) u_j - \]
\[ + i(n-1) R^D (T , X) \sigma_j + (n-1) R^D (N, X) \sigma_j - \]
\[ - f \; R^D (\pi_{0,1} \, \tau \, X , T) - R^D (\pi_{0,1} \, \phi \, X , T)
+ \] \[ + R^D (\pi_{0,1} \, X , N) + f \; R^D (\pi_{0,1} \, \phi
\, \tau \, X , N).  \] We obtain
\begin{equation} \sum_\alpha \{
(D_{W_\alpha} R^D )(W_{\overline{\alpha}} , X) +
(D_{W_{\overline{\alpha}}} R^D )(W_\alpha , X) \} = -
(\delta_b^{D^\delta} R^\delta )_X \, u_j + \label{e:b35}
\end{equation}
\[ + \{ R^D (N , (2n-3) X - f \, \phi \, \tau \, X ) + R^D (T ,
\phi \, X + f \, \tau \, X)\} \sigma_j    \] (cf. section 5 for
the definition of the operator $\delta_b^{D^\delta}$). Moreover
\[ \{ (D_{\xi} R^D )(\overline{\xi} , X) +
(D_{\overline{\xi}} R^D )(\xi , X) \} \sigma_j = \]
\[ = \frac{1}{2} \{ D_N \left( R^D (N , X) \sigma_j \right) + D_T
\left( R^D (T , X) \sigma_j \right) - \] \[ - R^D (N,X) D_N
\sigma_j - R^D (T , X) D_T \sigma_j - \] \[ - R^D (\nabla^g_N N ,
X) \sigma_j - R^D (\nabla^g_T T , X) \sigma_j - \]
\[ - R^D (N , \nabla^g_N X ) \sigma_j - R^D (T , \nabla^g_T X)
\sigma_j \}  . \] Substitution from (\ref{e:b25}), (\ref{e:b29})
and (\ref{e:b32})-(\ref{e:b33}) gives
\[ - R^D (\nabla^g_N N ,
X) - R^D (\nabla^g_T T , X) - R^D (N , \nabla^g_N X ) - R^D (T ,
\nabla^g_T X) =  \]
\[ = R^D (\nabla^H r , X)  - R^D (T , \nabla_T X) - R^D (N ,
\nabla_N X)  + \] \[ + \frac{1}{f} \; R^D (T , \phi \, X) + f \;
(\phi X)(r) \, R^D (T , N) + \left( \frac{1}{\varphi} + 2 r
\right) R^D (N , X). \] We conclude that
\begin{equation}
\{ (D_{\xi} R^D )(\overline{\xi} , X) + (D_{\overline{\xi}} R^D
)(\xi , X) \} \sigma_j = \label{e:b36}
\end{equation}
\[ = \frac{1}{2} \{ (D_N \; i_N R^D ) X + (D_T \; i_T R^D )X + R^D
(\nabla^H r , X) + \]
\[ + \frac{1}{f} \; R^D (T , \phi \, X) + f \;
(\phi X)(r) \, R^D (T , N) + \left( \frac{1}{\varphi} + 2 r
\right) R^D (N , X) \] (the covariant derivatives in the right
hand member of (\ref{e:b36}) are defined with respect to $D$ and
$\nabla$). Finally (by (\ref{e:b35})-(\ref{e:b36}))
\[ (\delta^D R^D )_X \,  \sigma_j = \frac{\varphi}{n+1} \left\{
- (\delta_b^{D^\delta} R^\delta )_X \, u_j \right. + \]
\[ \left. +  [ R^D (N , (2n-3) X - f \, \phi \, \tau \, X ) + R^D (T ,
\phi \, X + f \, \tau \, X)] \sigma_j \right\} - \] \[ - \frac{f
\, \varphi}{n+1} \left\{ (D_N \; i_N R^D ) X + (D_T \; i_T R^D )X
+ R^D (\nabla^H r , X) + \right. \] \[ + \left. R^D (T \, , \,
\frac{1}{f} \, \phi \, X + f \, (\phi\, X)(r) \, N) + \left(
\frac{1}{\varphi} + 2 r \right) R^D (N , X) \right\} \, \sigma_j
\, . \] Assume that $D$ is a Yang-Mills field on $(\Omega  , g)$,
i.e. $\delta^D R^D = 0$ in $\Omega$. Then, for $\varphi \to 0$ (as
$r$ and $\nabla^H r$ stay finite near $\partial \Omega$, cf.
\cite{kn:GrLe}, p. 164)
\[ (\delta_b^{D_b} R^{D_b} )_X u_j = 2(n-2) R^D (N, X)\sigma_j \]
where $D_b \equiv D^0$ is the boundary values of $D$. Therefore,
if $i_T R^{D_b} = 0$ then (cf. (\ref{e:pymb}) in section 5) $D_b$
is a pseudo Yang-Mills field on $\partial \Omega$ if and only if
$i_N R^D = 0$ on $H(\partial \Omega )$. Theorem \ref{t:3} is
proved. With the same techniques we may show that
\begin{corollary} Let $D \in \mathcal{C}(F, h)$ be a Yang-Mills field
on $(\Omega , g)$ such that $i_N R^D = 0$. Then the boundary
values $D_b$ of $D$ satisfy $\Lambda_\theta R^{D_b} = 0$.
\label{c:1}
\end{corollary}
\noindent Corollary \ref{c:1} shows that the axiom (\ref{e:axS})
(with $S = 0$) in the description of the Tanaka connection, as
well as (\ref{e:02}) in Theorem \ref{t:1}, are rather natural
occurrences. The proof is
\[ 0 = (\delta^D R^D )_T \sigma_j = \]
\[ = \frac{\varphi}{n+1} \sum_\alpha \{
(D_{W_\alpha} R^D )(W_{\overline{\alpha}} , T) +
(D_{W_{\overline{\alpha}}} R^D )(W_\alpha , T) \} \sigma_j  - \]
\[ - \frac{2 f \varphi}{n+1} \{ D_N (R^D (N,T) \sigma_j ) - R^D
(N,T) D_N \sigma_j - \]
\[ - R^D (\nabla^g_N N , T)  - R^D (N , \nabla^g_N T ) \sigma_j \}
\]
or (by (\ref{e:b13}), (\ref{e:b17}), (\ref{e:b30}) and
(\ref{e:b33}))
\[ 0 = \varphi \{ (\delta^{D^\delta} R^\delta )_T \; u_j + 2(n-1) R^D
(T , N) \sigma_j \} + \]
\[ + \frac{1}{2} \; \varphi \, \{ f \, [ R^D (T , \nabla^H r )  + R^D
(N , \phi \, \nabla^H r )]  + trace \, \pi_H R^D (\cdot \, , \;
\tau \, \cdot ) \} \sigma_j - \]
\[  + (n+1) f \{ 2 (D_N \; i_N
R^D )T + R^D (N , \phi \nabla^H r) - R^D (T , \nabla^H r) \}
\sigma_j - \]
\[ - \frac{2\varphi}{f} \; \Lambda_\theta R^D \sigma_j
- (n+1) f^2 \; \left\{ N(r) + \frac{4}{\varphi^2} - \frac{2
r}{\varphi} \right\} \, R^D (T , N)\sigma_j \, . \] When $\varphi
\to 0$ one observes $\varphi /f \to 1$ and $f^2 /\varphi^2 \to 1$
hence
\[ \left( \Lambda_\theta R^{D_b} \right) u_j = - 2 (n+1) R^D (T , N) \sigma_j \,
. \] Q.e.d.

\section{Yang-Mills fields and the Fefferman metric} We wish to
relate $\mathcal{PYM}$ to the Yang-Mills functional on $C(M)$.
Given a contact form $\theta$ on $M$ such that the Levi form
$L_\theta$ is positive definite, let $F_\theta$ be the
corresponding {\em Fefferman metric} (a Lorentz metric on $C(M)$).
We recall (cf. \cite{kn:Lee}) that
\begin{equation}
F_\theta = \pi^* \tilde{G}_\theta + 2 (\pi^* \theta ) \odot \sigma
, \label{e:4}
\end{equation}
\begin{equation}
\sigma = \frac{1}{n+2} \{ d \gamma + \pi^* ( i
\omega^\alpha_\alpha - \frac{i}{2} g^{\alpha\overline{\beta}}  d
g_{\alpha\overline{\beta}} - \frac{\rho}{4(n+1)} \, \theta ) \} .
\label{e:5}
\end{equation}
Here $\omega^\alpha_\beta$ are the connection $1$-forms of the
Tanaka-Webster connection of $(M , \theta )$, i.e. $\nabla T_\beta
= \omega^\alpha_\beta \otimes T_\alpha$, and
$g_{\alpha\overline{\beta}} = L_\theta (T_\alpha ,
T_{\overline{\beta}})$. Moreover $\rho =
g^{\alpha\overline{\beta}} R_{\alpha\overline{\beta}}$ is the {\em
pseudohermitian scalar curvature} (cf. e.g. \cite{kn:Dra}, p.
229). The $(0,2)$-tensor field $\tilde{G}_\theta$ is got by
extending the Levi form $G_\theta$ to the whole of $T(M)$.
Precisely, one requests that $\tilde{G}_\theta = G_\theta$ on
$H(M) \otimes H(M)$, while $\tilde{G}_\theta (X, T) = 0$, for any
$X \in T(M)$ (obviously $\tilde{G}_\theta$ is degenerate). Note
that when $M$ is compact $C(M)$ is compact, as well. It is
noteworthy that $\sigma$ (given by (\ref{e:5})) is a connection
$1$-form in $S^1 \to C(M) \to M$. Let $T^\uparrow$ be the
horizontal lift (with respect to $\sigma$) of the characteristic
direction of $d \theta$ and $S$ the tangent to the $S^1$-action.
Then $T^\uparrow - S$ is timelike, hence $(C(M), F_\theta )$ is
time oriented by $T^\uparrow - S$, i.e. $(C(M) , F_\theta )$ is a
space-time (see \cite{kn:BeEh}, p. 17). However, as $M$ is compact
$(C(M), F_\theta )$ is not chronological (cf. Proposition 2.6 in
\cite{kn:BeEh}, p. 23). \par Let $S \in \mathcal{X}(C(M))$ be the
tangent to the $S^1$-action (locally $S = ((n+2)/2) \partial
/\partial \gamma$). Then (by (\ref{e:4})) $F_\theta (S , S) = 0$.
Next (by Lemma \ref{l:3}) $\nabla^{C(M)}_S S = 0$, i.e. the
integral curves of $S$ are (null) geodesics of $(C(M), F_\theta
)$. Also $\mathcal{L}_S F_\theta = 0$, hence (cf. (28) in
\cite{kn:NuTr}, p. 185) $S$ generates a {\em shear-free
congruence} of null geodesics. The congruence is {\em symmetric}
if there is a vector field $X \in \mathcal{X}(C(M))$ such that
\[ \mathcal{L}_X (\pi^* \theta ) = t \; \pi^* \theta , \;\;
\mathcal{L}_X (\pi^* \theta^\alpha ) = w^\alpha_\beta \, \pi^*
\theta^\beta + \ell^\alpha \, \pi^* \theta , \] where $t$ is a
real function and $w^\alpha_\beta \, , \; \ell^\alpha$ are complex
functions on $C(M)$. We say that $X$ is a {\em symmetry} of the
congruence. We may look for $\mathcal{L}_X S$ in the form
\[ \mathcal{L}_X S = a^\alpha T_\alpha^\uparrow +
a^{\overline{\alpha}} T_{\overline{\alpha}}^\uparrow + b
T^\uparrow + f S. \] As $X$ is a symmetry $a^\alpha = 0$, $b = 0$
and $f = - 2 (\mathcal{L}_X \sigma )S$. Therefore
\begin{equation}
\mathcal{L}_X S = f \, S. \label{e:n1}
\end{equation} Also one may easily check (by using the local frame
$\{ T^\uparrow_\alpha , T^\uparrow_{\overline{\alpha}} ,
T^\uparrow , S \}$ of $T(C(M)) \otimes \mathbb{C}$) that
\begin{equation}
\mathcal{L}_S \, (\pi^* \theta ) = 0, \;\; \mathcal{L}_S \, (\pi^*
\theta^\alpha ) = 0. \label{e:n2}
\end{equation}
Using (\ref{e:n1})-(\ref{e:n2}) and $\mathcal{L}_Y \mathcal{L}_Z
\omega = \mathcal{L}_Z \mathcal{L}_Y \omega + \mathcal{L}_{[Y,Z]}
\omega$ (for any $Y,Z \in \mathcal{X}(C(M))$, $\omega \in \Omega^1
(C(M))$) we obtain
\begin{equation}
S(t) = 0, \;\; S(w^\alpha_\beta ) = 0, \;\; S(\ell^\alpha ) = 0.
\label{e:n3}
\end{equation}
For instance
\[ S(t) \pi^* \theta = \mathcal{L}_{[S,X]} \, \pi^* \theta = - f \; \mathcal{L}_S
\, \pi^* \theta - (\pi^* \theta )(S) \, d f = 0. \] Our
considerations draw inspiration from the calculations in
\cite{kn:Nur} (which are both purely local and confined to the
$3$-dimensional case ($n=1$)). For this reason some of the results
(e.g. Propositions \ref{p:n2} and \ref{p:n1}) are attributed to
\cite{kn:Nur} (the proofs are however new). (\ref{e:n3}) implies
that $t, \, w^\alpha_\beta , \, \ell^\alpha$ are vertical lifts of
functions on $M$. A vector field of the form $\rho \, S$, for some
function $\rho \neq 0$, is a {\em trivial} symmetry of the
congruence.
\begin{proposition}
$(${\rm P. Nurowski, \cite{kn:Nur}}$)$
\par\noindent
Each nontrivial symmetry of the shear-free congruence $($of null
geodesics on $C(M))$ projects on a unique symmetry of the CR
structure on $M$. \label{p:n1}
\end{proposition} \noindent
Indeed, if $X$ is a symmetry of $S$ then $X - 2 \sigma (X) S \in
{\rm Ker}(\sigma )$, hence there is a unique vector field
$\tilde{X} \in \mathcal{X}(M)$ such that
\[ \tilde{X}^\uparrow = X - 2 \sigma (X) \, S. \]
Then
\[ \pi^* (t \, \theta ) = \mathcal{L}_X (\pi^* \theta ) =
\mathcal{L}_{\tilde{X}^\uparrow} (\pi^* \theta ) + 2
\mathcal{L}_{\sigma (X) S} (\pi^* \theta ) =
\mathcal{L}_{\tilde{X}^\uparrow} ( \pi^* \theta ). \]
Consequently, for any $Z \in \mathcal{X}(M)$
\[ t \, \theta (Z) = \left( \pi^* (t \, \theta ) \right)
Z^\uparrow = (\mathcal{L}_{\tilde{X}^\uparrow} (\pi^* \theta ))
Z^\uparrow = \tilde{X} (\theta (Z)) - (\pi^* \theta )
[\tilde{X}^\uparrow \, , \, Z^\uparrow ]  \] hence, as $[\tilde{X}
, Z]^\uparrow$ is the ${\rm Ker}(\sigma )$-component of
$[\tilde{X}^\uparrow \, , \, Z^\uparrow ]$ (with respect to the
decomposition $T(C(M)) = {\rm Ker}(\sigma ) \oplus \mathbb{R} S$),
we obtain $\mathcal{L}_{\tilde{X}} \theta = t \, \theta$. It may
be shown in a similar manner that $\mathcal{L}_{\tilde{X}}
\theta^\alpha = w^\alpha_\beta \theta^\beta + \ell^\alpha \theta$,
i.e. $\tilde{X}$ is a symmetry of the CR structure. Q.e.d.
\par Let $E \to M$ be a complex vector bundle and $\hat{E} = \pi^*
E \to C(M)$ the pullback of $E$ via $\pi$. The {\em natural lift}
$\hat{u} : \pi^{-1}(U) \to \hat{E}$ of a section $u : U \to E$ is
given by $\hat{u}(z) = (x , u(\pi (z)))$, $z \in \pi^{-1}(U)$. If
$E$ carries a Hermitian metric $h$ then so does $\hat{E}$. Indeed
we may set $\hat{h}(\hat{e}_i , \hat{e}_j ) = h_{i\overline{j}}
\circ \pi$, where $h_{i\overline{j}} = h(e_i , e_j )$ and $\{ e_1
, \cdots , e_m \}$ is a (local) frame in $E$ on $U$. There is a
natural inner product $\langle \; , \; \rangle$ on $\Omega^2 ({\rm
Ad}(\hat{E}))$ induced by the inner product on scalar $2$-forms
\[ F_\theta^* (\alpha , \beta ) \; d \, {\rm vol}
(F_\theta ) = \alpha \wedge * \beta , \] $\alpha , \beta \in
\Gamma^\infty (\Lambda^2 T^* (C(M)))$, and by the Killing-Cartan
form of ${\bf u}(m)$, $m = {\rm rank}_{\mathbb{C}} E$,
respectively. Here $*$ is the Hodge operator associated with the
Fefferman metric $F_\theta$. Precisely, if $S, T \in \Omega^0
({\rm Ad}(\hat{E}))$ then
\[ \langle \alpha \otimes S , \beta \otimes T \rangle =
F^*_\theta (\alpha , \beta ) [S^i_j ] \cdot [T^i_j ], \] where $S
\hat{e}_j = S^i_j \hat{e}_i$, $T \hat{e}_j = T^i_j \hat{e}_i$ with
respect to a (local) orthonormal ($h(e_i , e_j ) = \delta_{ij}$)
frame $\{ e_j \}$ in $E$, and $A \cdot B = - {\rm trace}(A B)$,
$A, B \in {\bf u}(m)$. The Yang-Mills functional is given by
\[ \widehat{\mathcal{YM}} (\mathbb{D}) = \frac{1}{2} \int_{C(M)}
\langle R^{\mathbb{D}} , R^{\mathbb{D}} \rangle \; d \, {\rm
vol}(F_\theta ), \;\; \mathbb{D} \in \mathcal{C}(\hat{E} ,
\hat{h}). \]  Any $D \in \mathcal{C}(E , h)$ induces a connection
$\hat{D} = \pi^* D \in \mathcal{C}(\hat{E} , \hat{h})$ which is
described (in local coordinates) as follows. Let $(U , x^A )$ be a
local coordinate system on $M$. Then $(\pi^{-1}(U), \, \hat{x}^A
:= x^A \circ \pi , \, \gamma )$ are local coordinates on $C(M)$.
We set by definition
\[ \hat{D}_{\partial /\partial \hat{x}^A} \hat{e}_j =
(\Gamma^i_{Aj} \circ \pi ) \hat{e}_i \, , \;\; \hat{D}_{\partial
/\partial \gamma} \hat{e}_j = 0, \] where $D_{\partial /\partial
x^A} e_j = \Gamma^i_{Aj} e_i$. Our conventions as to the range of
indices are $A,B,C, \cdots \in \{ 1, \cdots , 2n+1 \}$ and $i,j,k,
\cdots \in \{ 1, \cdots , m \}$. We consider the linear map
\[ \pi^* : \Gamma^\infty (U, \Lambda^k T^* (M) \otimes E) \to
\Gamma^\infty (\pi^{-1}(U), \Lambda^k T^* (C(M)) \otimes \hat{E})
\] given by
\[ \pi^* (\omega^j \otimes e_j ) = (\pi^* \omega^j ) \otimes
\hat{e}_j \, , \;\; \omega^j \in \Omega^k (U), \] (pullback and
natural lifting). As $\{ \hat{e}_j \}$ is a local frame in
$\hat{E} \to C(M)$ it suffices to specify $\hat{D}$ on natural
lifts of sections in $E \to M$. Then $\hat{D}$ admits the
following coordinate-free description
\[ \hat{D} \hat{u} = \pi^* (D u), \;\;\; u \in \Omega^0 (E). \]
Clearly, if $D h = 0$ then $\hat{D} \hat{h} = 0$. Let us consider
the functional $\mathcal{PYM} : \mathcal{C}(E , h) \to [0, +
\infty )$ given by
\[ \mathcal{PYM}(D) = \frac{1}{2} \int_M \| \pi_H R^D \|^2 \,
\theta \wedge (d \theta )^n . \] Here $\pi_H : \Omega^2 ({\rm
Ad}(E)) \to \Omega^2 ({\rm Ad}(E))/\mathcal{J}_\theta^2$ is the
projection described in section 2. Of course, when an admissible
coframe $\{ \theta^\alpha \}$ is fixed $\Omega^{\bf \cdot}({\rm
Ad}(E))/\mathcal{J}^{\bf \cdot}_\theta$ may be identified with the
subalgebra $\Omega_H^{\bf \cdot}({\rm Ad}(E)) = \{ \omega \in
\Omega^{\bf \cdot}({\rm Ad}(E)) : i_T \, \omega = 0 \}$.
Integration along the fibre in $\widehat{\mathcal{YM}}(\pi^* D)$,
$D \in \mathcal{C}(E,h)$, leads to (\ref{e:1}) in Theorem
\ref{t:1}. Indeed, let us set
\[ R^i_{ABj} e_i = (R^D e_j )(\partial /\partial x^A \, ,
\, \partial /\partial x^B). \] Then  \[ R^{\hat{D}} \hat{e}_j = [
(R^i_{ABj} \circ \pi ) d \hat{x}^A \wedge d \hat{x}^B ] \otimes
\hat{e}_i \] hence \begin{equation} \label{e:inpiu} R^{\hat{D}}
\hat{e}_j = \pi^* (R^D e_j ). \end{equation} Given $\Omega =
\Omega^j \otimes e_j \in \Omega^2 (E)$, $\Omega^j = \Omega^j_{AB}
d x^A \wedge d x^B$, we set
\begin{equation}
\langle \pi^* \Omega , \pi^* \Omega \rangle = F^*_\theta (\pi^*
\Omega^j , \pi^* \Omega^{\overline{k}} ) (h_{j\overline{k}} \circ
\pi ). \label{e:6}
\end{equation}
Of course $\pi^* \Omega^j = (\Omega_{AB}^j \circ \pi ) d \hat{x}^A
\wedge d \hat{x}^B$ and the main technical difficulty in
calculating (\ref{e:6}) is the need for $F^{AB} = F^*_\theta (d
\hat{x}^A , d \hat{x}^B )$, where $[F^{AB}] = [F_{AB}]^{-1}$ and
$F_{AB} := F_\theta (\partial /\partial \hat{x}^A , \partial
/\partial \hat{x}^B )$. Let
\[ F_\theta : \left[ \begin{array}{cc} F_{AB} & F_{A, 2n+2} \\
F_{2n+2,B} & F_{2n+2, 2n+2} \end{array} \right] \] be the
components of the Fefferman metric with respect to $(\hat{x}^A ,
\gamma )$. Let us set $\partial /\partial x^A = \lambda_A^B T_B$,
$\lambda^B_A \in C^\infty (U)$. Here one either adopts the
convention $A,B,C, \cdots \in \{ 0, 1, \cdots , n , \overline{1},
\cdots , \overline{n} \}$ (with $T_0 = T$) or relabels the vector
fields $\{ T , T_\alpha , T_{\overline{\alpha}} : 1 \leq \alpha
\leq n \}$. Then (by (\ref{e:4}))
\[ F_{AB} = \tilde{G}_\theta (\frac{\partial}{\partial x^A} ,
\frac{\partial}{\partial x^B}) + \theta (\frac{\partial}{\partial
x^A}) \sigma (\frac{\partial}{\partial \hat{x}^B}) + \theta
(\frac{\partial}{\partial x^B}) \sigma (\frac{\partial}{\partial
\hat{x}^A}) = \] \[ = g_{\alpha\overline{\beta}} (\lambda_A^\alpha
\lambda_B^{\overline{\beta}} + \lambda_B^\alpha
\lambda_A^{\overline{\beta}} ) + \lambda_A^0 \sigma_B +
\lambda_B^0 \sigma_A \] where $\sigma_A = \sigma (\partial /
\partial \hat{x}^A )$. A calculation based on (\ref{e:5}) shows
that
\[ \sigma_A = \frac{1}{n+2} \{ i \lambda^B_A (
\Gamma_{B\alpha}^\alpha - \frac{1}{2} g^{\alpha\overline{\beta}}
T_B (g_{\alpha\overline{\beta}} )) - \frac{\rho}{4(n+1)} \,
\lambda^0_A \} \circ \pi \] where $\Gamma^\beta_{B\alpha}$ are
(among) the coefficients of the Tanaka-Webster connection of $(M,
\theta)$ (i.e. $\nabla_{T_B} T_\alpha = \Gamma_{B\alpha}^\beta
T_\beta$). Moreover (by (\ref{e:4}))
\[ F_{A,2n+2} = 2 [(\pi^* \theta ) \odot \sigma
](\frac{\partial}{\partial \hat{x}^A} , \frac{\partial}{\partial
\gamma} ) = \frac{1}{n+2} \, \lambda^0_A \, , \]
\[ F_{2n+2,2n+2} = 0. \]
Next, using $F^{ab} F_{bc} = \delta^a_c$ (with $a,b,c, \cdots \in
\{ 1, \cdots , 2n+2 \}$) we find
\begin{equation} \begin{cases} F^{AB} F_{BC} +
\displaystyle{\frac{\lambda_C^0}{n+2} \, F^{A,2n+2}} = \delta^A_C
\cr F^{AB} \lambda^0_B = 0 \cr F^{2n+2,B} F_{BC} +
\displaystyle{\frac{\lambda^0_C}{n+2} \, F^{2n+2,2n+2}} = 0 \cr
F^{2n+2,B} \lambda^0_B = n+2. \cr
\end{cases} \label{e:7}
\end{equation}
Let us set
\[ {{P_{AB}}^i}_j \; e_i = (R^D e_j )(T_A , T_B ) \] so that
$R^i_{ABj} = \lambda_A^C \lambda_B^D {{P_{CD}}^i}_j$. In the
sequel, for the sake of simplicity, we do not distinguish
notationally between $f \in C^\infty (M)$ and its vertical lift $f
\circ \pi$. Then
\[ \langle R^{\hat{D}} , R^{\hat{D}} \rangle = h^{j\overline{k}}
\langle R^{\hat{D}} \hat{e}_j , R^{\hat{D}} \hat{e}_k \rangle = \]
\[ =  h^{j\overline{k}} h_{r\overline{s}} F^*_\theta (R^r_{ABj} d
\hat{x}^A \wedge d \hat{x}^B , R^s_{CDk} d \hat{x}^C \wedge
d\hat{x}^D ) = \]
\[ = \frac{1}{2} h^{j\overline{k}} h_{r\overline{s}} R^r_{ABj}
R^{\overline{s}}_{CD\overline{k}} (F^{AC} F^{BD} - F^{AD} F^{BC}),
\] where $R^{\overline{i}}_{AB\overline{j}} =
\overline{R^i_{ABj}}$. We obtain \[ \langle R^{\hat{D}} ,
R^{\hat{D}} \rangle = \] \begin{equation} = \frac{1}{2}
h^{j\overline{k}} h_{r\overline{s}} \lambda^E_A \lambda^F_B
\lambda^{\overline{G}}_C \lambda^{\overline{H}}_D {{P_{EF}}^r}_j
{{P_{\overline{G} \, \overline{H}}}^{\overline{s}}}_{\overline{k}}
(F^{AC} F^{BD} - F^{AD} F^{BC}), \label{e:8}\end{equation} where
$\lambda^{\overline{B}}_A = \overline{\lambda^B_A}$ and
${{P_{\overline{A} \,
\overline{B}}}^{\overline{i}}}_{\overline{j}} =
\overline{{{P_{AB}}^i}_j}$. Note that $\lambda^0_A$ is real valued
while $\lambda^{\overline{\alpha}}_A = \lambda_A^{\alpha + n}$. To
calculate $\lambda^E_A \lambda^F_B \lambda^{\overline{G}}_C
\lambda^{\overline{H}}_D (F^{AC} F^{BD} - F^{AD} F^{BC})$ we need
the identities
\begin{equation} F^{AB} \lambda_A^\alpha
\lambda^{\overline{\beta}}_B = g^{\alpha\overline{\beta}} ,
\label{e:9}\end{equation}
\begin{equation}
F^{AB} \lambda_A^\alpha \lambda_B^\beta = 0.
\label{e:10}\end{equation} The proof of (\ref{e:9})-(\ref{e:10})
follows from (\ref{e:7}). Indeed (\ref{e:7}) may be written
\[ F^{AB} g_{\alpha\overline{\beta}} (\lambda_B^\alpha
\lambda_C^{\overline{\beta}} + \lambda_C^\alpha
\lambda^{\overline{\beta}}_B ) + F^{AB} \lambda^0_C \sigma_B +
\frac{1}{n+2} \, F^{A, 2n+2} \lambda^0_C = \delta^A_C , \]
\[ F^{AB} \lambda^0_B = 0, \]
\[ F^{2n+2,B} g_{\alpha\overline{\beta}} (\lambda_B^\alpha
\lambda_C^{\overline{\beta}} + \lambda_C^\alpha
\lambda^{\overline{\beta}}_B ) + (n+2) \sigma_C + \] \[ + F^{2n+2,
B} \lambda^0_C \sigma_B + \frac{1}{n+2} \, F^{2n+2,2n+2}
\lambda^0_C = 0, \]
\[ F^{2n+2,B} \lambda^0_B = n+2. \]
If $\mu := \lambda^{-1}$ then (by the first of the previous four
identities)
\[ \mu^A_D = (\frac{1}{n+2} \, F^{A,2n+2} + F^{AB} \sigma_B )
\delta^0_D + F^{AB} g_{\alpha\overline{\beta}} (\lambda^\alpha_B
\delta^{\beta + n}_D + \lambda^{\overline{\beta}}_B
\delta^\alpha_D ) \] yielding
\begin{equation}
\begin{cases} \mu^A_0 = \displaystyle{\frac{1}{n+2} \, F^{A, 2n+2} + F^{AB}
\sigma_B} \cr \mu^A_\alpha = F^{AB} g_{\alpha\overline{\beta}}
\lambda^{\overline{\beta}}_B \cr \mu^A_{\beta + n} = F^{AB}
g_{\alpha\overline{\beta}} \lambda^\alpha_B . \cr
\end{cases}\label{e:11}
\end{equation}
The second and third of the identities (\ref{e:11}) lead to
(\ref{e:9}) and (\ref{e:10}), respectively. A calculation based on
(\ref{e:9})-(\ref{e:10}) shows that (\ref{e:8}) may be written \[
\langle R^{\hat{D}} , R^{\hat{D}} \rangle = \] \begin{equation}
\label{e:12} = P^{\overline{\alpha} \, \overline{\beta} \,
\overline{k} j} P_{\overline{\alpha} \, \overline{\beta} \,
\overline{k} j} + P^{\alpha \beta \overline{k} j} P_{\alpha \beta
\overline{k} j} + P^{\overline{\alpha} \beta \overline{k} j}
P_{\overline{\alpha} \beta \overline{k} j} + P^{\alpha
\overline{\beta} \, \overline{k} j} P_{\alpha \overline{\beta} \,
\overline{k} j} \end{equation} where $P_{AB\overline{k}j} =
h_{j\overline{s}} {P_{AB\overline{k}}}^{\overline{s}}$ and
$P^{AB\overline{k}j} = h^{r\overline{k}} {{P^{AB}}_r}^j$. Also
${{P^{\alpha\beta}}_r}^j = g^{\alpha\overline{\lambda}}
g^{\beta\overline{\mu}} {P_{\overline{\lambda} \, \overline{\mu}\,
r}}^j$, etc. As $R^D e_j = ({{P_{AB}}^i}_j \theta^A \wedge
\theta^B ) \otimes e_i$ it follows that
\[ \langle \pi_H R^D , \pi_H R^D \rangle = h^{j\overline{k}}
\langle (\pi_H R^D )e_j , (\pi_H R^D ) e_k \rangle = \]
\[ =  h^{j\overline{k}} h_{r\overline{s}} G^*_\theta
({{P_{\alpha\beta}}^r}_j \; \theta^\alpha \wedge \theta^\beta + 2
{{P_{\alpha\overline{\beta}}}^r}_j \; \theta^\alpha \wedge
\theta^{\overline{\beta}} + {{P_{\overline{\alpha} \,
\overline{\beta}}}^r}_j \; \theta^{\overline{\alpha}} \wedge
\theta^{\overline{\beta}} , \] \[ {{P_{\lambda\mu}}^s}_k \;
\theta^\lambda \wedge \theta^\mu + 2
{{P_{\lambda\overline{\mu}}}^s}_k \; \theta^\lambda \wedge
\theta^{\overline{\mu}} + {{P_{\overline{\lambda} \,
\overline{\mu}}}^s}_k \; \theta^{\overline{\lambda}} \wedge
\theta^{\overline{\mu}}) = \]
\[ = \frac{1}{2} \, {P_{\alpha\beta}}^{\overline{k}j}
P_{\overline{\lambda} \, \overline{\mu} \, \overline{k} j}
(g^{\alpha\overline{\lambda}} g^{\beta\overline{\mu}} -
g^{\alpha\overline{\mu}} g^{\beta\overline{\lambda}} ) + \]
\[ + 2 g^{\alpha\overline{\lambda}} g^{\overline{\beta}\mu}
{P_{\alpha\overline{\beta}}}^{\overline{k} j}
P_{\overline{\lambda} \mu \overline{k} j} + \]
\[ + \frac{1}{2}\, {P_{\overline{\alpha} \,
\overline{\beta}}}^{\overline{k} j} P_{\lambda \mu \overline{k} j}
(g^{\overline{\alpha}\lambda} g^{\overline{\beta}\mu} -
g^{\overline{\alpha} \mu} g^{\overline{\beta} \lambda}) = \]
\[ = P^{\overline{\lambda}\, \overline{\mu} \, \overline{k} j}
P_{\overline{\lambda}\, \overline{\mu} \, \overline{k} j} + 2
P^{\overline{\lambda} \mu \overline{k} j} P_{\overline{\lambda}
\mu \overline{k} j} + P^{\lambda\mu \overline{k} j} P_{\lambda\mu
\overline{k} j} \] hence (by (\ref{e:inpiu}) and (\ref{e:12}))
\[ \langle R^{\hat{D}} , R^{\hat{D}} \rangle = ( \| \pi_H R^D \|
\circ \pi )^2 . \] Finally we may integrate over $C(M)$ and use
the identity
\[ \int_{C(M)} (f \circ \pi ) d {\rm vol}(F_\theta ) = 2 \pi
\int_M f \, \theta \wedge (d \theta )^n \, , \;\; f \in C^\infty
(M). \] The identity\footnote{The symbol $\pi$ in the right hand
side denotes the irrational number $\pi \in \mathbb{R}$.}
(\ref{e:1}) in Theorem 1 is proved. Assume now that $\hat{D} =
\pi^* D$ is a Yang-Mills field on $C(M)$. Let $D^t = D + t
\varphi$, $\varphi \in \Omega^1 ({\rm Ad}(E))$, be a variation of
$D$. Then
\begin{equation} \pi^* D^t =  \hat{D} + t \pi^* \varphi .
\label{e:new}
\end{equation}
A word on the conventions in (\ref{e:new}). As seen earlier in
this section, there is a natural map $\pi^* : \Omega^1 ({\rm Ad}
(E)) \to \Omega^1 (\pi^* {\rm Ad}(E))$. Yet ${\rm Ad}(\pi^* E)
\approx \pi^* {\rm Ad} (E)$ (a vector bundle isomorphism) hence
$\pi^* \varphi$ is an ${\rm Ad}(\pi^* E)$-valued $1$-form on
$C(M)$. Then (by (\ref{e:1})) \[ 0 = \frac{d}{dt} \{
\widehat{\mathcal{YM}}(\hat{D} + t \pi^* \varphi ) \}_{t = 0} = \]
\[ =  \frac{d}{dt} \{ \widehat{\mathcal{YM}}(\pi^* D^t )\}_{t=0} = 2
\pi \frac{d}{dt} \{ \mathcal{PYM}(D^t )\}_{t=0} , \] i.e. $D$ is a
pseudo Yang-Mills field on $M$. The converse requires the first
variation formula for the functional $\mathcal{PYM}$ (as well as
the fact that the Yang-Mills equations on $C(M)$ project on $M$
via $\pi$ to give the Euler-Lagrange equations of the variational
principle $\delta \; \mathcal{PYM} = 0$, cf. section 4). To
establish (\ref{e:2}) we need the following
\begin{lemma} Let $M$ be a nondegenerate CR manifold, $\theta$ a contact
form on $M$, and $d \, {\rm vol}(g_\theta )$ the canonical volume
form associated to the Webster metric $g_\theta$. Then $\theta
\wedge (d \theta )^n = \pm c_n \; d \, {\rm vol}(g_\theta )$ where
$c_n = (-1)^s 2^n n !$, provided that the Levi form $L_\theta$ has
$s$ negative eigenvalues. \label{l:1}
\end{lemma}
This corrects the constant $c_n$ from \cite{kn:Ura1}, p. 546. If
the Levi form $L_\theta$ has $r$ positive and $s$ negative
eigenvalues ($r+s=n$) then $g_\theta$ is a semi-Riemannian metric
of signature $(2r+1, 2s)$. Let $\mathcal{O}$ be a fixed
orientation of $M$. To prove Lemma \ref{l:1}, let $G_{AB}$ be the
components of the Webster metric with respect to a chart $(U , x^A
) \in \mathcal{O}$, so that $d \, {\rm vol}(g_\theta ) =
\sqrt{|\det (G_{AB})|} d x^1 \wedge \cdots \wedge d x^{2n+1}$. Let
$\{ T_\alpha \}$ be a local frame of $T_{1,0}(M)$ and $\mu \in
GL(2n+1, \mathbb{C})$ such that $T_A = \mu_A^B
\partial /\partial x^B$. Then $d {\rm vol}(g_\theta ) =
\sqrt{|\det (G_{AB})|} \det (\mu ) \theta^{01\cdots n \overline{1}
\cdots \overline{n}}$, where $\theta^{01\cdots n \overline{1}
\cdots \overline{n}}$ is short for $\theta \wedge \theta^1 \wedge
\cdots \wedge \theta^n \wedge \theta^{\overline{1}} \wedge \cdots
\wedge \theta^{\overline{n}}$, hence $\overline{\det (\mu )} =
(-1)^{n^2} \det (\mu )$ (as $d {\rm vol}(g_\theta )$ is a real
form). It follows that $\sqrt{|\det (G_{AB})|} = (-1)^s |\det (\mu
)|^{-1} \det (g_{\alpha\overline{\beta}})$. A calculation shows
that
\[ \theta \wedge (d \theta )^n = 2^n i^{n^2} n! \det
(g_{\alpha\overline{\beta}}) \, \theta^{01 \cdots n \overline{1}
\cdots \overline{n}} \] (cf. also \cite{kn:Lee2}) and then $\theta
\wedge (d \theta )^n = \pm c_n d\, {\rm vol}(g_\theta )$. The sign
is $+1$ if $\mathcal{O}$ and the orientation of $H(M)$ (induced by
its complex structure $J$) agree. Lemma \ref{l:1} is proved. Let
us prove (\ref{e:2}). As
\[ R^D e_j = (\pi_H R^D ) e_j - 2 {{P_{0A}}^i}_j (\theta \wedge \theta^A )
\otimes e_i \] it follows that
\[ \| R^D \|^2 = \| \pi_H R^D \|^2 + 4 g^{A\overline{B}}
{{P_{0A}}^i}_j {{P_{0B}}^j}_i \, . \] Yet $i_T R^D = \theta^B
\otimes [{{P_{0B}}^i}_j] \in \Omega^1 ({\rm Ad} (E))$ hence $\|
R^D \|^2 = \| \pi_H R^D \|^2 + 4 \| i_T R_D \|^2$. At this point
we may integrate over $M$ with respect to $\theta \wedge (d \theta
)^n$ and use Lemma \ref{l:1}. To proof of the last statement in i)
of Theorem \ref{t:1} is delegated to the next section.

\section{The first variation formula} Let $E \to M$ be a vector
bundle and $D$ a connection in $E$. We shall need the differential
operator $d^D : \Omega^k (E) \to \Omega^{k+1}(E)$ given by
\[ (d^D \varphi )(X_1 , \cdots , X_{k+1}) = \sum_{i=1}^{k+1}
(-1)^{i+1} D_{X_i} (\varphi (X_1 , \cdots , \hat{X}_i , \cdots ,
X_{k+1})) + \]
\[ + \sum_{1 \leq i < j \leq k+1} (-1)^{i+j} \varphi ([X_i , X_j ]
, X_1 , \cdots , \hat{X}_i , \cdots , \hat{X}_j ,  \cdots ,
X_{k+1}) \] for any $\varphi \in \Omega^k (E)$ and any $X_i \in
T(M)$, $1 \leq i \leq k$. Here a hat indicates, as usual, the
suppression of a term. Let $D \in \mathcal{C}(E,h)$ and let us
denote by the same symbol the connection induced by $D$ in ${\rm
Ad}(E) \to M$. The operator $\delta^D$ in (\ref{e:3}) is the
formal adjoint of $d^D : \Omega^1 ({\rm Ad}(E)) \to \Omega^2 ({\rm
Ad}(E))$ with respect to the inner product
\begin{equation}
(\varphi , \psi ) = \int_M \langle \varphi , \psi \rangle \;
\theta \wedge (d \theta )^n \, ,  \;\; \varphi , \psi \in \Omega^k
(E). \label{e:ip}
\end{equation}
Let $\varphi \in \Omega^1 ({\rm Ad}(E))$. A standard calculation
shows that $R^{D+t\varphi} = R^D + t \, d^D \varphi + t^2 \,
[\varphi \wedge \varphi ]$ (where $[\varphi \wedge \psi ]_{X,Y} =
[\varphi_X , \psi_Y ] - [\varphi_Y , \psi_X ]$, $X,Y \in T(M)$,
$\varphi , \psi \in \Omega^1 ({\rm Ad}(E))$) hence
\[ \| \pi_H \, R^{D+t \varphi} \|^2 = \| \pi_H \, R^D \|^2 + 2 t \;
\langle \pi_H \, R^D \, , \, \pi_H \, d^D  \varphi \rangle + O(t^2
) \] and
\[ \frac{d}{dt} \{ \mathcal{PYM}(D + t \varphi )\}_{t=0} =
\frac{1}{2} \int_M \frac{d}{dt} \{ \| \pi_H R^{D + t \varphi} \|^2
\}_{t=0} \; \theta \wedge (d \theta )^n = \] \[ = \int_M \langle
\pi_H \, R^D \, , \, d^D \varphi \rangle \; \theta \wedge (d
\theta )^n = \int_M \langle \delta^D \pi_H R^D \, , \, \varphi
\rangle \; \theta \wedge (d \theta )^n . \] Then $\frac{d}{dt} \{
\mathcal{PYM}(D + t \varphi )\}_{t=0} = 0$ yields
\begin{equation} \delta^D \, \pi_H \, R^D = 0.
\label{e:16}
\end{equation}
Let $D \in \mathcal{C}(E , h)$ such that $i_T R^D = 0$. Then (by
(\ref{e:3}) and (\ref{e:16})) $D$ is a pseudo Yang-Mills field if
and only if $D$ is a Yang-Mills field, and the last statement in
part i) of Theorem \ref{t:1} follows from Theorem 2.3 in
\cite{kn:Ura1}, p. 551. \par Let us consider the operator
$\delta^D_b : \Omega^{k+1} (E) \to \Omega^k(E)$ given by
\[ (\delta^D_b \varphi )(X_1 , \cdots , X_k ) = - \sum_{a=1}^{2n}
(D_{E_a} \varphi )(E_a , X_1 , \cdots , X_k ), \] for any $\varphi
\in \Omega^{k+1} (E)$ and $X_i \in T(M)$, $1 \leq i \leq k$, where
$\{ E_a : 1 \leq a \leq 2n \}$ is a local $G_\theta$-orthonormal
frame of $H(M)$. Clearly, if $\varphi \in \Omega^k_H (E)$ then
$\delta^D \varphi = \delta^D_b \varphi$ and $i_T \delta^D_b
\varphi = 0$. Consequently, if $i_T R^D = 0$ then the equations
(\ref{e:16}) may also be written
\begin{equation}
\delta^D_b R^D = 0.
\label{e:pymb}
\end{equation}
Now we attack the problem whether the pullback $\hat{D} = \pi^* D$
of a pseudo Yang-Mills field $D$ on $M$ is a Yang-Mills field on
$C(M)$. As argued in the previous section, this doesn't follow
directly from (\ref{e:1}). In turn, the Yang-Mills equations on
$C(M)$ are related to (\ref{e:16}) due to
\begin{equation} (\delta^{\hat{D}} R^{\hat{D}} )(X^\uparrow )
\hat{u} = \left( (\delta^D_b R^D )(X) u + R^D (T , J X) u
\right)\hat{\,} \; , \label{e:18}
\end{equation}
\begin{equation} (\delta^{\hat{D}} R^{\hat{D}} )(T^\uparrow ) \hat{u} =
\left( (\delta^D_b R^D )(T) u \right)\hat{\,} \;\; - \label{e:19}
\end{equation}
\[ - \frac{i}{n+2} ( R^{\alpha\overline{\beta}} -
\frac{\rho}{2(n+1)} \, g^{\alpha\overline{\beta}} ) \left( R^D
(T_\alpha , T_{\overline{\beta}} ) u \right)\hat{\,} \; , \]
\begin{equation}
(\delta^{\hat{D}} R^{\hat{D}} )(S) \hat{u} = 2 \left(
(\Lambda_\theta R^D ) u \right)\hat{\,} \; ,\label{e:20}
\end{equation}
for any $X \in H(M)$ and $u \in \Omega^0 (E)$. Here $X^\uparrow$
is the horizontal lift of $X$ with respect to the connection
$1$-form $\sigma$ in $S^1 \to C(M) \to M$. Let $D \in
\mathcal{C}(E,h)$ be a pseudo Yang-Mills field with $i_T R^D = 0$.
Then (by (\ref{e:18})-(\ref{e:20})) $\delta^{\hat{D}} R^{\hat{D}}
= 0$ if and only if (\ref{e:01})-(\ref{e:02}) hold. This completes
the proof of Theorem \ref{t:1}.
\par
It remains that we prove (\ref{e:18})-(\ref{e:20}). The formal
adjoint $\delta^{\mathbb{D}}$ of $d^{\mathbb{D}} : \Omega^1 ({\rm
Ad}(\pi^* E)) \to \Omega^2 ({\rm Ad}(\pi^* E))$ is given by
\[ (\delta^{\mathbb{D}} \psi )(Y) v = -
\sum_{j=1}^{2n+2} \epsilon_j (\mathbb{D}_{X_j} \psi )(X_j , Y) v =
\] \[ = - \sum_{j=1}^{2n+2} \left( \mathbb{D}_{X_j} \psi (X_j , Y) v -
\psi (\nabla^{C(M)}_{X_j} X_j , Y) v - \right. \] \[ \left. - \psi
(X_j , \nabla^{C(M)}_{X_j} Y )v - \psi (X_j , Y)\mathbb{D}_{X_j} v
\right) ,
\] for any $\psi \in \Omega^2 ({\rm Ad} (\pi^* E))$, $Y
\in T(C(M))$ and $v \in \Omega^0 (\pi^* E)$, where $\{ X_j : 1
\leq j \leq 2n+2 \}$ is a local orthonormal (i.e. $F_\theta (X_j ,
X_k ) = \epsilon_j \delta_{jk}$, $\epsilon_1 = \cdots =
\epsilon_{2n+1} = - \epsilon_{2n+2} = 1$) frame of $T(C(M))$ and
$\nabla^{C(M)}$ is the Levi-Civita connection of $(C(M), F_\theta
)$. As $S^1 \to C(M) \to M$ is a principal bundle, the projection
$\pi$ is a submersion. However, if $S = ((n+2)/2) \partial
/\partial \gamma$ then $F_\theta (S , S) = 0$, i.e. $S$ is null,
so that $\pi$ is not a semi-Riemannian submersion (in the sense of
\cite{kn:O'Nei1}, p. 212, as the fibers of $\pi$ are degenerate
submanifolds). Nevertheless, we may relate $\nabla^{C(M)}$ to the
Tanaka-Webster connection $\nabla$ of $(M , \theta )$, very much
in the spirit of \cite{kn:O'Nei2}. Precisely, we may state
\begin{lemma} For any $X,Y \in H(M)$
\[ \nabla^{C(M)}_{X^\uparrow} Y^\uparrow = (\nabla_X Y)^\uparrow -
(d \theta )(X,Y) T^\uparrow - (A(X,Y) + (d \sigma )(X^\uparrow ,
Y^\uparrow )) S, \]
\[ \nabla^{C(M)}_{X^\uparrow} T^\uparrow  = (\tau X + \phi
X)^\uparrow , \]
\[ \nabla^{C(M)}_{T^\uparrow} X^\uparrow = (\nabla_T X + \phi
X)^\uparrow + 2(d \sigma )(X^\uparrow , T^\uparrow ) S,  \]
\[ \nabla^{C(M)}_{X^\uparrow} S = \nabla^{C(M)}_S X^\uparrow = (J
X)^\uparrow , \]
\[ \nabla^{C(M)}_{T^\uparrow} T^\uparrow = V^\uparrow , \;\;
\nabla^{C(M)}_S S = 0, \]
\[ \nabla^{C(M)}_S T^\uparrow = \nabla^{C(M)}_{T^\uparrow} S = 0,
\]
where $\phi : H(M) \to H(M)$ is given by $G_\theta (\phi X , Y) =
(d \sigma )(X^\uparrow , Y^\uparrow )$, and $V \in H(M)$ is given
by $G_\theta (V , Y) = 2 (d \sigma )(T^\uparrow , Y^\uparrow )$.
\label{l:3}
\end{lemma}
{\em Proof of Lemma $\ref{l:3}$}. Let us recall that \[ 2 F_\theta
(\nabla^{C(M)}_{\tilde{X}} \tilde{Y} , \tilde{Z} ) = \tilde{X}
(F_\theta (\tilde{Y} , \tilde{Z})) + \tilde{Y}(F_\theta (\tilde{X}
, \tilde{Z})) - \tilde{Z}(F_\theta (\tilde{X}, \tilde{Y})) + \]
\begin{equation} + F_\theta ([\tilde{X}, \tilde{Y}], \tilde{Z}) + F_\theta
([\tilde{Z} , \tilde{X}], \tilde{Y}) + F_\theta (\tilde{X} ,
[\tilde{Z} , \tilde{Y}]), \label{e:21}\end{equation} for any
$\tilde{X} , \tilde{Y}, \tilde{Z} \in T(C(M))$. In particular for
$\tilde{X} = X^\uparrow , \, \tilde{Y} = Y^\uparrow , \, \tilde{Z}
= Z^\uparrow$, for any $X,Y,Z \in H(M)$
\[ F_\theta (\nabla^{C(M)}_{X^\uparrow} Y^\uparrow , Z^\uparrow )
= g_\theta (\nabla^M_X Y , Z), \] where $\nabla^M$ is the
Levi-Civita connection of $(M , g_\theta )$. Here one used the
fact that $[X,Y]^\uparrow$ is the horizontal component of
$[X^\uparrow ,Y^\uparrow ]$, with respect to $\sigma$ (cf. e.g.
\cite{kn:KoNo}, Vol. I, p. 65). The Levi-Civita connection
$\nabla^M$ and the Tanaka-Webster connection $\nabla$ of $(M ,
\theta )$ are related by (cf. (1) in \cite{kn:BaDr}, p. 238)
\begin{equation}
\nabla^M = \nabla - (d \theta + A) \otimes T + \tau \otimes \theta
+ 2 \theta \odot J, \label{e:22}
\end{equation}
where $A(X,Y) = g_\theta (X , \tau Y)$. Recall that $A$ is
symmetric and $\tau$ trace-less (cf. \cite{kn:Web}). As $H(M)$ is
$\nabla$-parallel $\pi_H \nabla^M_X Y = \nabla_X Y$, where $\pi_H
: T(M) \to H(M)$ is the projection associated with the direct sum
decomposition $T(M) = H(M) \oplus \mathbb{R} T$. Therefore, by
taking into account the decomposition $T(C(M)) = {\rm Ker}(\sigma
) \oplus {\rm Ker}(d \pi ) = H(M)^\uparrow \oplus \mathbb{R}
T^\uparrow \oplus \mathbb{R} S$
\begin{equation}
\nabla^{C(M)}_{X^\uparrow} Y^\uparrow = (\nabla_X Y)^\uparrow +
\lambda T^\uparrow + \mu S , \label{e:23}
\end{equation}
for some $\lambda , \mu \in C^\infty (C(M))$, depending on $X, Y$.
We may determine $\lambda$, $\mu$ by taking the inner product with
$S$, $T^\uparrow$, respectively. To this end let us first observe
that
\[ F_\theta (\nabla^{C(M)}_{X^\uparrow} Y^\uparrow , S) = - (d
\theta )(X, Y). \] Here we used again (\ref{e:21}) together with
the fact that $[X^\uparrow , S] = 0$ (cf. e.g. \cite{kn:KoNo},
Vol. I, p. 79). Similarly
\[ 2 F_\theta (\nabla^{C(M)}_{X^\uparrow} Y^\uparrow , T^\uparrow
) = \sigma ([X^\uparrow , Y^\uparrow ]) - \] \[ - T(g_\theta
(X,Y)) + g_\theta ([T,X], Y) + g_\theta (X , [T,Y]) \] and
\[ 2 g_\theta (\nabla^M_X Y , T) = - T(g_\theta (X,Y)) + \theta
([X,Y]) + \]
\[ + g_\theta ([T,X], Y) + g_\theta (X , [T,Y]) \]
hence
\[ 2 F_\theta (\nabla^{C(M)}_{X^\uparrow} Y^\uparrow , T^\uparrow
) = 2 \theta (\nabla^M_X Y) - \theta ([X,Y]) + \sigma ([X^\uparrow
, Y^\uparrow ]) \] or (by (\ref{e:22}))
\[ F_\theta (\nabla^{C(M)}_{X^\uparrow} Y^\uparrow , T^\uparrow
) = - A(X,Y) - (d \sigma )( X^\uparrow , Y^\uparrow ). \] Summing
up, (\ref{e:23}) leads to the first identity in Lemma \ref{l:3}.
The proof of the remaining identities in Lemma \ref{l:3} may be
obtained in a similar manner. Let us go back to the proof of
(\ref{e:18})-(\ref{e:20}). Let $\{ E_a : 1 \leq a \leq 2n \}$ be a
local orthonormal frame of the Levi distribution $H(M)$. Then $\{
E_a^\uparrow , T^\uparrow \pm S \}$ is a local orthonormal frame
of $T(C(M))$ with respect to the Feferman metric $F_\theta$. We
make use of $i_S (\hat{D} \hat{u}) = 0$ and $i_S (R^{\hat{D}}
\hat{u}) = 0$, for any $u \in \Omega^0 (E)$. Then (by Lemma
\ref{l:3})
\[ (\delta^{\hat{D}} R^{\hat{D}} )(X^\uparrow ) \hat{u} = -
\sum_{a=1}^{2n} (\hat{D}_{E_a^\uparrow} R^{\hat{D}} )(E_a^\uparrow
, X^\uparrow ) - \]
\[ - (\hat{D}_{T^\uparrow + S} R^{\hat{D}})(T^\uparrow + S ,
X^\uparrow ) + (\hat{D}_{T^\uparrow - S} R^{\hat{D}}) (T^\uparrow
- S , X^\uparrow ) = \]
\[ = \{ (\delta^D_b R^D )(X) u + 2 R^D (T , J X) u \}\hat{\,} -
\sum_{a=1}^{2n} \{ (d \theta )(E_a , X) \, R^D (E_a , T) u
\}\hat{\,} \] and
\[ \sum_{a=1}^{2n} (d \theta )(E_a , X) E_a = - J X \]
hence (\ref{e:18}) is proved. Similarly
\begin{equation}
(\delta^{\hat{D}} R^{\hat{D}})(T^\uparrow ) \hat{u} = ((\delta^D_b
R^D )(T) u )\hat{\,} + \sum_{a=1}^{2n} \{ R^D (E_a , \tau E_a +
\phi E_a )u \}\hat{\,} . \label{e:0}
\end{equation}
Now, on one hand
\[ \sum_{a=1}^{2n} R^D (E_a , \tau E_a ) u =
g^{\alpha\overline{\beta}} \{ R^D (T_\alpha , \tau
T_{\overline{\beta}} ) u + R^D (T_{\overline{\beta}} , \tau
T_\alpha ) u \} = \]
\[ = A^{\alpha\gamma} R^D (T_\alpha , T_\gamma ) u +
A^{\overline{\beta} \overline{\gamma}} R^D (T_{\overline{\beta}} ,
T_{\overline{\gamma}}) u = 0 \] (as $A_{\alpha\beta} =
A_{\beta\alpha}$) with the corresponding simplification of
(\ref{e:0}). On the other hand
\[ \sum_{a=1}^{2n} R^D (E_a , \phi E_a ) u
= \phi^{\alpha\gamma} R^D (T_\alpha , T_\gamma )u + \phi^{\alpha
\overline{\gamma}} R^{D}(T_\alpha , T_{\overline{\gamma}}) u + \]
\[ + \phi^{\overline{\beta}\gamma} R^D (T_{\overline{\beta}} ,
T_\gamma ) u + \phi^{\overline{\beta}\overline{\gamma}} R^D
(T_{\overline{\beta}} , T_{\overline{\gamma}}) u , \] where $\phi
T_\alpha = {\phi_\alpha}^\beta T_\beta +
{\phi_\alpha}^{\overline{\beta}} T_{\overline{\beta}}$,
$\phi^{\alpha\beta} = g^{\alpha\overline{\gamma}}
{\phi_{\overline{\gamma}}}^{\beta}$, etc. Let us take the exterior
derivative of (\ref{e:5}) so that to obtain
\[ (n+2) d \sigma = \pi^* (i \, d {\omega_\alpha}^\alpha -
\frac{i}{2} \, d g^{\alpha\overline{\beta}} \wedge d
g_{\alpha\overline{\beta}} - \frac{1}{4(n+1)} \, d (\rho \theta )
). \] Using the identities $d g_{\alpha\overline{\beta}} =
g_{\alpha\overline{\gamma}}
{\omega_{\overline{\beta}}}^{\overline{\gamma}} +
{\omega_\alpha}^\gamma g_{\gamma\overline{\beta}}$ (a consequence
of $\nabla g_\theta = 0$) and $d g^{\alpha\overline{\beta}} = -
g^{\gamma\overline{\beta}} g^{\alpha\overline{\rho}} d
g_{\overline{\rho}\gamma}$ (a consequence of
$g^{\alpha\overline{\beta}} g_{\overline{\beta}\gamma} =
\delta^\alpha_\gamma$) it follows that
\[  d g^{\alpha\overline{\beta}} \wedge d
g_{\alpha\overline{\beta}} = \omega_{\alpha\overline{\beta}}
\wedge \omega^{\alpha\overline{\beta}} +
\omega_{\overline{\alpha}\beta} \wedge
\omega^{\overline{\alpha}\beta} = 0. \] Also (cf. e.g.
\cite{kn:Web})
\[ d {\omega_\alpha}^\alpha = R_{\lambda\overline{\mu}} \,
\theta^\alpha \wedge \theta^{\overline{\mu}} +
(W^\alpha_{\alpha\lambda} \theta^\lambda -
W^\alpha_{\alpha\overline{\mu}} \theta^{\overline{\lambda}} )
\wedge \theta \] where $R_{\lambda\overline{\mu}}$ is the
pseudohermitian Ricci curvature and $W^\alpha_{\alpha\lambda}$
(respectively $W^\alpha_{\alpha\overline{\mu}}$) are certain
contractions of the covariant derivatives of
$A^\alpha_{\overline{\beta}}$. It follows that
\[ (n+2) G_\theta (\phi X , Y) = i (R_{\alpha\overline{\beta}} \,
\theta^\alpha \wedge \theta^{\overline{\beta}})(X,Y) -
\frac{\rho}{4(n+1)} \, (d \theta )(X,Y), \] for any $X,Y \in
H(M)$. Therefore
\[ \phi^{\overline{\alpha}\beta} = \frac{i}{2(n+2)}
(R^{\overline{\alpha}\beta} - \frac{\rho}{2(n+1)} \,
g^{\overline{\alpha}\beta} ), \;\; \; \phi^{\alpha\beta} = 0. \]
We may conclude that
\[ \sum_{A=1}^{2n} R^D (E_a , \phi E_a ) u = - \frac{i}{n+2}
(R^{\alpha\overline{\beta}} - \frac{\rho}{2(n+1)} \,
g^{\alpha\overline{\beta}} ) \, R^D (T_\alpha ,
T_{\overline{\beta}}) u \] and (\ref{e:0}) leads to (\ref{e:19}).
Finally (again by Lemma \ref{l:3})
\[ (\delta^{\hat{D}} R^{\hat{D}} )(S) \hat{u} = \sum_{a=1}^{2n}
(R^D (E_a , J E_a ) u )\hat{\,} = \]
\[ = - 2i \{ g^{\alpha\overline{\beta}} R^D (T_\alpha ,
T_{\overline{\beta}}) u \}\hat{\,} = 2 \{ (\Lambda_\theta R^D )u
\}\hat{\,} \] and (\ref{e:20}) is proved.

\section{The second variation formula}
Let $\{ D^t : |t| < \epsilon \}$ be a smooth family of connections
in $E$, where $D = D^0$ is a pseudo Yang-Mills field. We write
$D^t = D + A^t$, where $A^t \in \Omega^1 ({\rm Ad}(E))$ for each
$|t| < \epsilon$. The curvature $R^t$ of $D^t$ is then given by
\[ R^t = R^D + d^D A^t + \frac{1}{2} \, [A^t \wedge A^t ] \]
(cf. e.g. (6.2) in \cite{kn:BoLa}, p. 212). Next, let us set
$\varphi = \{ d A^t /d t\}_{t = 0}$ and $\psi = \{ d^2 A^t /d
t^2\}_{t=0}$ and observe that
\[ \| \pi_H R^t \|^2 = \| \pi_H R^D \|^2 + 2 t \langle \pi_H R^D ,
d^D \varphi \rangle + \]
\[ + t^2 \{ 2 \langle \pi_H R^D , d^D \psi \rangle + \langle \pi_H
R^D , [\varphi \wedge \varphi ] \rangle + \| \pi_H d^D \varphi
\|^2 \} + O(t^3 ). \] Integrating by parts and using $\delta^D
\pi_H R^D = 0$ we obtain
\begin{equation}
\frac{d^2}{d t^2} \{ \mathcal{PYM}(D^t )\}_{t=0} = \label{e:26}
\end{equation}
\[ = \int_M \{ \|
\pi_H d^D \varphi \|^2 + \langle \pi_H R^D , [\varphi \wedge
\varphi ] \rangle \} \, \theta \wedge (d \theta )^n . \] We shall
need the (zero order) operator $\mathcal{R}^D : \Omega^1 ({\rm
Ad}(E)) \to \Omega^1 ({\rm Ad}(E))$ given by
\[ \mathcal{R}^D (\varphi )_X = \sum_{A=1}^{2n+1} [R^D_{E_A , X}
\, , \, \varphi_{E_A} ], \] for any $X \in T(M)$, $\varphi \in
\Omega^1 ({\rm Ad}(E))$, where $\{ E_A \}$ is a local orthonormal
frame of $(T(M), g_\theta )$. Then (cf. (6.7) in \cite{kn:BoLa},
p. 213)
\begin{equation}
\langle [\varphi \wedge \varphi ] \, , \, R^D \rangle = \langle
\varphi \, , \, \mathcal{R}^D (\varphi ) \rangle . \label{e:27}
\end{equation}
Let us set
\[ \mathcal{R}^D_b (\varphi )_X = \sum_{a=1}^{2n} [R^D_{E_a , X} \, ,
\, \varphi_{E_a}] , \;\; \mathcal{R}^D_0 (\varphi )_X = [R^D_{T,X}
\, , \, \varphi_T ] , \] where $\{ E_a \}$ is a local orthonormal
frame of $(H(M), G_\theta )$, so that $\mathcal{R}^D =
\mathcal{R}^D_b + \mathcal{R}^D_0$. Taking into account the
identities
\[ [\varphi \wedge \varphi ] e_j = 4 (\varphi^i_k \wedge
\varphi^k_j ) \otimes e_i \, , \]
\[ (\theta \wedge i_T R^D ) e_j = ({{P_{0A}}^i}_j \, \theta \wedge
\theta^A ) \otimes e_i \, , \]
\[ \mathcal{R}^D_0 (\varphi ) e_j = \{ \varphi^k_j (T)
{{P_{0A}}^i}_k - {{P_{0A}}^k}_j \, \varphi^i_k (T) \} \, \theta^A
\otimes e_i \, , \] where $\varphi \, e_j = \varphi^i_j \otimes
e_i$, $\varphi^i_j \in \Omega^1 (U)$, we may conduct the
calculations
\[ \langle [\varphi \wedge \varphi ] \, , \, \theta \wedge i_T R^D
\rangle  = h^{j\overline{r}} \langle [\varphi \wedge \varphi ] e_j
\, , \, (\theta \wedge i_T R^D ) e_r \rangle = \]
\[ = 4 h^{j\overline{r}} h_{i\overline{s}} \; g^*_\theta
(\varphi^i_k \wedge \varphi^k_j \, , \, {{P_{0A}}^s}_r \, \theta
\wedge \theta^A ) = \]
\[ = 2 h^{j\overline{r}} h_{i\overline{s}} g^{A\overline{B}} \{
\varphi^i_k (T) \varphi^k_j (T_A ) - \varphi^i_k (T_A )
\varphi^k_j (T) \}
{{P_{0\overline{B}}}^{\overline{s}}}_{\overline{r}}  \; , \] and
\[ \langle \varphi \, , \, \mathcal{R}^D_0 (\varphi )\rangle =
h^{j\overline{r}} h_{i\overline{s}} \varphi^i_j (T_A ) \{
\varphi^{\overline{k}}_{\overline{r}} (T)
{{P_{0\overline{B}}}^{\overline{s}}}_{\overline{k}} -
{{P_{0\overline{B}}}^{\overline{k}}}_{\overline{r}} \,
\varphi^{\overline{s}}_{\overline{k}} (T) \} g^{A\overline{B}} =
\]
\[ = \varphi^i_j (T_A ) g^{A\overline{B}}
\{ h^{j\overline{k}} h_{i\overline{s}}
\varphi^{\overline{r}}_{\overline{k}} (T) - h^{j\overline{r}}
h_{i\overline{k}} \varphi^{\overline{k}}_{\overline{s}} (T) \} \;
{{P_{0\overline{B}}}^{\overline{s}}}_{\overline{r}} \; . \] Assume
now that $\{ e_j \}$ is orthonormal $(h_{i\overline{j}} =
\delta_{ij})$, so that $\varphi^{\overline{i}}_{\overline{j}} = -
\varphi^j_i$ (as $\varphi$ is ${\rm Ad}(E)$-valued). Then
\[ \langle [\varphi \wedge \varphi ] \, , \, \theta \wedge i_T R^D
\rangle = 2 \sum_{i,j} \{ \varphi^i_k (T) \varphi^k_j (T_A ) -
\varphi^i_k (T_A ) \varphi^k_j (T) \} g^{A\overline{B}} \,
{{P_{0\overline{B}}}^{\overline{i}}}_{\overline{j}} \]
\[ \langle \varphi \, , \, \mathcal{R}^D_0 (\varphi )\rangle =
\sum_{r,s} \{ \varphi^s_i (T) \varphi^i_r (T_A ) - \varphi^s_j
(T_A ) \varphi^j_r (T) \} g^{A\overline{B}}
{{P_{0\overline{B}}}^{\overline{s}}}_{\overline{r}} \] hence
\[ \langle [\varphi \wedge \varphi ] \, , \, \theta \wedge i_T R^D
\rangle = 2 \langle \varphi \, , \, \mathcal{R}^D_0 (\varphi
)\rangle . \] Finally, let us take into account (\ref{e:27}) and
the identity
\[ R^D = \pi_H R^D + 2 \, \theta \wedge i_T R^D . \]
We obtain
\[ \langle [\varphi \wedge \varphi ] \, , \, \pi_H R^D \rangle =
\langle \varphi \, , \, \mathcal{R}^D (\varphi ) - 4 \,
\mathcal{R}^D_0 (\varphi ) \rangle , \] so that (\ref{e:26})
becomes
\begin{equation}
\frac{d^2}{d t^2} \{ \mathcal{PYM}(D^t )\}_{t=0} = \label{e:28}
\end{equation}
\[ = \int_M \langle \delta^D \pi_H d^D \varphi + \mathcal{R}^D
(\varphi ) - 4 \mathcal{R}^D_0 (\varphi ) \, , \, \varphi \rangle
\, \theta \wedge (d \theta )^n . \] We now restrict our variations
to those whose first order part $\varphi$ satisfies $i_T \varphi =
0$ and $\delta^D \varphi = 0$. Also, let us introduce the first
order differential operator $d_b^D : \Omega^1 ({\rm Ad}(E)) \to
\Omega^2_H ({\rm Ad}(E))$ given by $d^D_b \equiv \pi_H \circ d^D$.
Then $\delta^D \pi_H d^D \varphi = \delta^D_b d^D_b \varphi =
\Delta^D_b \varphi$ and $\mathcal{R}^D_0 (\varphi ) = 0$, so that
(\ref{e:28}) yields (\ref{e:t2}) in Theorem \ref{t:2}. Here
$\Delta^D_b \equiv d^D_b \delta^D_b + \delta^D_b d^D_b$ is the
{\em generalized sublaplacian}. The Riemannian counterpart
$\mathcal{S}^D = \Delta^D + \mathcal{R}^D$ (cf. \cite{kn:BoLa}, p.
213, where $\Delta^D$ is the generalized Hodge-de Rham laplacian)
of $\mathcal{S}^D_b =\Delta^D_b + \mathcal{R}^D_b$ in Theorem
\ref{t:2} is an elliptic operator, hence its restriction to ${\rm
Ker}(\delta^D ) \subset \Omega^1 ({\rm Ad}(E))$ has a discrete
spectrum tending to $+\infty$ and the eigenspace corresponding to
each eigenvalue of $\mathcal{S}^D$ is finite dimensional. This
allows one to employ concepts from Morse theory (cf. Definition
6.10 in \cite{kn:BoLa}, p. 213) in order to discuss stability and
weak stability of Yang-Mills fields (cf. \cite{kn:BoLa}, p. 214).
The CR analog of this phenomenon is that $\Delta_b^D :
\Omega^{0,1}({\rm Ad}(E)) \to \Omega^{0,1}({\rm Ad}(E))$ is
subelliptic of order $1/2$, where $\Omega^{0,q}({\rm Ad}(E)) =
\Gamma^\infty (\Lambda^{0,q}(M) \otimes {\rm Ad}(E))$. A complex
valued $q$-form $\eta$ on $M$ is of {\em type} $(0,q)$, or a
$(0,q)$-{\em form}, if $T_{1,0}(M) \, \rfloor \, \eta = 0$ and
$i_T \eta = 0$. We denote by $\Lambda^{0,q}(M) \to M$ the relevant
bundle and set $\Omega^{0,q}(M) = \Gamma^{\infty}
(\Lambda^{0,q}(M))$. Let $M$ be a strictly pseudoconvex CR
manifold (not necessarily compact). It is the proper place to
recall that a formally self adjoint second order differential
operator $L : C^\infty (M) \to C^\infty (M)$ is {\em subelliptic}
(of order $0 < \epsilon \leq 1$) at a point $x \in M$ if there is
a neighborhood $U$ of $x$ such that
\[ \| u \|_\epsilon^2 \leq C \left( (L u , u) + \| u \|^2 \right) \]
for any $u \in C^\infty_0 (U)$, where $\| u \|_\epsilon$ is the
Sobolev norm $u$ of order $\epsilon$, $\| u \| = (u,u)^{1/2}$, and
\begin{equation}
(u,v) = \int_M u v \; d {\rm vol}(g_\theta ) \label{e:31}
\end{equation}
is the ordinary $L^2$ inner product. $L$ is subelliptic (of order
$\epsilon$) if it is subelliptic at any $x \in M$. A typical
example is the {\em sublaplacian}
\[ \Delta_b u = - {\rm div} (\nabla^H u ), \;\;\; u \in C^\infty
(M), \] where $\nabla^H u \equiv \pi_H \nabla u$, $\nabla u$ is
the gradient of $u$ with respect to the Webster metric $g_\theta$,
and the divergence is defined with respect to the volume form
$\omega \equiv \theta \wedge (d \theta )^n$, i.e.
\[ \mathcal{L}_X \omega = {\rm div}(X) \; \omega , \]
for any $X \in \mathcal{X}(M)$, where $\mathcal{L}_X$ is the Lie
derivative. It is easily seen that $\Delta_b u = - \sum_{a
=1}^{2n} E_a^* E_a u$, for any local orthonormal frame $\{ E_a \}$
of $H(M)$ hence, by a well known lemma of E.V. Radkevic,
\cite{kn:Rad}, it follows that $\| u \|^2_{1/2} \leq C \left(
(\Delta_b u , u) + \| u \|^2 \right)$, for any $u \in C^\infty_0
(U)$, i.e. $\Delta_b$ is subelliptic of order $1/2$. Here $E^*_a$
is the formal adjoint of $E_a$ with respect to the inner product
(\ref{e:31}). In the next section we relate $\Delta^D_b$ to the
Kohn-Rossi operator $\square_b$ and explain the subellipticity of
$\square_b$ on $(0,1)$-forms.

\section{Subellipticity of $\Delta^D_b$}
Let $\{ T_\alpha \}$ be a local frame of $T_{1,0}(M)$. We start by
computing
\[ (\delta^D_b d^D_b \varphi ) T_\lambda = - \sum_{a=1}^{2n}
(D_{E_a} d^D_b \varphi )(E_a , T_\lambda ), \] for any $\varphi
\in \Omega^1_H ({\rm Ad}(E))$. Let us take into account the
identities
\[ g^{\alpha\overline{\beta}} (D_{T_\alpha} d^D_b \varphi
)(T_{\overline{\beta}} , T_\lambda ) = \] \[ =
g^{\alpha\overline{\beta}} \{ D_{T_\alpha} ((d^D_b \varphi
)(T_{\overline{\beta}} , T_\lambda )) -
\Gamma^{\overline{\gamma}}_{\alpha\overline{\beta}} (d^D_b \varphi
)(T_{\overline{\gamma}} , T_\lambda ) -
\Gamma^\gamma_{\alpha\lambda} (d^D_b \varphi
)(T_{\overline{\beta}} , T_\gamma ) \} , \]
\[ (d^D_b \varphi )(T_\alpha , T_{\overline{\beta}}) e_j = \{
(\nabla_{T_\alpha} \varphi^i_j ) T_{\overline{\beta}} -
(\nabla_{T_{\overline{\beta}}} \varphi^i_j ) T_\alpha + \]
\[ + 2 (\omega^i_k \wedge \varphi^k_j + \varphi^i_k \wedge
\omega^k_j )(T_\alpha , T_{\overline{\beta}} ) \} e_i \, , \]
where $\omega^i_j$ are the connection $1$-forms of $D$ with
respect to $\{ e_i \}$. Also, let us set
\[ \nabla_{\overline{\mu}} \nabla_\lambda \psi_{\overline{\alpha}}
\equiv (\nabla_{T_{\overline{\mu}}} \nabla \psi )(T_\lambda ,
T_{\overline{\alpha}}) = T_{\overline{\mu}} (\nabla_\lambda
\psi_{\overline{\alpha}} ) - \Gamma^\rho_{\overline{\mu}\lambda}
\nabla_\rho \psi_{\overline{\alpha}} -
\Gamma^{\overline{\rho}}_{\overline{\mu}\, \overline{\alpha}}
\psi_{\overline{\rho}} \, , \] for any $\psi =
\psi_{\overline{\alpha}} \theta^{\overline{\alpha}} \in
\Omega^{0,1}(M)$, where $\nabla_\alpha \psi_{\overline{\beta}}
\equiv (\nabla_{T_\alpha} \psi ) T_{\overline{\beta}}$. We obtain
\[ g^{\alpha\overline{\beta}} (D_{T_\alpha} d^D_b \varphi
)(T_{\overline{\beta}} , T_\lambda ) e_j = -
g^{\alpha\overline{\beta}} \{ \nabla_\alpha \nabla_\lambda
\varphi^i_{j\overline{\beta}} - \nabla_\alpha
\nabla_{\overline{\beta}} \varphi^i_{j\lambda} + \]
\[ + \omega^k_j (T_{\overline{\beta}}) A^i_k (T_\alpha ,
T_\lambda ) - \omega^k_j (T_\lambda ) A^i_k (T_\alpha ,
T_{\overline{\beta}} ) +  \] \[ + \omega^i_k (T_\lambda ) B^k_j
(T_\alpha , T_{\overline{\beta}}) - \omega^i_k
(T_{\overline{\beta}}) B^k_j (T_\alpha , T_\lambda ) +
\]
\[ + \omega^i_k (T_\alpha )
[(\nabla_{T_\lambda} \varphi^k_j )T_{\overline{\beta}} -
(\nabla_{T_{\overline{\beta}}} \varphi^k_j ) T_\lambda ] + \]
\[ + \omega^k_j (T_\alpha ) [(\nabla_{T_{\overline{\beta}}}
\varphi^i_k ) T_\lambda - (\nabla_{T_\lambda} \varphi^i_k )
T_{\overline{\beta}} ] +
\]
\[ + \varphi^k_{j\overline{\beta}} C^i_k (T_\alpha , T_\lambda
) - \varphi^k_{j\lambda} C^i_k (T_\alpha , T_{\overline{\beta}}) +
\] \[ + \varphi^i_{k \lambda} D^k_j (T_\alpha , T_{\overline{\beta}})
- \varphi^i_{k\overline{\beta}} D^k_j (T_\alpha ,
T_{\overline{\beta}}) \} \otimes e_i \] where
\[ A^i_j \equiv \nabla \varphi^i_j + \omega^i_k \otimes
\varphi^k_j \; , \;\;\; B^i_j \equiv \nabla \varphi^i_j -
\omega^k_j \otimes \varphi^i_k \; , \]
\[ C^i_j \equiv \nabla \omega^i_j  + \omega^i_k \otimes
\omega^k_j \; , \;\;\; D^i_j \equiv \nabla \omega^i_j - \omega^k_j
\otimes \omega^i_k \; , \] and $\varphi^i_{jA} = \varphi^i_j (T_A
)$, so that $\nabla_\alpha \nabla_{\overline{\beta}}
\varphi^i_{j\lambda}$ (respectively $\nabla_\alpha \nabla_\beta
\varphi^i_{j\overline{\lambda}}$) is the second order covariant
derivative of the $(1,0)$-form $\pi_{1,0} \varphi^i_j$
(respectively of the $(0,1)$-form $\pi_{0,1} \varphi^i_j$)
($\pi_{1,0} : \Omega^1 (M) \to \Omega^{1,0}(M)$ and $\pi_{0,1} :
\Omega^1 (M) \to \Omega^{0,1}(M)$ are the natural projections).
The previous identity is rather involved, yet one is interested in
the second order terms alone. Together with the similar expression
for $g^{\alpha\overline{\beta}} (D_{T_{\overline{\beta}}} d^D_b
\varphi )(T_\alpha , T_\lambda ) e_j$ this leads to
\begin{equation} (\delta^D_b d^D_b \varphi )(T_\lambda ) e_j =
g^{\alpha\overline{\beta}} (\nabla_\alpha \nabla_\lambda
\varphi^i_{j\overline{\beta}} + \nabla_{\overline{\beta}}
\nabla_\lambda \varphi^i_{j\alpha} - \label{e:32}
\end{equation}
\[ - \nabla_\alpha \nabla_{\overline{\beta}}  \varphi^i_{j\lambda} -
\nabla_{\overline{\beta}} \nabla_\alpha \varphi^i_{j\lambda}) \,
e_i + lower \; order \; terms. \] By {\em lower order terms}
(l.o.t.) we mean a linear combination of $\nabla_A \varphi^i_{jB}$
and $\varphi^i_{jB}$ (with $C^\infty (U)$-coefficients). Next, we
need to compute
\[ (d^D_b \delta^D_b \varphi ) T_\lambda = D_{T_\lambda}
(\delta^D_b \varphi ) = - \sum_{a=1}^{2n} D_{T_\lambda} ((D_{E_a}
\varphi ) E_a ). \] We have
\[ [D_{T_\lambda} (g^{\alpha\overline{\beta}} (D_{T_\alpha}
\varphi ) T_{\overline{\beta}})] e_j = g^{\alpha\overline{\beta}}
\{ \nabla_\lambda \nabla_\alpha \varphi^i_{j\overline{\beta}} +
\]
\[ - \omega^k_j (T_\lambda ) A^i_k (T_\alpha , T_{\overline{\beta}})
+ \omega^i_k (T_\lambda ) B^k_j (T_\alpha , T_{\overline{\beta}})
- \]
\[ - \omega^j_k (T_\alpha ) (\nabla_{T_\lambda} \varphi^i_k ) T_{\overline{\beta}} +
\omega^i_k (T_\alpha ) (\nabla_{T_\lambda} \varphi^k_j
)T_{\overline{\beta}} + \]
\[ + \varphi^k_{j\overline{\beta}} C^i_k (T_\lambda , T_\alpha ) -
\varphi^i_{k\overline{\beta}} D^k_j (T_\lambda , T_\alpha ) \} e_i
\, . \] Together with the similar expression for $D_{T_\lambda}
(g^{\alpha\overline{\beta}} (D_{T_{\overline{\beta}}} \varphi )
T_\alpha )$ this yields
\begin{equation}
(d^D_b \delta^D_b \varphi )(T_\lambda ) e_j = -
g^{\alpha\overline{\beta}} \{ \nabla_\lambda \nabla_\alpha
\varphi^i_{j\overline{\beta}} + \nabla_\lambda
\nabla_{\overline{\beta}} \varphi^i_{j\alpha} \} e_i + \; l.o.t.
\label{e:33}
\end{equation}
We shall need the commutation formulae
\[ \nabla_\alpha \nabla_\beta \eta_{\overline{\gamma}} -
\nabla_\beta \nabla_\alpha \eta_{\overline{\gamma}} = -
\eta_{\overline{\rho}} \,
{{R_{\overline{\gamma}}}^{\overline{\rho}}}_{\alpha\beta} \, , \]
\[ \nabla_{\overline{\beta}} \nabla_\alpha \eta_\gamma -
\nabla_\alpha \nabla_{\overline{\beta}} \eta_\gamma = 2 i
g_{\alpha\overline{\beta}} \nabla_0 \eta_\gamma - \eta_\rho \,
{{R_\gamma}^\rho}_{\overline{\beta}\alpha} \, , \] where the
convention for the curvature components (of the Tanaka-Webster
connection) is $R(T_A , T_B ) T_C = {{R_C}^D}_{AB} T_D$. Then (by
(\ref{e:32})-(\ref{e:33}))
\[ (\Delta^D_b \varphi )(T_\lambda ) e_j =
g^{\alpha\overline{\beta}} ( \nabla_\alpha \nabla_\lambda
\varphi^i_{j\overline{\beta}} - \nabla_\lambda \nabla_\alpha
\varphi^i_{j\overline{\beta}} + \] \[ + \nabla_{\overline{\beta}}
\nabla_\lambda \varphi^i_{j\alpha} - \nabla_\lambda
\nabla_{\overline{\beta}} \varphi^i_{j\alpha} - \nabla_\alpha
\nabla_{\overline{\beta}} \varphi^i_{j\lambda} -
\nabla_{\overline{\beta}} \nabla_\alpha \varphi^i_{j\lambda} ) e_i
+ \; l.o.t. = \]
\[ = \{ - \varphi^i_j (T_{\overline{\rho}})
g^{\alpha\overline{\beta}}
{{R_{\overline{\beta}}}^{\overline{\rho}}}_{\alpha\lambda} + 2
\sqrt{-1} \nabla_0 \varphi^i_{j\lambda} - \varphi^i_j (T_\rho )
g^{\alpha\overline{\beta}}
{{R_\alpha}^\rho}_{\overline{\beta}\lambda} - \] \[ - 2
g^{\alpha\overline{\beta}} \nabla_{\overline{\beta}} \nabla_\alpha
\varphi^i_{j\lambda} + 2 \sqrt{-1} n \nabla_0 \varphi^i_{j\lambda}
- \varphi^i_j (T_\rho ) g^{\alpha\overline{\beta}}
{{R_\lambda}^\rho}_{\overline{\beta}\alpha} \} e_i + \; l.o.t. \]
At this point we need the Kohn-Rossi operator $\square_b$ on
$\Omega^{0,1}(M)$. We start by extending $\overline{\partial}_b$
(originally defined on functions, cf. section 3) to $(0,1)$-forms.
Precisely, if $\eta \in \Omega^{0,1}(M)$ then
$\overline{\partial}_b \eta$ is the unique $(0,2)$-form on $M$
coinciding with $d \eta$ on $T_{0,1}(M) \otimes T_{0,1}(M)$. Next,
let us set $\square_b \equiv \overline{\partial}^{\, *}_b
\overline{\partial}_b + \overline{\partial}_b
\overline{\partial}^{\, *}_b$, where $\overline{\partial}^{\,
*}_b$ is the formal adjoint of $\overline{\partial}_b$ with
respect to the $L^2$ inner product $(\alpha , \beta ) = \int_M
g^*_\theta (\alpha , \overline{\beta}) \, \omega$. A
straightforward calculation leads to
\begin{equation}
\square_b \eta = ( - g^{\alpha\overline{\beta}} \, \nabla_\alpha
\nabla_{\overline{\beta}} \eta_{\overline{\gamma}} - 2 i \nabla_0
\eta_{\overline{\gamma}} + \eta_{\overline{\rho}} \,
{R^{\overline{\rho}}}_{\overline{\gamma}} ) \,
\theta^{\overline{\gamma}} , \label{e:35}
\end{equation}
for any $\eta = \eta_{\overline{\gamma}}\,
\theta^{\overline{\gamma}} \in \Omega^{1,0}(U)$. Consequently
\[ \overline{\square}_b (\pi_{1,0} \varphi^i_j ) = \left( \square_b (\pi_{1,0}
\varphi^i_j )^- \right)^- = \] \[ = ( - g^{\alpha\overline{\beta}}
\nabla_{\overline{\beta}} \nabla_\alpha \varphi^i_{j\lambda} + 2 i
\nabla_0 \varphi^i_{j\lambda} + \varphi^i_j (T_\rho )
{R^\rho}_\lambda ) \theta^\lambda \] hence
\[ (\Delta^D_b \varphi )(T_\lambda ) \theta^\lambda \otimes e_j =
2 \overline{\square}_b (\pi_{1,0} \varphi^i_j ) e_i + \{ 2(n-1)
\sqrt{-1} \, \nabla_0 \varphi^i_{j\lambda} - \]
\[ - \varphi^i_j (T_{\overline{\rho}}) \,
g^{\alpha\overline{\beta}}
{{R_{\overline{\beta}}}^{\overline{\rho}}}_{\alpha\lambda} -
\varphi^i_j (T_\rho ) [2 {R^\rho}_\lambda +
g^{\alpha\overline{\beta}}(
{{R_\alpha}^\rho}_{\overline{\beta}\lambda} +
{{R_\lambda}^\rho}_{\overline{\beta}\alpha}) ] \} \,
\theta^\lambda \otimes e_i + \; l.o.t. \] To compute the curvature
terms we need the identities
\[ R_{\lambda\overline{\mu}} =
{{R_\lambda}^\alpha}_{\alpha\overline{\mu}} \, , \]
\[ {{R_\alpha}^\rho}_{\lambda\mu} = 2i (A_{\mu\alpha}
\delta^\rho_\lambda - A_{\lambda\alpha} \delta^\rho_\mu ) \, , \]
\[ {{R_\alpha}^\rho}_{\overline{\lambda} \, \overline{\mu}} = 2 i
(g_{\alpha\overline{\lambda}} A^\rho_{\overline{\mu}} -
g_{\alpha\overline{\mu}} A^\rho_{\overline{\lambda}} ) \, , \]
\[ R_{\alpha\overline{\beta}\lambda\overline{\mu}} =
R_{\lambda\overline{\beta}\alpha\overline{\mu}} \, , \;\;
R_{\alpha\overline{\beta}\lambda\overline{\mu}} = -
R_{\overline{\beta}\alpha\lambda\overline{\mu}} \, , \] following
essentially by the techniques developed in \cite{kn:Web}. Indeed
we have
\[ g^{\alpha\overline{\beta}}
{{R_{\overline{\beta}}}^{\overline{\rho}}}_{\alpha\lambda} = -2 i
(n-1) A^{\overline{\rho}}_\lambda \, , \]
\[ 2 {R^\rho}_\lambda +
g^{\alpha\overline{\beta}}(
{{R_\alpha}^\rho}_{\overline{\beta}\lambda} +
{{R_\lambda}^\rho}_{\overline{\beta}\alpha}) = 0, \] hence
\begin{equation} (\Delta^D_b \varphi )(T_\lambda ) \theta^\lambda \otimes
e_j = 2 \overline{\square}_b (\pi_{1,0} \varphi^i_j ) \otimes e_i +
\label{e:34}
\end{equation}
\[ + 2(n-1) \sqrt{-1} \, [(\nabla_T \varphi^i_j ) T_\lambda +
A_\lambda^{\overline{\rho}} \, \varphi^i_j (T_{\overline{\rho}})]
\theta^\lambda \otimes e_i + \; l.o.t. \] Recall that $\tau
T_\alpha = A_\alpha^{\overline{\beta}} T_{\overline{\beta}}$. Then
(\ref{e:34}) together with the similar identity $(\Delta^D_b
\varphi )(T_{\overline{\lambda}}) e_j = 2 \{ (\square_b \pi_{0,1}
\varphi^i_j )_{\overline{\lambda}} - (n-1) \sqrt{-1}\, (\nabla_0
\varphi^i_{j\overline{\lambda}} + A^\rho_{\overline{\lambda}}
\varphi^i_{j\rho} ) \} e_i + \, l.o.t.$ leads to
\[ (\Delta^D_b \varphi ) \otimes e_j = 2 \{ \square_b (\pi_{0,1}
\varphi^i_j ) + \overline{\square}_b (\pi_{1,0} \varphi^i_j ) +
\] \[ + (n-1) (\nabla_T \varphi^i_j + \varphi^i_j \circ \tau ) \circ J \}
\otimes e_i \; + l.o.t.  \] and therefore to (\ref{e:t2.2}) when
$\varphi$ is a $(0,1)$-form. Finally, let us show that $\square_b$
is subelliptic on (scalar) $(0,1)$-forms. As $\nabla \omega = 0$,
the sublaplacian may be computed as
\[ \Delta_b f = - trace \{ T_A \mapsto \nabla_{T_A} \nabla^H f \}
= - \nabla_\alpha f^\alpha - \nabla_{\overline{\alpha}}
f^{\overline{\alpha}} \, , \] for any $C^\infty$ function $f : M
\to \mathbb{C}$, where $f^\alpha = g^{\alpha\overline{\beta}}
T_{\overline{\beta}}(f)$. As $\nabla g_\theta = 0$
\[ g^{\alpha\overline{\beta}} \nabla_\alpha
\nabla_{\overline{\beta}} \, \eta_{\overline{\lambda}} =
\nabla_\alpha (g^{\alpha\overline{\beta}}
\nabla_{\overline{\beta}} \, \eta_{\overline{\lambda}} ) =
\nabla_\alpha (\eta_{\overline{\lambda}})^\alpha - \nabla_\alpha
(g^{\alpha\overline{\beta}}
\Gamma^{\overline{\rho}}_{\overline{\beta} \, \overline{\lambda}}
\eta_{\overline{\rho}} ). \] A similar expression holds for
$g^{\alpha\overline{\beta}} (\nabla_{\overline{\beta}}
\nabla_\alpha \eta_{\overline{\lambda}})$. Adding up the two
identities leads to
\[ g^{\alpha\overline{\beta}} (\nabla_\alpha
\nabla_{\overline{\beta}} \, \eta_{\overline{\lambda}} +
\nabla_{\overline{\beta}} \nabla_\alpha \eta_{\overline{\lambda}})
= - \Delta_b \eta_{\overline{\lambda}} + \; l.o.t. \] or (by the
commutation formulae for the second order derivatives and
(\ref{e:35}))
\[  2 (\square_b \eta )_{\overline{\lambda}} - 2 i (n-1) \nabla_0
\eta_{\overline{\lambda}} - \eta_{\overline{\rho}}
{R^{\overline{\rho}}}_{\overline{\lambda}} = \] \[ = \Delta_b
\eta_{\overline{\lambda}} + g^{\alpha\overline{\beta}} \{
\nabla_\alpha (\Gamma^{\overline{\rho}}_{\overline{\beta} \,
\overline{\lambda}} \eta_{\overline{\rho}} ) +
\nabla_{\overline{\beta}} (\Gamma^{\overline{\rho}}_{\alpha
\overline{\lambda}} \eta_{\overline{\rho}} )\} . \] Hence $\square_b$
is subelliptic on $\Omega^{0,1}(M)$, i.e. $\square_b \eta$ is
locally given by a subelliptic operator acting on the coefficients
of $\eta$, plus lower order terms. In particular (by
(\ref{e:t2.2})) $(\Delta^D_b \varphi )\otimes e_j = \Delta_b
(\varphi^i_j (T_{\overline{\alpha}})) \theta^{\overline{\alpha}}
\otimes e_i +$ l.o.t. Theorem \ref{t:2} is completely proved.

\begin{appendix}
\section{The Graham-Lee connection}
Let $\Omega = \{ \varphi < 0 \} \subset \mathbb{C}^n$ be a
strictly pseudoconvex domain and $\mathcal{F}$ the foliation by
level sets of $\varphi$ of a one-sided neighborhood $V$ of
$\partial \Omega$ (as in section 3 of this paper). Let $\{
W_\alpha : 1 \leq \alpha \leq n-1 \}$ be a local frame of
$T_{1,0}(\mathcal{F})$, so that $\{ W_\alpha , \xi \}$ is a local
frame of $T^{1,0}(V)$. Let $g_\theta$ be the tensor field given by
\begin{equation}
g_\theta (X,Y) = (d \theta )(X, J Y), \;\; g_\theta (X,T) = 0,
\;\; g_\theta (T , T) = 1, \label{e:A.1}
\end{equation}
for any $X,Y \in H(\mathcal{F})$. Then $g_\theta$ is a tangential
Riemannian metric for $\mathcal{F}$, i.e. a Riemannian metric in
$T(\mathcal{F}) \to V$. We consider as well
\[ L_\theta (Z,  \overline{W}) \equiv - i (d \theta )(Z, \overline{W}), \;\;\; Z,W
\in T_{1,0}(\mathcal{F}). \]
Note that $L_\theta$ and (the
$\mathbb{C}$-linear extension of) $g_\theta$ coincide on
$T_{1,0}(\mathcal{F}) \otimes T_{0,1}(\mathcal{F})$. We set
$g_{\alpha\overline{\beta}} = g_\theta (W_\alpha ,
W_{\overline{\beta}})$. Let $\{ \theta^\alpha : 1 \leq \alpha \leq
n-1 \}$ be the (locally defined) complex $1$-forms on $V$
determined by \[ \theta^\alpha (W_\beta ) = \delta^\alpha_\beta \,
, \;\; \theta^\alpha (W_{\overline{\beta}}) = 0, \;\;
\theta^\alpha (T) = 0, \;\; \theta^\alpha (N) = 0. \] Then $\{
\theta^\alpha , \, \theta^{\overline{\alpha}} , \, \theta , \, d
\varphi \}$ is a local frame of $T(V) \otimes \mathbb{C}$ and we
may look for $d \theta$ in the form
\[ d \theta = B_{\alpha\beta} \, \theta^\alpha \wedge \theta^\beta
+ B_{\alpha\overline{\beta}} \, \theta^\alpha \wedge
\theta^{\overline{\beta}} + B_{\overline{\alpha} \,
\overline{\beta}} \, \theta^{\overline{\alpha}} \wedge
\theta^{\overline{\beta}} + \] \[ + (B_\alpha \, \theta^\alpha +
B_{\overline{\alpha}} \, \theta^{\overline{\alpha}}) \wedge \theta
+ (C_\alpha \, \theta^\alpha + C_{\overline{\alpha}}\,
\theta^{\overline{\alpha}}) \wedge d \varphi + D \, d\varphi
\wedge \theta . \]
 As $d \theta = i \partial \overline{\partial}
\varphi \in \Omega^{1,1}(U)$ it follows that $B_{\alpha\beta} =
0$, $B_{\overline{\alpha}\, \overline{\beta}} = 0$. Also
\[ g_{\alpha\overline{\beta}} = g_\theta (W_\alpha  ,
W_{\overline{\beta}}) = - i (d \theta )(W_\alpha ,
W_{\overline{\beta}}) = - \frac{i}{2} \,
B_{\alpha\overline{\beta}} \] i.e. $B_{\alpha \overline{\beta}} =
2 i g_{\alpha\overline{\beta}}$. Next
\[ \frac{1}{2} B_\alpha = (d \theta )(W_\alpha , T) = i \partial
\overline{\partial} \varphi (W_\alpha , T) = 0 \] as $T = i(\xi -
\overline{\xi})$ (and $\xi$ is orthogonal to
$T_{1,0}(\mathcal{F})$ with respect to $\partial
\overline{\partial} \varphi$), i.e. $B_\alpha = 0$,
$B_{\overline{\alpha}} = 0$. Similarly $C_\alpha = 0$,
$C_{\overline{\alpha}} = 0$. Finally
\[ D = (d \theta )(N, T) = i \, \partial \overline{\partial}
\varphi (N , T) = 2 \, \partial \overline{\partial} \varphi (\xi ,
\overline{\xi}) = r \] i.e. $D = r$. We obtain the identity
\begin{equation}
d \theta = 2 i g_{\alpha\overline{\beta}} \, \theta^\alpha \wedge
\theta^{\overline{\beta}} + r \, d \varphi \wedge \theta .
\label{e:A.2}
\end{equation}
As an immediate consequence
\begin{equation}
i_T \, d \theta = - \frac{r}{2} \, d \varphi , \label{e:A.3}
\end{equation}
\begin{equation}
i_N \, d \theta = r \, \theta . \label{e:A.4}
\end{equation}
For instance (by (\ref{e:A.2}))
\[ (d \theta )(X,T) = \frac{r}{2} \{ (d \varphi )(X) - (d \varphi
)(T) \theta (X) \} , \] for any $X \in T(\mathcal{F})$, hence (as
$(d \varphi )(T) = 0$) one derives (\ref{e:A.3}). As an
application of (\ref{e:A.2}) we decompose $[T,N]$ (according to
$T(V) \otimes \mathbb{C} = T_{1,0}(\mathcal{F}) \oplus
T_{0,1}(\mathcal{F}) \oplus \mathbb{C} T \oplus \mathbb{C} N$).
This is a bit trickier, as shown below. By (\ref{e:A.3})
\[ \theta ([T,N]) = - 2 (d \theta )(T, N) = r \, d \varphi
(N) = 2 r. \] Next
\[ 2(d \theta )(W_\alpha , [T,N]) = 2 W_\alpha (r) - \theta
([W_\alpha , [T,N]]) = \;\; (Jacobi's \;\; identity) \]
\[ = 2 W_\alpha (r) + \theta ([T, [N,W_\alpha ]]) + \theta ([N,
[W_\alpha , T]]) = \]
\[ = 2 W_\alpha (r) + 2(d \theta )(T, [W_\alpha , N]) - T(\theta
([W_\alpha , N])) + \]
\[ + 2 (d \theta )(N , [T, W_\alpha ]) - N(\theta ([T , W_\alpha
])) \] hence (by (\ref{e:A.3})-(\ref{e:A.4}))
\[  (d \theta )(W_\alpha , [T,N]) = W_\alpha (r). \] We conclude
that
\begin{equation}
[T,N] = i \, W^\alpha (r) W_\alpha - i \, W^{\overline{\alpha}}(r)
W_{\overline{\alpha}} + 2 r T, \label{e:A.5}
\end{equation}
where $W^\alpha (r) = g^{\alpha\overline{\beta}}
W_{\overline{\beta}}(r)$ and $W^{\overline{\alpha}} (r) =
\overline{W^\alpha (r)}$.
\par
Let $\nabla$ be a linear connection on $V$. Let us consider the
$T(V)$-valued $1$-form $\tau$ on $V$ defined by
\[ \tau (X) = T_\nabla (T , X), \;\;\; X \in T(V), \]
where $T_\nabla$ is the torsion tensor field of $\nabla$. We say
$T_\nabla$ is {\em pure} if
\begin{equation}
T_\nabla (Z, W) = 0, \;\; T_\nabla (Z, \overline{W}) = 2 i
L_\theta (Z , \overline{W}) T, \label{e:A.6}
\end{equation}
\begin{equation}
T_\nabla (N , W) = r \, W + i \, \tau (W), \label{e:A.7}
\end{equation}
for any $Z,W \in T_{1,0}(\mathcal{F})$, and
\begin{equation}
\tau (T_{1,0}(\mathcal{F})) \subseteq T_{0,1}(\mathcal{F}),
\label{e:A.8}
\end{equation}
\begin{equation}
\tau (N) = -  \, J \, \nabla^H r -  2 r \, T. \label{e:A.9}
\end{equation}
Here $\nabla^H r$ is defined by $\nabla^H r = \pi_H \nabla r$ and
$g_\theta (\nabla r , X) = X(r)$, $X \in T(\mathcal{F})$. Also
$\pi_H : T(\mathcal{F}) \to H(\mathcal{F})$ is the projection
associated to the direct sum decomposition $T(\mathcal{F}) =
H(\mathcal{F}) \oplus \mathbb{R} T$. Appendix A is aimed at the
following
\begin{theorem} $(${\rm C.R. Graham \& J.M. Lee,
\cite{kn:GrLe}}$)$ \par\noindent There is a unique linear
connection $\nabla$ on $V$ such that {\rm i)}
$T_{1,0}(\mathcal{F})$ is parallel with respect to $\nabla$, {\rm
ii)} $\nabla L_\theta = 0$, $\nabla T = 0$, $\nabla N = 0$, and
{\rm iii)} $T_\nabla$ is pure. \label{t:A.1}
\end{theorem} \noindent  $\nabla$ given by
Theorem \ref{t:A.1} is the {\em Graham-Lee connection}. Compare to
Proposition 1.1 in \cite{kn:GrLe}, p. 701-702. The axiomatic
description in Theorem \ref{t:A.1} is however new (cf. also
Theorem 2 in \cite{kn:DrNi}). We first establish
\begin{lemma} Let $\phi : T(\mathcal{F}) \to T(\mathcal{F})$ be
the bundle morphism given by $\phi (X) = J X$, for any $X \in
H(\mathcal{F})$, and $\phi (T) = 0$. Then
\[ \phi^2 = - I + \theta \otimes T, \]
\[ g_\theta (X , T) = \theta (X), \]
\[ g_\theta (\phi X , \phi Y) = g_\theta (X,Y) - \theta (X) \theta (Y), \]
for any $X,Y \in T(\mathcal{F})$. Moreover, if $\nabla$ is a
linear connection on $V$ satisfying the axioms {\rm (i)-(iii)} in
Theorem $\ref{t:A.1}$ then
\begin{equation}
\phi \circ \tau + \tau \circ \phi = 0 \label{e:A.10}
\end{equation}
along $T(\mathcal{F})$. Consequently $\tau$ may be computed as
\begin{equation}
\tau (X) = - \frac{1}{2} \phi (\mathcal{L}_T \phi ) X,
\label{e:A.11}
\end{equation} for any $X \in H(\mathcal{F})$.
\label{l:A.1}
\end{lemma}
{\em Proof}. For any $X \in T(\mathcal{F})$
\[ \phi (X) = \phi (\pi_H X + \theta (X) T) = J (\pi_H X) \in H(\mathcal{F}), \]
\[ \phi^2 (X) = J^2 (\pi_H X) = - \pi_H X = - X + \theta (X) T. \]
The second statement in Lemma \ref{l:A.1} follows from definitions
(cf. (\ref{e:A.1})). The third identity follows from
\[ g_\theta (\phi X , \phi Y) = (d \theta )(\phi \, \pi_H \, X ,
\phi^2\,  \pi_H\,  Y) = \] \[ = g_\theta (\pi_H \, Y , \pi_H \, X)
= g_\theta (Y,X) - \theta (X) g_\theta (Y, T). \] Let us prove
(\ref{e:A.10}). As $\tau (T_{1,0}(\mathcal{F})) \subseteq
T_{0,1}(\mathcal{F})$ (cf. axiom (\ref{e:A.8})) there are complex
valued functions $A_\alpha^{\overline{\beta}}$ such that $\tau
(W_\alpha ) = A_\alpha^{\overline{\beta}} W_{\overline{\beta}}$.
Then \[ (\tau \circ \phi + \phi \circ \tau ) W_\alpha = i \tau
(W_\alpha ) + A_\alpha^{\overline{\beta}} \phi
(W_{\overline{\beta}}) = 0. \] It remains that we check
(\ref{e:A.11}). As $T_{1,0}(\mathcal{F})$ is parallel with respect
to $\nabla$ and $\nabla$ is a real operator it follows that
$T_{0,1}(\mathcal{F})$ is parallel, hence both $H(\mathcal{F})$
and its complex structure $\left. J_{\mathcal{F}} \equiv
J\right|_{H(\mathcal{F})}$ are parallel. Moreover, as $\nabla T =
0$, it follows that $\phi$ is parallel, as well. Let $X \in
H(\mathcal{F})$. Then (by (\ref{e:A.10}))
\[ \phi \tau X = - T_\nabla (T , \phi X) = - \nabla_T \phi X + [T,
\phi X]. \] Applying $\phi$ in both sides gives (as $\nabla \phi =
0$)
\[ \tau X = - \nabla_T X - \phi [T , \phi X] = - [T,X] - \tau X -
\phi [T , \phi X] \] or
\[ 2 \tau X = - \mathcal{L}_T X - \phi \mathcal{L}_T \phi X. \]
Q.e.d.
\par
{\em Proof of Theorem $\ref{t:A.1}$}. To establish uniqueness,
note first that, for any $X = X^{1,0} + X^{0,1} + \theta (X) T \in
T(\mathcal{F})$ (with $X^{1,0} \in T_{1,0}(\mathcal{F})$, $X^{0,1}
= \overline{X^{1,0}}$) one has (by $\nabla N = 0$)
\[ \nabla_N X = [N,X] + T_\nabla (N , X) = \;\; (by \; (\ref{e:A.6})-(\ref{e:A.7}),
\; (\ref{e:A.9}))\]
\[ = [N,X] + r X^{1,0} + i \tau X^{1,0} + r X^{0,1} - i \tau
X^{0,1} + \theta (X) \{ J \nabla^H r + 2 r T \} \] that is
\begin{equation}
\nabla_N X = r X + \tau \phi X - [X,N] +  \theta (X) \{ J \nabla^H
r + r T \} , \label{e:A.12}
\end{equation}
for any $X \in T(\mathcal{F})$. In view of (\ref{e:A.11})
$\nabla_N X$ is determined. As $\nabla N = 0$, $\nabla T = 0$ it
remains that we compute $\nabla_X Z$, for $X \in T(\mathcal{F})$
and $Z \in T_{1,0}(\mathcal{F})$. Note that $\nabla T = 0$,
$\nabla L_\theta = 0$ and $\nabla J_{\mathcal{F}} = 0$ yield
$\nabla g_\theta = 0$, i.e.
\[ X(g_\theta (Y,Z)) = g_\theta (\nabla_X Y , Z ) + g_\theta (Y ,
\nabla_X Z), \] for any $X,Y,Z \in T(\mathcal{F})$. The well known
Christoffel process then leads to
\begin{equation}
2 g_\theta (\nabla_X Y , Z) = X(g_\theta (Y,Z)) + Y(g_\theta (X,
Z)) - Z(g_\theta (X,Y)) + \label{e:A.13}
\end{equation}
\[ + g_\theta ([X,Y], Z) + g_\theta (T_\nabla (X,Y) , Z) + \]
\[ + g_\theta ([Z,X], Y) + g_\theta (T_\nabla (Z,X), Y) + \]
\[ + g_\theta (X , [Z,Y]) + g_\theta (X , T_\nabla (Z,Y)). \]
Note that (again by the purity axioms)
\begin{equation}
T_\nabla (X,Y) = 2 (d \theta )(X,Y) T  + 2 (\theta \wedge \tau
)(X,Y), \label{e:A.14}
\end{equation}
for any $X,Y \in T(\mathcal{F})$. Indeed (by (\ref{e:A.6}))
\[ T_\nabla (X,Y) = - 2 g_\theta (X , \phi Y) T + 2(\theta \wedge
\tau )(X,Y). \] Moreover
\[ g_\theta (X, \phi \, Y) = g_\theta (\pi_H \, X , \phi \, \pi_H\, Y) +
\theta (X) g_\theta (T , \phi\,  \pi_H \, Y) =  \] \[ = - (d
\theta )(\pi_H X , \pi_H Y) = - (d \theta )(X,Y) + \] \[ + \theta
(X) (d \theta )(T, Y) + \theta (Y) (d \theta )(X,T). \] Finally
(by (\ref{e:A.3})) $(d \theta )(X,T) = 0$, $X \in T(\mathcal{F})$,
and (\ref{e:A.14}) is proved. Replacing the torsion terms (from
(\ref{e:A.14}) into (\ref{e:A.13})) leads to
\begin{equation}
2 g_\theta (\nabla_X Z , \overline{W}) = X(g_\theta (Z ,
\overline{W})) + Z(g_\theta (X , \overline{W})) -
\overline{W}(g_\theta (X, Z)) + \label{e:A.15}
\end{equation}
\[ + g_\theta ([X, Z], \overline{W}) + g_\theta ([\overline{W} ,X],Z)
+ g_\theta (X , [\overline{W}, Z]), \] for any $X \in
T(\mathcal{F})$ and $Z,W \in T_{1,0}(\mathcal{F})$, as (by
(\ref{e:A.8}))
\[ \theta (X) \{ g_\theta (\tau Z , \overline{W}) - g_\theta (\tau
\overline{W} , Z) \} = 0. \] The uniqueness statement in Theorem
\ref{t:A.1} is proved. The following explicit expressions of (the
various components of) $\nabla$ are also available. By
(\ref{e:A.6})
\begin{equation}
\nabla_Z \overline{W} = \pi_{0,1} [Z, \overline{W}], \;\;\; Z,W
\in T_{1,0}(\mathcal{F}), \label{e:A.16}
\end{equation}
where $\pi_{0,1} : T(\mathcal{F}) \otimes \mathbb{C} \to
T_{0,1}(\mathcal{F})$ is the projection. Of course
\[ \nabla_{\overline{Z}} W = \overline{\nabla_Z \overline{W}}. \]
Moreover (by $\nabla L_\theta = 0$ and (\ref{e:A.16}))
\[
L_\theta (\nabla_Z W , \overline{V}) = Z(L_\theta (W ,
\overline{V})) - L_\theta (W , \pi_{0,1}[Z , \overline{V}]), \]
for any $Z,W,V \in T_{1,0}(\mathcal{F})$, i.e.
\begin{equation}
\nabla_Z W = g^{\alpha\overline{\beta}} \{ Z(L_\theta (W ,
W_{\overline{\beta}})) - L_\theta (W , \pi_{0,1} [Z ,
W_{\overline{\beta}}])\} W_\alpha \label{e:A.17}
\end{equation}
and
\[ \nabla_{\overline{Z}} \overline{W} = \overline{\nabla_Z W}. \]
Next (by (\ref{e:A.7})-(\ref{e:A.8}))
\begin{equation}
\nabla_N Z = r Z + \pi_{1,0} [N,Z], \label{e:A.18}
\end{equation}
for any $Z \in T_{1,0}(\mathcal{F})$, and
\[ \nabla_N \overline{Z} = \overline{\nabla_N Z}. \]
Finally
\begin{equation}
\nabla_T Z = - \frac{1}{2} \, \phi (\mathcal{L}_T \phi )Z - [Z,T],
 \label{e:A.19}
\end{equation}
\[ \nabla_T \overline{Z} = \overline{\nabla_T Z}, \;\;
Z \in T_{1,0}(\mathcal{F}). \] To establish the existence
statement in Theorem \ref{t:A.1} let $\nabla$ be the linear
connection on $V$ defined by (\ref{e:A.16})-(\ref{e:A.19}) and
$\nabla T = 0$, $\nabla N = 0$. Let us check (i)-(iii) in Theorem
\ref{t:A.1}. Clearly \[ \nabla_Z W , \;\; \nabla_{\overline{Z}} W
, \;\; \nabla_N W \in T_{1,0}(\mathcal{F}), \] for any $W \in
T_{1,0}(\mathcal{F})$, by the very definitions (cf.
(\ref{e:A.16})-(\ref{e:A.18})). Moreover (by (\ref{e:A.19}))
\[ \nabla_T Z = \frac{1}{2} \left( \mathcal{L}_T Z - i \phi
\mathcal{L}_T Z \right) \] as (by (\ref{e:A.3})) $\mathcal{L}_T Z
\in H(\mathcal{F}) \otimes \mathbb{C}$. Therefore
\[ \phi \nabla_T Z = i \, \nabla_T Z \]
that is [as $T_{1,0}(\mathcal{F})$ is the eigenspace corresponding
to the eigenvalue $i$ of (the $\mathbb{C}$-linear extension to
$H(\mathcal{F})$ of) $\phi$] $\nabla_T Z \in
T_{1,0}(\mathcal{F})$. We conclude that $\nabla$ obeys to (i). Let
us check purity. By (\ref{e:A.17})
\[ L_\theta (\nabla_Z W - \nabla_W Z - [Z,W], \overline{V}) = \]
\[ = Z(L_\theta (W , \overline{V})) - L_\theta (W , \pi_{0,1} [Z,
\overline{V}]) - \]
\[ - W(L_\theta (Z, \overline{V})) + L_\theta
(Z , \pi_{0,1} [W , \overline{V}]) - L_\theta ([Z,W],
\overline{V}) = \]
\[ = 3 (d^2 \theta )(Z, W, \overline{V}) = 0. \]
Therefore $T_\nabla (Z,W) = 0$. Next (by (\ref{e:A.16}) and
$\overline{\pi_{0,1} X} = \pi_{1,0} \overline{X}$, $X \in
T(\mathcal{F}) \otimes \mathbb{C}$)
\[ T_\nabla (Z , \overline{W}) = \pi_{0,1} [Z, \overline{W}] -
\overline{\pi_{0,1} [W, \overline{Z}]} - [Z, \overline{W}] = \]
\[ = - \theta ([Z, \overline{W}]) T = 2 i L_\theta (Z ,
\overline{W}) T. \] Moreover (by (\ref{e:A.18}))
\[ T_\nabla (N,Z) = \nabla_N Z - [N,Z] = r Z - \pi_{0,1} [N,Z]. \]
Also (by (\ref{e:A.19}))
\[ \tau (Z) = - \frac{1}{2} \phi (\mathcal{L}_T \phi ) Z, \;\; Z
\in T_{1,0}(\mathcal{F}), \] so that on one hand (\ref{e:A.8}) is
satisfied, and on the other
\[ \tau (Z) = - \frac{1}{2} \{ i \, \phi [T,Z] + [T,Z]\} = - \pi_{0,1} [T,Z] =
\] \[ = i \left\{ \pi_{0,1} [\overline{\xi} , Z] - \pi_{0,1} [\xi
, Z] \right\} = i \, \pi_{0,1}[\overline{\xi}, Z] \] i.e.
\begin{equation}
\tau (Z) = i \, \pi_{0,1} [N, Z], \;\; Z \in T_{1,0}(\mathcal{F}).
\label{e:A.20}
\end{equation}
Here we made use of $T = i(\xi - \overline{\xi})$, $N = \xi +
\overline{\xi}$ and $\pi_{0,1} [\xi , Z] = 0$. Then (\ref{e:A.20})
yields (\ref{e:A.7}). Finally $\nabla T = \nabla N = 0$ and
(\ref{e:A.5}) yield (\ref{e:A.9}) and we conclude that $\nabla$
obeys to (iii). It remains that we check $\nabla L_\theta = 0$.
Clearly $\nabla_Z L_\theta = 0$, $Z \in T_{1,0}(\mathcal{F})$ (by
(\ref{e:A.16})-(\ref{e:A.17})). Next (by (\ref{e:A.19}) and
(\ref{e:A.8}))
\[ (\nabla_T L_\theta )(Z , \overline{W}) = (\mathcal{L}_T
L_\theta )(Z , \overline{W}), \;\; Z, W \in T_{1,0}(\mathcal{F}),
\]
and
\[ (\mathcal{L}_T L_\theta ) (Z, \overline{W}) = - i \{ T(d \theta
(Z , \overline{W})) - d \theta ([T,Z] , \overline{W}) - d \theta
(Z , [T, \overline{W}]) \} = \]
\[ = \frac{i}{2} \{ T(\theta ([Z, \overline{W}])) - \theta
([[T,Z], \overline{W}]) - \theta ([Z, [T, \overline{W}]]) \} = \]
(by applying the Jacobi identity to the term $\theta ([[T,Z],
\overline{W}])$)
\[ = \frac{i}{2} \{ T(\theta ([Z, \overline{W}])) + \theta ([[Z ,
\overline{W}], T]) \} = \] \[ = i (d \theta )(T , [Z,
\overline{W}]) = - \frac{i \; r}{2} \, (d \varphi )([Z,
\overline{W}]) = 0 \] (by (\ref{e:A.3}) and $[Z , \overline{W}]
\in T(\mathcal{F}) \otimes \mathbb{C}$). Hence $\nabla_T L_\theta
= 0$. Finally (by (\ref{e:A.18}))
\[ (\nabla_N L_\theta )(Z, \overline{W}) = - 2 r g_\theta (Z,
\overline{W}) + (\mathcal{L}_N g_\theta )(Z , \overline{W}), \;\;
Z,W \in T_{1,0}(\mathcal{F}), \] and $\nabla_N L_\theta = 0$
follows from (\ref{e:A.23}) in Lemma \ref{l:A.2} below.
\begin{lemma} The following identities hold for any $X \in
T({\mathcal F})$
\begin{equation} T_\nabla (N , X) = r X + \tau (\phi X) +
\theta (X) \{ \phi \nabla^H r + r T \} ,
\label{e:A.21}\end{equation}
\begin{equation}
[N , \phi X] - \phi [N , X] = 2 \tau (X) - \theta (X) \nabla^H r ,
\label{e:A.22}
\end{equation}
Moreover
\begin{equation}
(\mathcal{L}_N g_\theta )(X,Y) = 2 r g_\theta (X,Y) + 2 (d \theta
)(X , \tau (Y)), \label{e:A.23}
\end{equation}
for any $X,Y \in H(\mathcal{F})$. \label{l:A.2}
\end{lemma}
{\em Proof}. (\ref{e:A.21}) follows from (\ref{e:A.7}). Let us
replace $X$ by $\phi X$ in (\ref{e:A.21})
\[ \nabla_N \phi X - [N, \phi X] = r \phi X - \tau (X) \]
and subtract the identity got from (\ref{e:A.21}) by applying
$\phi$ to both sides. Since $\nabla \phi = 0$ we obtain
(\ref{e:A.22}). The proof of (\ref{e:A.23}) is a consequence of
(\ref{e:A.4}), (\ref{e:A.22}), and the Jacobi identity
\[ (\mathcal{L}_N g_\theta )(X,Y)
= N((d \theta )(X , \phi Y)) - (d \theta )([N,X], \phi Y) + \]
\[ + 2 (d \theta )(X , \tau (Y)) - (d \theta )(X , [N , \phi Y]) =
\]
\[ = - \frac{1}{2} \, N(\theta ([X , \phi Y])) + \frac{1}{2} \,
\theta ([[N,X], \phi Y]) + \]
\[ + 2 (d \theta )(X , \tau (Y)) + \frac{1}{2} \; \theta ([X, [N,
\phi Y ]]) = \]
\[ = - \frac{1}{2} \, N(\theta ([X , \phi Y])) - \frac{1}{2} \, \theta ([[X, \phi Y],
N]) - \theta ([X , \tau (Y)]) =
\]
\[ = - (d \theta )(N , [X , \phi Y]) - \theta ([X, \tau (Y)]) = \]
\[ = - r \; \theta ([X, \phi Y]) - \theta ([X , \tau (Y)]) = \] \[
= 2 r (d \theta )(X , \phi Y) + 2 (d \theta )(X , \tau (Y)) = \]
\[ = 2 r g_\theta (X,Y) + 2 (d \theta )(X , \tau (Y)). \] Q.e.d.
Theorem \ref{t:A.1} is proved.
\par
As to the local calculations, if $\varphi^\alpha_\beta$ are the
connection $1$-forms of the Graham-Lee connection (i.e. $\nabla
W_\beta = \varphi^\alpha_\beta \otimes W_\alpha$) then we may look
for $d \theta^\alpha$ in the form
\begin{equation} d \theta^\alpha
= B^\alpha_{\beta\gamma} \, \theta^\beta \wedge \theta^\gamma +
B^\alpha_{\beta\overline{\gamma}} \, \theta^\beta \wedge
\theta^{\overline{\gamma}} + B^\alpha_{\overline{\beta} \,
\overline{\gamma}} \, \theta^{\overline{\beta}} \wedge
\theta^{\overline{\gamma}} + \label{e:A.24}
\end{equation}
\[ + (B^\alpha_\beta \, \theta^\beta + B^\alpha_{\overline{\beta}}
\, \theta^{\overline{\beta}} ) \wedge \theta + (C^\alpha_\beta \,
\theta^\beta + C^\alpha_{\overline{\beta}} \,
\theta^{\overline{\beta}} ) \wedge d \varphi + D \, d \varphi
\wedge \theta . \] Indeed, applying this identity to the pair
$(W_\beta , W_\gamma )$ (respectively to $(W_\beta ,
W_{\overline{\gamma}})$ and $(W_\beta , T)$) gives
\[ B^\alpha_{\beta\gamma} - B^\alpha_{\gamma\beta} =
\varphi^\alpha_\beta (W_\gamma ) - \varphi^\alpha_\gamma (W_\beta
), \;\; B^\alpha_{\beta\overline{\gamma}} = \varphi^\alpha_\beta
(W_{\overline{\gamma}}), \;\; B^\alpha_\beta =
\varphi^\alpha_\beta (T). \] Similarly (applying (\ref{e:A.24}) to
$(W_{\overline{\beta}}, T)$, $(W_\beta , N)$ and
$(W_{\overline{\beta}} , N)$ respectively)
\[ B^\alpha_{\overline{\beta}} = - A^\alpha_{\overline{\beta}} \,
, \;\; C^\alpha_\beta = \frac{1}{2} (\varphi^\alpha_\beta (N) - r
\, \delta^\alpha_\beta ), \;\; C^\alpha_{\overline{\beta}} =
\frac{i}{2} \, A^\alpha_{\overline{\beta}} \, . \] Finally (by
(\ref{e:A.5}))
\[ D \, d \varphi (N) = 2 (d \theta^\alpha )(N , T) = -
\theta^\alpha ([N,T]) = i \, W^\alpha (r). \] Summing up
\begin{equation}
d \theta^\alpha = \theta^\beta \wedge \varphi^\alpha_\beta - i \,
\partial \varphi \wedge \tau^\alpha + \frac{i}{2} \, W^\alpha (r)
\, d \varphi \wedge \theta + \frac{r}{2} \, d \varphi \wedge
\theta^\alpha \, , \label{e:A.25}
\end{equation}
where $\tau^\alpha \equiv A^\alpha_{\overline{\beta}} \,
\theta^{\overline{\beta}}$.
\par Given a linear connection $\nabla$ on $V$ we set $\alpha
(X,Y) \equiv \Pi \nabla_X Y$, for any $X,Y \in T(\mathcal{F})$. If
$\nabla$ is the Graham-Lee connection then (by the proof of
Theorem \ref{t:A.1}) $\alpha = 0$. One may identify, as usual, the
normal bundle $\nu (\mathcal{F}) = T(V)/T(\mathcal{F})$ with
$\mathbb{R} N$. If $\Pi^\bot : T(V) \to T(\mathcal{F})$ is the
projection, let us set $\nabla^{\mathcal{F}} \equiv \Pi^\bot
\nabla$. It is easily seen that $\nabla^{\mathcal{F}}$ is the
Tanaka-Webster connection of each $M_\delta$ (i.e. the pointwise
restriction of the Graham-Lee connection to a leaf of
$\mathcal{F}$ is the Tanaka-Webster connection of the leaf). In
particular $\tau : T(\mathcal{F}) \to T(\mathcal{F})$ is the
pseudohermitian torsion of each leaf (hence $g_\theta (\tau \, X ,
Y) = g_\theta (X, \tau \, Y)$, for any $X,Y \in T(\mathcal{F})$).

\end{appendix}

\end{document}